\documentstyle[12pt]{article}
\newcommand{\sect}[1]{\section{#1}\setcounter{equation}{0}}

\font\mbn=msbm10 scaled \magstep1
\font\mbs=msbm7 scaled \magstep1
\font\mbss=msbm5 scaled \magstep1
\newfam\mbff
\textfont\mbff=\mbn
\scriptfont\mbff=\mbs
\scriptscriptfont\mbff=\mbss\def\mbf{\fam\mbff}
\def\Re{{\mbf R}}
\def\Q{{\mbf Q}}
\def\Z{{\mbf Z}}
\def\Co{{\mbf C}}

\def\P{{\mbf P}}
\def\Di{{\mbf D}}

\def\F{{\mbf F}}
\def\N{{\mbf N}}
\newtheorem{Th}{Theorem}[section]
\newtheorem{Lm}[Th]{Lemma}
\newtheorem{C}[Th]{Corollary}
\newtheorem{D}[Th]{Definition}
\newtheorem{Proposition}[Th]{Proposition}
\newtheorem{R}[Th]{Remark}
\newtheorem{E}[Th]{Example}
\author{Alexander Brudnyi\thanks{Research supported in part by NSERC.
}\\
Department of Mathematics and Statistics\\
University of Calgary, Calgary\\
Canada}
\title{On The Center Problem For Ordinary Differential Equations}
\date{} 
\begin{document} 
\maketitle
\begin{abstract}
{The classical Center-Focus problem posed by H. Poincar\'{e} in 1880's
asks about the characterization of planar polynomial vector fields such
that all their integral trajectories are closed curves whose interiors
contain a fixed point, a {\em center}.
In this paper we describe a new general approach to the Center Problem.}
\end{abstract}
\begin{description}
\item[]
\item[] \hspace*{-1em} {\Large \bf Contents}
\item[{\bf 1.}]\ \ \hspace*{0.0mm} {\bf Introduction.}\hspace*{112mm} {\bf 2}
\item[{\bf 2.}]\ \ \hspace*{0.0mm}
{\bf Monodromy in the Center Problem.} \hspace*{64.7mm} {\bf 6}
\item[{\bf 3.}]\ \ \hspace*{0.0mm}
{\bf Algebraic Model for the Center Problem.} 
\hspace*{52.5mm} {\bf 13}
\item[{\bf 4.}]\ \ \hspace*{0.0mm} {\bf Proofs of Results of Section 2.1.}
\hspace*{71.9mm} {\bf 17}
\item[{\bf 5.}]\ \ \hspace*{0.0mm}
{\bf Proofs of Results of Sections 2.2 and 2.3.}
\hspace*{54.5mm} {\bf 28}
\item[{\bf 6.}]\ \ \hspace*{0.0mm} {\bf Proof of Theorem 1.9.}
\hspace*{91.5mm} {\bf 28}
\item[{\bf 7.}]\ \ \hspace*{0.0mm}
{\bf Proofs of Theorem 2.12 and Corollaries 2.10, 2.15 
and 2.17.} \hspace*{17mm} {\bf 32}
\item[{\bf 8.}]\ \ \hspace*{0.0mm}
{\bf Proofs of Theorems 2.18, 2.20 and Corollaries 2.21 and
2.22.} \hspace*{15.0mm} {\bf 37}
\item[{\bf 9.}]\ \ \hspace*{0.0mm} {\bf Proofs of Results of Section 3.1.}
\hspace*{72mm} {\bf 38}
\item[{\bf 10.}] \hspace*{-0.04mm} {\bf Proofs of Results of Section 3.2.}
\hspace*{71.7mm} {\bf 45}
\item[{\bf 11.}] \hspace*{-0.04mm} {\bf  Proofs of Results of Section 3.3.} 
\hspace*{71.7mm} {\bf 50}
\end{description}
{\sect{\hspace*{-1em}. Introduction.}
{\bf 1.1.} Let $L^{\infty}(S^{1})$ be the Banach space of bounded 
complex-valued functions on the unit circle $S^{1}$.
We will also identify elements of $L^{\infty}(S^{1})$ 
with bounded $2\pi$-periodic functions on $\Re$.
Let us consider the ordinary differential equation
\begin{equation}\label{e2}
\frac{dv}{dx}=\sum_{i=1}^{\infty}a_{i}(x)v^{i+1},
\end{equation}
where all $a_{i}\in L^{\infty}(S^{1})$. If the coefficients of (\ref{e2}) 
grow not very fast for a sufficiently small initial value
one can solve this equation by the Picard iteration method 
so that the (generalized) solution is a Lipschitz function on $[0,2\pi]$.
(I. e., in this way we obtain a function $v$ for which
(\ref{e2}) holds almost everywhere on $[0,2\pi]$.) Moreover, there
is a unique solution with the prescribed initial value $v(0)$. We say that
equation (\ref{e2}) determines a {\em center} if for any sufficiently
small initial values $v(0)$ the solution of (\ref{e2}) satisfies 
$v(0)=v(2\pi)$. 
The Center Problem for equation (\ref{e2}) is to find 
conditions on the coefficients $a_{i}$ under which this equation determines
a center. This problem arises naturally in the framework of the 
qualitative theory of ordinary differential equations created
by H. Poincar\'{e}.
The main purpose of the theory can be described as follows: without
explicitly solving a given differential equation, using only certain 
properties of its right-hand side, to try and give an as complete as possible
description of the geometry of the solution curves of this equation (where
they are defined). In this paper we present a new general
approach to the Center Problem for equation (\ref{e2}). \\
{\bf 1.2.} The Center Problem for (\ref{e2}) is 
closely related to the classical Center-Focus problem first studied by
Poincar\'{e} [P] and further developed by Lyapunov [L], Bendixson [B] and 
Frommer [F]. In this case one considers the system of ODEs in the plane
\begin{equation}\label{e1}
\frac{dx}{dt}=X(x,y),\ \ \ \ \frac{dy}{dt}=Y(x,y)\
\end{equation}
where $X$, $Y$ are real polynomials of degree $d$, or more
generally, real analytic functions defined in an open neighbourhood of
the origin $P:=(0,0)\in\Re^{2}$. One assumes that $P$ is an {\em equilibrium
point}, i.e. $X(0,0)=Y(0,0)=0$. 
The equilibrium point is called a {\em center} if 
there exists an open neighbourhood $U$ of $P$ that does
not contain another equilibrium points such that any trajectory of the 
vector field (\ref{e1}) that intersects $U\setminus\{P\}$ in some point
is closed. Suppose 
$$
X(x,y)=ax+by+\sum_{k=2}^{\infty}X_{k}(x,y),\ \ \ \ \
Y(x,y)=cx+dy+\sum_{k=2}^{\infty}Y_{k}(x,y)
$$
where $X_{k}$, $Y_{k}$ are real homogeneous polynomials of degree $k$. In the
aforementioned papers the case of a non-degenerated equilibrium point was
studied, i.e. the matrix
$$
A=\left(
\begin{array}{cc}
a&b\\
c&d
\end{array}
\right)
$$
assumed to be invertible. It was proved by Poincar\'{e} that a necessary 
condition for $P$ to be a center is that $A$ has pure imaginary eigenvalues. 
In this case, making a linear change of variables in (\ref{e1}) and then a 
linear reparametrization of trajectories one reduces (\ref{e1}) to an 
equivalent system:
\begin{equation}\label{e2'}
\frac{dx}{dt}=-y+F(x,y),\ \ \ \ \ \frac{dy}{dt}=x+G(x,y)
\end{equation}
where $F,G$ are real analytic functions in an open neighbourhood of $P$
whose Taylor expansions at $P$ do not contain
constant and linear terms. For $F,G$ polynomials of a given degree, 
the classical
Poincar\'{e} Center-Focus Problem\footnote{This name originates from the 
fact that in the case of system (\ref{e2'}) the point $P$ can be either a 
center or a focus.} asks about conditions on the coefficients of
$F$ and $G$ under which all trajectories of (\ref{e2'}) 
situated in a small open
neighbourhood of the origin are closed. (A similar problem can be posed for
the general case.) Poincar\'{e} proved that $P$ is a
center if and only if the coefficients of $F$ and $G$ satisfy a certain 
infinite
system of algebraic equations $E_{1}, E_{2},\dots $
such that the coefficients of $E_{n}$ are functions in the solutions of
$E_{1},\dots, E_{n-1}$. Thus in order to solve $E_{n}$ one first should solve
all the previous equations. In the present paper we give
another completely new characterization of centers in terms of an
infinite system of algebraic equations 
$\widetilde E_{1},\widetilde E_{2},\dots$ such that each $\widetilde E_{n}$ 
can be solved {\em independently} of the others (see Theorem \ref{center}). 
This allows to establish
non-existence of a center using only one successfully chosen equation
$\widetilde E_{n}$. Unfortunately, trying to apply either this or Poincar\'{e}
criterion in the converse direction, gives rise, in the general case,
to almost insurmountable difficulties. 
Therefore, as it was emphasized by
Poincar\'{e}, it is the matter of great importance to find typical 
situations for which all equations of the system determining centers are 
satisfied. As an example,  Poincar\'{e} described a case related to 
certain symmetries for $F$ and $G$. 
Some results of our paper can be viewed as a further development in this
direction (see Sections 2.4 and 2.5).

Another necessary and
sufficient condition for $P$ to be a center in (\ref{e2'}) was obtained by 
Lyapunov [L]. His result says that $P$ is a center if and only if system 
(\ref{e2'}) has an independent from $t$ real analytic integral $F(x,y)=C$. 
Unfortunately, trying to apply the Lyapunov condition gives rise, 
in general, to difficulties comparable with those 
for Poincar\'{e}'s type criterions.

It is important to observe that passing to polar coordinates in (\ref{e2'}) 
we rewrite this in the form
\begin{equation}\label{e3}
\frac{dr}{d\theta}=\frac{P}{1+Q}r
\end{equation}
where 
$P(r,\phi):=\frac{xF(x,y)+yG(x,y)}{r^{2}},\  
Q(r,\phi):=\frac{xG(x,y)-yF(x,y)}{r^{2}},\ (x=r\cos\phi,\ y=r\sin\phi). 
$
If the moduli of the coefficients of $F$ and $G$ are small enough
we can expand the right-hand side of (\ref{e3}) as a series in $r$ to obtain 
an equation of type (\ref{e2}) whose coefficients are trigonometric 
polynomials.
This reduces the Center Problem for (\ref{e2'}) to the Center Problem for a
class of equations (\ref{e2}) whose coefficients depend 
polynomially on the coefficients of (\ref{e2'}).\\
{\bf 1.3.} Let us briefly discuss the content of the present paper.

Let $X_{i}:=L^{\infty}(S^{1})$ be the space of all coefficients $a_{i}$ from
(\ref{e2}), and $X$ be the complex vector space of sequences 
$a=(a_{1},a_{2},\dots)\in\prod_{i\geq 1}X_{i}$ satisfying 
\begin{equation}\label{bound}
\sup_{x\in S^{1}}|a_{i}(x)|\leq l^{i},\ \ \ i=1,2,\dots ,
\end{equation}
for some positive $l$ (depending on $a$).
From Picard iteration it follows that for any $a\in X$ the corresponding
equation (\ref{e2}) is locally solvable for sufficiently small initial
values. Let ${\cal C}\subset X$ be the center set of equation (\ref{e2}),
that is, the set of those $a\in X$ for which the corresponding equations 
(\ref{e2}) determine centers. 
One of our main results (see Theorem \ref{te1}) shows 
that $a\in {\cal C}$ if and only if the monodromy of a certain linear
differential equation $\frac{dF}{dx}=AF$ defined on $S^{1}$
is trivial. Here $A$ is a function on the circle with values in 
the associative 
algebra ${\cal A}(S_{1},S_{2})[[t]]$ of formal power series in $t$ whose
coefficients are complex non-commutative polynomials in variables $S_{1}$
and $S_{2}$ satisfying the relation
$$
[S_{1},S_{2}]:=S_{1}S_{2}-S_{2}S_{1}=-S_{2}^{2}\ .
$$

Next, let us consider the iterated integrals
$$
I_{i_{1},\dots,i_{k}}(a):=
\int\cdots\int_{0\leq s_{1}\leq\cdots\leq s_{k}\leq 2\pi}
a_{i_{k}}(s_{k})\cdots a_{i_{1}}(s_{1})\ ds_{k}\cdots ds_{1}
$$
defined on $X$ (for $k=0$ we assume that this equals 1).
They can be thought of as $k$-linear holomorphic functions on $X$. By
the Ree formula [R] (see also our Section 6.3) the linear space generated by
all such functions is an algebra. A linear combination of iterated
integrals of order $\leq k$  is called an {\em iterated polynomial of degree 
$k$}. 

As an important corollary of Theorem \ref{te1} we obtain an explicit 
description of the center set ${\cal C}$ as a subset of $X$
determined by a system of polynomial equations $c_{n}(a)=0$, $n=1,2,\dots$,
where 
$$
c_{n}(a)=\sum_{i_{1}+\dots + i_{k}=n}c_{i_{1},\dots,i_{k}}
I_{i_{1},\dots,i_{k}}(a)\ , \ \ \ {\rm and}
$$
$$
c_{i_{1},\dots,i_{k}}=
(n-i_{1}+1)\cdot (n-i_{1}-i_{2}+1)\cdot (n-i_{1}-i_{2}-i_{3}+1)\cdots 1\ .
$$

On the other hand, we present a natural ``parametrization'' of the center 
set (see Theorem \ref{param}).
This makes it possible to construct elements from ${\cal C}$ explicitly.

We say that equation (\ref{e2}) corresponding to $a\in X$
determines a {\em universal center}, if $I_{i_{1},\dots,i_{k}}(a)=0$ for all
positive integers $i_{1},\dots, i_{k}$ and $k\geq 1$.
The set ${\cal U}$ of universal centers is, in a sense, a stable part of the
center set ${\cal C}$. In our paper we describe some classes of equations
(\ref{e2}) which determine universal centers. 
Some of our results reveal a connection with the so-called 
{\em composition condition} whose role and importance for the Center Problem
was studied in [Y], [BFY1], [BFY2], [AL] for the special case of Abel
differential equations (cf. Corollary \ref{lorentz} below).

Next, we focus our attention on the most important
algebraic aspects of the Center Problem that were not previously known in
this area.

To this end we introduce a natural multiplication $*:X\times X\rightarrow X$
of elements of $X$ similar to the multiplication of 
continuous paths in the Homotopy Theory (see Section 3.1). 
We say that $a,b\in X$ are
{\em equivalent} (written, $a\sim b$) if $a*b^{-1}\in {\cal U}$. In fact, 
we will show then that $\sim$ is an equivalence relation, that is 
$X$ partitions 
into mutually disjoint equivalence classes. Throughout this paper,
the set of all these  classes will be denoted by $G(X)$. 
An important observation is that the multiplication $*$ induces a similar
multiplication $\cdot :G(X)\times G(X)\rightarrow G(X)$ such that
the pair $(G(X),\cdot)$ is a group. Moreover, we will prove that
the iterated integrals are constant on any equivalence class. 
Thus without loss of generality they may be considered as functions on 
$G(X)$. Also, it could be extracted from the construction of $G(X)$ that
these functions separate points on $G(X)$. 

Further, let us provide $G(X)$ with the weakest topology $\tau$ in which all 
iterated integrals (considered as functions on $G(X)$) are continuous. 
One of the central results of Section 3.1 says that $(G(X),\cdot,\tau)$ is a 
separable topological group satisfying the following properties
\begin{itemize}
\item[(1)]
$(G(X),\tau)$ is contractible
to a point, arcwise connected, locally arcwise and simply connected;
\item[(2)]
$G(X)$ is
{\em residually torsion free nilpotent} (that is, the set of all
finite-dimensional unipotent representations separates elements of $G(X)$).
\end{itemize}

To formulate our next result we introduce the set $G_{c}[[r]]$ of 
complex power series of the form
$f(r)=r+\sum_{i=1}^{\infty}d_{i}r^{i+1}$ each convergent in some open 
neighbourhood of $0\in\Co$. Consider the functions 
$d_{i}:G_{c}[[r]]\rightarrow\Co$ such that $d_{i}(f)$ equals the 
$i+1$-st coefficient of the Taylor expansion of $f$ at $0$. Let
$\tau'$ be the weakest topology on $G_{c}[[r]]$ in which all the functions
$d_{i}$ are continuous. Also, we consider the multiplication on
$G_{c}[[r]]$ defined by the composition of series. 
Then, as in the case of $G(X)$,
 we will prove that $(G_{c}[[r]],\circ,\tau')$ is a 
separable topological group satisfying the properties
\begin{itemize}
\item[(3)]
$(G_{c}[[r]],\tau')$ is contractible to a 
point, arcwise connected, locally arcwise and simply connected;
\item[(4)]
$G_{c}[[r]]$ is residually torsion free nilpotent.
\end{itemize}

Now, for any $a\in X$ by $v(x;r;a)$, $x\in [0,2\pi]$, we denote the 
Lipschitz solution of equation (\ref{e2}) corresponding to $a$ with 
initial value $v(0;r;a)=r$. It is clear that for every $x$ we have 
$v(x;r;a)\in G_{c}[[r]]$. The function $P(a)(r):=v(2\pi;r;a)$ is called {\em
the first return map} of equation (\ref{e2}). In Section 3.2 we will show 
that
$$
P(a*b)=P(b)\circ P(a)\ .
$$
This fundamental property together with the fact that $P(a)(r)\equiv r$ for 
any $a\in {\cal U}$ imply
that there exists a map $\widehat P:G(X)\rightarrow G_{c}[[r]]$
such that $\widehat P([a]):=P(a)$ where $[a]$ denotes the equivalence class
containing $a\in X$. Then we will prove that
\begin{itemize}
\item[(5)]
$\widehat P$ is a surjective homomorphism of topological groups;
\item[(6)] The kernel $\widehat {\cal C}\subset G(X)$ of $\widehat P$
coincides with the image of the center
set ${\cal C}\subset X$ in $G(X)$;
\item[(7)]
$(\widehat {\cal C},\tau)$ is contractible to a 
point, arcwise connected,  locally arcwise and simply connected.
\end{itemize}

Let us consider the quotient group $Q(X)=G(X)/\widehat {\cal C}$. By
$\pi:G(X)\rightarrow Q(X)$ we denote the quotient homomorphism. 
It follows from the preceding discussion that every function 
$c_{i}$ in the definition of the center set ${\cal C}$
satisfies $c_{i}(a)=c_{i}(b)$ for $a\sim b$. Therefore every $c_{i}$ can be 
considered as a continuous function on $G(X)$. In fact, we will show that 
these functions are constant on any fibre of the map $\pi$ and hence they 
determine functions $\overline{c}_{i}:Q(X)\rightarrow\Co$. 
Let $\tau''$ be the
weakest topology on $Q(X)$ in which all the functions $\overline{c}_{i}$
are continuous. Then we will prove that
\begin{itemize}
\item[(8)]
$(Q(X),\tau'')$ is a topological group.
\item[(9)]
The homomorphism $\widehat P$ determines an isomorphism  
$\overline{P}:Q(X)\rightarrow G_{c}[[r]]$ of topological groups defined by 
$\widehat P=\overline{P}\circ\pi$.
\end{itemize}

Also, we will show that there is a continuous map 
$T:G_{c}[[r]]\rightarrow G(X)$ such that $\widehat P\circ T=id$. In 
particular, we obtain that
\begin{itemize}
\item[(10)] The map 
$\widetilde T:G_{c}[[r]]\times\widehat {\cal C}\rightarrow G(X)$,
$\widetilde T(f,g):=T(f)\cdot g$, is a homeomorphism.
\end{itemize}

Finally, in Section 3.3 we study some algebraic properties of the group 
$G(X)$ and its specific subgroups defined over certain
subfields of the filed $\Co$.\\
{\em Acknowledgment.} I would like to thank Professor Y. Yomdin for
very fruitful discussions during his visit to Calgary. I also would like to
thank Professors L. Bos and Y. Brudnyi for the valuable help and the support
of this work.

The next two sections contain a more leisurely description of the main
results of the paper.
\sect{\hspace*{-1em}. Monodromy in the Center Problem.}
{\bf 2.1.} Our first result shows that the Center Problem for (\ref{e2})
is equivalent to the triviality of the monodromy of a certain linear 
system of ODEs.

In order to formulate this result one introduces the associative algebra
${\cal A}(S_{1},S_{2})$ with unit $I$
of non-commutative polynomials with complex coefficients in variables 
$S_{1}$ and $S_{2}$ satisfying the relation 
\begin{equation}\label{relat}
[S_{1},S_{2}]:=S_{1}S_{2}-S_{2}S_{1}=-S_{2}^{2}\ .
\end{equation}
Then ${\cal A}(S_{1},S_{2})[[t]]$ stands for the associative algebra of 
formal power series in $t$ whose coefficients are elements from
${\cal A}(S_{1},S_{2})$. For an element $a=(a_{1},a_{2},\dots)\in X$ (see
Section 1.3) let us consider the equation on $S^{1}$
\begin{equation}\label{formal}
\frac{dF}{dx}=(\sum_{i=1}^{\infty}a_{i}(x)t^{i}S_{1}S_{2}^{i-1})F\ .
\end{equation}
This can be solved locally by the Picard iteration method
to obtain a local solution $F(x)$ as a function in $x$ with
values in the group $G(S_{1},S_{2})[[t]]$ of invertible elements of
${\cal A}(S_{1},S_{2})[[t]]$ whose coefficients are Lipschitz. As usual,
the monodromy of (\ref{formal}) is a homomorphism 
$\rho:\Z\rightarrow G(S_{1},S_{2})[[t]]$ where $\Z$ is the fundamental
group of the unit circle $S^{1}$. It is the only obstruction to extending 
local solutions of (\ref{formal}) to global ones.
\begin{Th}\label{te1}
Equation (\ref{e2}) determines a center if and only if the monodromy of
the corresponding equation (\ref{formal}) is trivial.
\end{Th}

From this result we obtain the following characterization of centers for
(\ref{e2}). As before, let $X$ be the complex vector space of sequences
$a=(a_{1},a_{2},\dots)$ satisfying (\ref{bound}).
The functions of the form
$$
I_{i_{1},\dots,i_{k}}(a):=
\int\cdots\int_{0\leq s_{1}\leq\cdots\leq s_{k}\leq 2\pi}
a_{i_{k}}(s_{k})\cdots a_{i_{1}}(s_{1})\ ds_{k}\cdots ds_{1}
$$
will be called the {\em basic iterated integrals} on $X$.

Now, for a sufficiently small $r$, let $v(x;r;a)$, $x\in [0,2\pi]$, be 
the Lipschitz solution of equation (\ref{e2}) corresponding to $a\in X$
with initial value $v(0;r;a)=r$. Then $P(a)(r):=v(2\pi;r;a)$ is the
first return map of this equation.
\begin{Th}\label{center}
For sufficiently small initial values $r$ the first return map $P(a)$ is an
absolutely convergent power series 
$P(a)(r)=r+\sum_{n=1}^{\infty}c_{n}(a)r^{n+1}$,
where
$$
c_{n}(a)=\sum_{i_{1}+\dots + i_{k}=n}c_{i_{1},\dots,i_{k}}
I_{i_{1},\dots,i_{k}}(a)\ , \ \ \ {\rm and}
$$
$$
c_{i_{1},\dots,i_{k}}=
(n-i_{1}+1)\cdot (n-i_{1}-i_{2}+1)\cdot (n-i_{1}-i_{2}-i_{3}+1)\cdots 1\ .
$$

The center set ${\cal C}\subset X$ of equation (\ref{e2}) is determined
by the system of polynomial equations $c_{n}(a)=0$, $n=1,2,\dots$ .
\end{Th}
\begin{C}\label{dimens}
(a)\ \ \ $c_{n}(a)=I_{n}(a)+f_{n}(a)$ where 
$I_{n}(a):=\int_{0}^{2\pi}a_{n}(s)\ ds$ and $f_{n}$ is an iterated 
polynomial of degree $n$ in $a_{1},\dots, a_{n-1}$;\\
(b)\ \ \ The set
$$
{\cal C}_{n}=\{a\in X\ :\ c_{1}(a)=c_{2}(a)=\dots =c_{n}(a)=0\}
$$ 
is a closed complex submanifold of $X$ of codimension $n$
containing $0\in X$;\\
(c)\ \ \ The tangent space to ${\cal C}_{n}$
at $0$ is determined by equations $I_{1}(a)=\dots=I_{n}(a)=0$.
\end{C}

Another characterization of centers of (\ref{e2}) is given by the
following theorem.
\begin{Th}\label{param}
An element $a=(a_{1},a_{2},\dots)\in X$ belongs to the center set 
${\cal C}$ if and only if
there is a sequence $u_{1},u_{2},\dots$ of $2\pi$-periodic Lipschitz 
functions such that $u_{i}(0)=0$ for any $i$ and  
$$
\sum_{i=1}^{\infty}a_{i}(x)t^{i+1}=-
\frac{\sum_{k=1}^{\infty}u_{k}'(x)t^{k+1}}
{1+\sum_{k=1}^{\infty}(k+1)u_{k}(x)t^{k}}
$$
as formal power series in $t$.
\end{Th}

This result gives a ``parametrization'' of the center set 
${\cal C}\subset X$ but leaves open the question on existence 
of such a sequence $u_{1},u_{2}\dots $ for a specific $a\in X$. \\
{\bf 2.2.} Based on Theorem \ref{te1} we split the Center Problem for
(\ref{e2}) into two parts. 

To this end we introduce the associative algebra
${\cal A}(X_{1},X_{2})$ with unit $I$ of non-commutative polynomials
with complex coefficients in free non-commutative variables $X_{1}$ and 
$X_{2}$.
By ${\cal A}(X_{1},X_{2})[[t]]$ we denote the associative algebra of formal
power series in $t$ whose coefficients are elements from 
${\cal A}(X_{1},X_{2})$. 
Then there is a natural homomorphism 
$\phi:{\cal A}(X_{1},X_{2})[[t]]\rightarrow {\cal A}(S_{1},S_{2})[[t]]$
uniquely defined by the relations $\phi(X_{1})=S_{1}$ and 
$\phi(X_{2})=S_{2}$. Next, for $a=(a_{1},a_{2},\dots)\in X$ as in 
(\ref{formal}) let us consider the equation on the circle $S^{1}$
\begin{equation}\label{formal1}
\frac{dF}{dx}=(\sum_{i=1}^{\infty}a_{i}(x)t^{i}X_{1}X_{2}^{i-1})F\ .
\end{equation}
Again one can solve this equation locally by Picard iteration. Then
the monodromy $\widetilde\rho:\Z\rightarrow G(X_{1},X_{2})[[t]]$ of
(\ref{formal1}) 
is the homomorphism of the fundamental group of $S^{1}$ into the 
group $G(X_{1},X_{2})[[t]]$ of invertible elements of 
${\cal A}(X_{1},X_{2})[[t]]$. It is clear that 
\begin{equation}\label{decom}
\rho=\phi\circ\widetilde\rho\ .
\end{equation}
Then from Theorem \ref{te1} follows
\begin{C}\label{univ}
If the monodromy $\widetilde\rho$ is trivial, 
then the corresponding equation (\ref{e2}) determines a center.
\end{C}

We will say that equation (\ref{e2}) determines a {\em universal center}, if
the monodromy of the corresponding equation (\ref{formal1}) is trivial.
The set ${\cal U}$ of universal centers is, in a sense, a stable part of the
center set ${\cal C}$. As we will see, ${\cal U}\neq {\cal C}$, in general.
The next two sections describe some classes of equations
(\ref{e2}) which determine universal centers. \\
{\bf 2.3.} Let $\omega(x):=\sum_{i=1}^{\infty}a_{i}(x)t^{i}X_{1}X_{2}^{i-1}$. 
As before we identify functions on the circle with $2\pi$-periodic functions 
on $\Re$. Recall that the fundamental solution of (\ref{formal1}) is a 
map $F:\Re\rightarrow G(X_{1},X_{2})[[t]]$ such that $F(0)=I$ and 
$F'(x)=\omega(x)\cdot F(x)$. It can be
presented by Picard iteration (cf. [Na]) in the form
\begin{equation}\label{e3'}
F(x):=I+\sum_{k=1}^{\infty}
\int\cdots\int_{0\leq s_{1}\leq\cdots\leq s_{k}\leq x}
\omega(s_{k})\cdots\omega(s_{1})\ ds_{k}\cdots ds_{1}
\end{equation}

Also, it is easy to see that
\begin{equation}\label{e4'}
F(x)=\sum_{i=0}^{\infty}f_{i}(x;X_{1},X_{2})t^{i}
\end{equation}
where $f_{0}=I$ and the other $f_{i}$ are homogeneous polynomials of degree
$i$ in $X_{1}$ and $X_{2}$ whose coefficients are locally Lipschitz  
functions in $x\in\Re$. Then the monodromy of (\ref{formal1}) is defined as
$\widetilde\rho(n):=F(2\pi n)=F(2\pi)^{n}$, $n\in\Z$\ .
\begin{Proposition}\label{integr}
The monodromy of (\ref{formal1}) is trivial, i.e. $F(2\pi)=I$, if and 
only if for all positive integers $i_{1},\dots, i_{k}$ and $k\geq 1$
\begin{equation}\label{vanish}
\int\cdots\int_{0\leq s_{1}\leq\cdots\leq s_{k}\leq 2\pi}
a_{i_{k}}(s_{k})\cdots a_{i_{1}}(s_{1})\ ds_{k}\cdots ds_{1}=0\ .
\end{equation}
In particular then,
$$
\widetilde a_{i}(x):=\int_{0}^{x}a_{i}(s)\ ds,\ \ \ i=1,2,\dots ,
$$
are $2\pi$-periodic Lipschitz functions.
\end{Proposition}

Next, let us consider equations 
\begin{equation}\label{neq}
F'(x)=\omega_{n}(x)\cdot F(x)\ \ \ {\rm with}\ \ \
\omega_{n}:=\sum_{i=1}^{n}a_{i}(x)t^{i}X_{1}X_{2}^{i-1}\ .
\end{equation}
By $\widetilde\rho_{n}:\Z\rightarrow G(X_{1},X_{2})[[t]]$ we denote
the monodromy of (\ref{neq}). 
\begin{Proposition}\label{integr1}
The monodromy of (\ref{formal1}) is trivial if and only if all
$\widetilde\rho_{n}$ are trivial. Moreover, the triviality
of $\widetilde\rho_{n}$ is equivalent to the fulfilment of equations
(\ref{vanish}) for any integers $1\leq i_{1},\dots,i_{k}\leq n$,
and any $k$.
\end{Proposition}
{\bf 2.4.} Let us introduce the Lipschitz map 
$A_{n}:S^{1}\rightarrow\Co^{n}$,
$A_{n}(x)=(\widetilde a_{1}(x),\dots,\widetilde a_{n}(x))$,
and set $\Gamma_{n}:=A_{n}(S_{1})$. We require
\begin{D}\label{de1}
The polynomially convex hull $\widehat K$ of a compact set
$K\subset\Co^{n}$ is the set of points $z\in\Co^{n}$
such that if $p$ is any holomorphic polynomial in $n$ variables
$$
|p(z)|\leq\max_{x\in K}|p(x)|\ .
$$
\end{D}

It is well known (see e.g. [AW]) that $\widehat K$ is compact, and if
$K$ is connected then $\widehat K$ is connected.
\begin{Th}\label{contr}
If the monodromy $\widetilde\rho_{n}$ is trivial then for 
any domain $U$ containing $\widehat\Gamma_{n}$ the path 
$A_{n}:S^{1}\rightarrow U$ is contractible in $U$ to a point.
\end{Th}

Since $A_{n}$ is Lipschitz, $\Gamma_{n}$ is of a finite linear measure.
Then according to the remarkable result of
Alexander [A], $\widehat\Gamma_{n}\setminus\Gamma_{n}$ is a (possibly
empty) pure one-dimensional analytic subset of $\Co^{n}\setminus\Gamma_{n}$.
In particular, since the covering dimension of $\Gamma_{n}$ is 1,
the covering dimension of $\widehat\Gamma_{n}$ is 2. However,
in general we do not know how to use this to get from
Theorem \ref{contr} more information about $\Gamma_{n}$. Thus we restrict
our presentation to several special cases.
\begin{C}\label{coro1}
Suppose $\Gamma_{n}$ is triangulable and $\widehat\Gamma_{n}=\Gamma_{n}$.
If the monodromy $\widetilde\rho_{n}$ is trivial then 
the path $A_{n}:S^{1}\rightarrow\Gamma_{n}$ is contractible inside 
$\Gamma_{n}$ to a point. Moreover, the contractibility of $A_{n}$ is 
equivalent to the 
factorization $A_{n}=A_{1n}\circ A_{2n}$ where $A_{2n}: 
S^{1}\rightarrow G_{n}$ is
a continuous map into a finite tree $G_{n}\subset\Re^{N_{n}}$, and 
$A_{1n}:G_{n}\rightarrow\Gamma_{n}$ is a finite continuous map.
\end{C}

The converse result requires a much
stronger than just triangulability condition on $\Gamma_{n}$ that is
described by 
\begin{D}\label{liptriang}
A compact curve $C\subset\Re^{N}$ is called Lipschitz triangulable if\\
(a)\ \ $C=\cup_{j=1}^{s}C_{i}$ and for $i\neq j$ the intersection
$C_{i}\cap C_{j}$ consists of at most one point;\\
(b)\ \ There are Lipschitz embeddings $f_{i}:[0,1]\rightarrow\Re^{N}$ such
that $f_{i}([0,1])=C_{i}$;\\
(c)\ \ The inverse maps $f_{i}^{-1}:C_{i}\rightarrow\Re$ are locally 
Lipschitz on $C_{i}\setminus (f_{i}(0)\cup f_{i}(1))$.
\end{D}
\begin{Th}\label{pr1'}
Suppose $\Gamma_{n}$ is Lipschitz triangulable,
$\widehat\Gamma_{n}=\Gamma_{n}$, and 
$A_{n}:S^{1}\rightarrow\Gamma_{n}$ is contractible inside
$\Gamma_{n}$ to a point. Suppose also that  $A_{n}^{-1}(x)$ is 
countable for any $x\in\Gamma_{n}$. Then the corresponding monodromy 
$\widetilde\rho_{n}$ is trivial.
\end{Th}
\begin{R}\label{re1}
{\rm  We will see from the proof that under the hypotheses of 
Theorem \ref{pr1'} the map $A_{2n}:S^{1}\rightarrow\Re^{N_{n}}$ is 
Lipschitz, $G_{n}$ is Lipschitz triangulable, and 
$A_{1n}:G_{n}\rightarrow\Gamma_{n}$ is locally Lipschitz outside a finite 
set.}
\end{R}
\begin{E}\label{ex1}
{\rm (1) If $A_{n}:S^{1}\rightarrow\Co^{n}$ is non-constant analytic, then
$\Gamma_{n}$ is Lipschitz triangulable, and $A_{n}^{-1}(x)$
is finite for any $x\in\Gamma_{n}$;\\
(2) To satisfy the second hypothesis of Corollary \ref{coro1}, one can 
suppose, e.g., that $\Gamma_{n}$ belongs to a compact $K_{n}$ in a 
$C^{1}$-smooth manifold $M_{n}$ with no 
complex tangents such that $\widehat K_{n}=K_{n}$. For instance, one can 
take any compact $K_{n}$ in $M_{n}=\Re^{n}$ (for the proof see e.g. 
[AW, Th.17.1]).}
\end{E}

Next, we consider maps into one-dimensional complex spaces. Suppose that
one of the following two conditions is satisfied:
\\
(1)\ \ $\Gamma_{n}\subset X$ where $X$ is a closed one-dimensional 
complex analytic subset of a domain $U\subset\Co^{n}$ such that 
$U=\cup_{j}K_{j}$
with $K_{j}\subset\subset K_{j+1}$, and $\widehat K_{j}=K_{j}$ for any
$j$;\\
(2)\ \ $\Gamma_{n}\subset X$ where $X\subset\Co^{n}$ is a connected 
one-dimensional complex space with  $dim_{\Co}H^{1}(X,\Co)<\infty$, and 
for each $\delta\in H^{1}(X,\Co)$ there is a holomorphic 1-form with 
polynomial coefficients  $\omega_{\delta}$ such that 
$\delta(c)=\int_{c}\omega_{\delta}$ for any $c\in H_{1}(X,\Co)$.
\begin{C}\label{coro2}
Let $\Gamma_{n}\subset X$ where $X$ satisfies either condition 
(1) or (2).
Suppose that there is a continuous map 
$\widetilde A_{n}:R\rightarrow X$ of an open neighbourhood 
$R\subset\Co$ of $S^{1}$ such that $\widetilde A_{n}|_{S^{1}}=A_{n}$, and
$\widetilde A_{n}^{-1}(x)$ is finite for any $x\in\widetilde A_{n}(R)$.
Then the monodromy $\widetilde\rho_{n}$ is
trivial if and only if $A_{n}=A_{1n}\circ A_{2n}$ where 
$A_{2n}:S^{1}\rightarrow\Di$ is a continuous map into the unit disk 
$\Di\subset\Co$, locally Lipschitz outside a finite set,
and $A_{1n}:\Di\rightarrow\Co^{n}$ is a finite holomorphic
map.
\end{C}
\begin{R}\label{rema2}
{\rm (1) In the proof of Corollary \ref{coro2} we give separate arguments
for assumptions (1) and (2). However, it is possible to show that condition
(1) can be reduced to the case of condition (2).\\
(2) We will also show that if  $\widetilde A_{n}$ is holomorphic then
$A_{2n}$ can be extended to a holomorphic map of an open neighbourhood
of $S^{1}$ into $\Di$.\\
(3) Note that the hypotheses of the corollary
are valid for $X$ a complex algebraic curve.\\
(4) In the case when $A_{n}:S^{1}\rightarrow X$ is a Lipschitz embedding 
and $\Gamma_{n}=A_{n}(S^{1})$ lies outside singularities of 
$X$, one can obtain 
under the hypotheses of Corollary \ref{coro2} that $A_{n}=A_{1n}\circ A_{2n}$ 
where $A_{2n}: S^{1}\rightarrow\Di$ is a Lipschitz embedding
and $A_{1n}:\Di\rightarrow\Co^{n}$ is a holomorphic map, one-to-one outside
a finite set.}
\end{R}
\begin{C}\label{lorentz}
Suppose that the coefficients $a_{1},\dots, a_{n}$ in (\ref{neq}) are 
trigonometric polynomials. Then the monodromy $\widetilde\rho_{n}$ is 
trivial if 
and only if there is a trigonometric polynomial $q$ and 
polynomials $p_{1},\dots, p_{n}\in\Co[z]$ such that 
$$
\widetilde a_{i}(x)=p_{i}(q(x)),\ \ \ x\in S^{1},\ \ \ 1\leq i\leq n\ .
$$
\end{C}
{\bf 2.5.} Finally, we formulate two finiteness results showing that in
certain cases the triviality of some $\widetilde\rho_{N}$ implies that of
$\widetilde\rho$.
\begin{Th}\label{tef1}
Suppose that the coefficients $a_{1},a_{2},\dots,$ in equation (\ref{formal1})
are continuous functions such that each 
$\widetilde a_{i}:=\int_{0}^{x}a_{i}(s)ds$ is $2\pi$-periodic. 
Let $A_{k}:S^{1}\rightarrow\Co^{k}$, 
$A_{k}(x):=(\widetilde a_{1}(x),\dots,\widetilde a_{k}(x))$, and
$\Gamma_{k}=A_{k}(S^{1})$. Suppose
that there is an integer $n$ such that $\Gamma_{n}$ is a 
piecewise smooth curve, the number of critical values of $A_{n}$ is
finite, and $A_{n}^{-1}(x)$ is finite for any 
$x\in\Gamma_{n}$. If also $\widehat\Gamma_{k}=\Gamma_{k}$ for any 
$k$, then there is an integer $N\geq 1$ such that the triviality of 
$\widetilde\rho_{N}$ implies the triviality of $\widetilde\rho$.
\end{Th}
\begin{R}\label{rea}
{\rm The hypotheses of the theorem are fulfilled if, e.g., all $a_{i}$ are
real continuous functions and some $\widetilde a_{n}$ is a non-zero analytic 
function.}
\end{R} 
\begin{Th}\label{tef2}
Suppose that there are $2\pi$-periodic Lipschitz functions
$b_{1},\dots,b_{k}$ such that each coefficient $a_{i}$ in equation 
(\ref{formal1}) is the uniform limit of functions of the form
$\sum_{j=1}^{k}p_{ji}(b_{1},\dots,b_{k})\cdot b_{j}'$, where
$p_{ji}\in\Co[z_{1},\dots,z_{k}]$ are holomorphic polynomials and
$b_{j}'$ is the derivative of $b_{j}$.
Then the triviality of the monodromy of
the equation $F'(x)=(\sum_{i=1}^{k}b_{i}'(x)t^{i}X_{1}X_{2}^{i-1})\cdot F(x)$
implies that the monodromy $\widetilde\rho$ for (\ref{formal1}) is trivial.
\end{Th}
\begin{C}\label{planar}
Let now $H(x,y)\in\Co[x,y]$ be a homogeneous polynomial. For any holomorphic
functions $P_{1}, P_{2}$  defined in an open neighbourhood 
of $0\in\Co$ we define $A(x,y):=P_{1}(H(x,y))$, and 
$B(x,y):=P_{2}(H(x,y))$. Then the vector field
$$
\left\{
\begin{array}{l}
\displaystyle
\dot x=-y-xy\frac{\partial A(x,y)}{\partial x}+
x^{2}\frac{\partial A(x,y)}{\partial y}-yB(x,y)\\
\\
\displaystyle
\dot y=x-y^{2}\frac{\partial A(x,y)}{\partial x}+
xy\frac{\partial A(x,y)}{\partial y}+xB(x,y)
\end{array}
\right .
$$
determines a center.
\end{C}

The next result deals with symmetries of vector fields. Since the
center problem is invariant under rotations, it suffices to
consider a particular case of the symmetry. This result first was
obtained by Poincar\'{e} [P] by a different argument.
\begin{C}\label{sim}
Let $F(x,y),G(x,y)$ be analytic functions defined in an open neighbourhood
of $0\in\Re^{2}$ such that their Taylor expansions at 0 do not contain
constant and linear terms. Suppose 
$F(x,-y)=-F(x,y)$ and $G(x,-y)=G(x,y)$. Then the vector field \
$\dot x=-y+F(x,y),\ \dot y=x+G(x,y)$ determines a center.
\end{C}
\begin{E}\label{cherkas}
{\rm The field (\ref{e1}) with homogeneous $F(x,y)$ and 
$G(x,y)$ of degree 2 can be written in the Dulac-Kapteyn form  
with 5 parameters instead of 6 by using an appropriate rotation of the plane:}
\begin{equation}\label{dva}
\left\{
\begin{array}{l}
\dot x=-y-\lambda_{3}x^{2}+(2\lambda_{2}+\lambda_{5})xy+\lambda_{6}y^{2}\\
\dot y=x+\lambda_{2}x^{2}+(2\lambda_{3}+\lambda_{4})xy-\lambda_{2}y^{2}\ .
\end{array}
\right .
\end{equation}
{\rm It has been established by H. Dulac that the center set for
(\ref{dva}) consists of 4 components described as follows (cf. [Si]).}\\
(1)\ Lotka-Volterra component:\  $\lambda_{3}=\lambda_{6}$;\\
(2)\ Symmetric component:\ $\lambda_{2}=\lambda_{5}=0$;\\
(3)\ Hamiltonian component:\ $\lambda_{4}=\lambda_{5}=0$;\\
(4)\ Darboux component:\ $\lambda_{5}=\lambda_{4}+5\lambda_{3}-5\lambda_{6}=
\lambda_{3}\lambda_{6}-2\lambda_{6}^{2}-\lambda_{2}^{2}=0$.

{\rm Next, passing to polar coordinates in (\ref{dva}) we obtain}
\begin{equation}\label{twodim}
\frac{dr}{d\phi}=\frac{f(\phi)r^{2}}{1+g(\phi)r}
\end{equation}
{\rm where  $f(\phi):=\frac{xF(x,y)+yG(x,y)}{r^{3}}$, 
$g(\phi):=\frac{xG(x,y)-yF(x,y)}{r^{3}}$, $x=r\cos\phi$, $y=r\sin\phi$.
Let us investigate the universal center
conditions for equation (\ref{twodim}). If (\ref{twodim}) determines a 
universal center then the first integrals $I(f)$ and
$I(fg)$ of $f$ and $fg$ satisfy the relations of Corollary \ref{lorentz}. 
This easily implies
that $I(f)$ and $g$ satisfy the same relations. Going back to the functions 
$F$ and $G$
we obtain a planar polynomial vector field of the form of Corollary 
\ref{planar}. Now $H$ is a linear homogeneous polynomial,
and $P_{1},P_{2}$ are linear polynomials without constant terms. 
A simple computation shows that in this case for the equation (\ref{dva}) we
have}\\
(A)\ $\lambda_{3}-\lambda_{6}=4\lambda_{2}+\lambda_{5}=
4\lambda_{3}+\lambda_{4}=0$ (a linear two-dimensional subspace in the 
Lotka-Volterra component);\\
(B)\ $\lambda_{2}=\lambda_{5}$ (Symmetric component).
\end{E}
\begin{R}\label{che}
{\rm We can further simplify equation (\ref{twodim}) by
the application of the Cherkas transformation [Che],
$r(\phi)=\frac{\rho(\phi)}{1-g(\phi)\rho(\phi)}$. Then 
we get the Abel differential equation}
\begin{equation}\label{ab}
\frac{d\rho}{d\phi}=p(\phi)\rho^{2}+q(\phi)\rho^{3}
\end{equation}
{\rm where $p(\phi)=f(\phi)+g'(\phi)$ and $q(\phi)=-f(\phi)g(\phi)$. 
Moreover, since the Cherkas transformation and its inverse are regular at
$r=\rho=0$, the center sets of (\ref{twodim}) and (\ref{ab}) coincide.
Let $P(\phi)=\int_{0}^{\phi}p(s)ds$ and $Q(\phi)=\int_{0}^{\phi}q(s)ds$.
In case equation (\ref{ab}) determines a center, $P$ and $Q$ are 
trigonometric polynomials.
It was shown by Blinov [Bl] that the Hamiltonian and the Symmetric 
components after
performing the Cherkas transformation determine universal centers for
(\ref{ab}) and the other two components do not. Here
for the Hamiltonian component we have $Q=-\frac{3}{16}P^{2}+
\frac{\lambda_{2}}{2}P$, and for the Symmetric component $P$ and $Q$ are
polynomials in $\sin(\phi)$. Then the required result follows from
Corollary \ref{lorentz}.} 

{\rm Thus the Hamiltonian component (which is a universal center for
(\ref{ab})) can be obtained as a result of an additional non-linear change
of variables in (\ref{twodim}). }
\end{R}

The results of Section 2 will be proved in Sections 4-8.
\sect{\hspace*{-1em}. Algebraic Model for the Center Problem.}
{\bf 3.1.} Now we consider the Center Problem for (\ref{e2}) from an
algebraic point of view. 

Let us introduce a multiplication $*:X\times X\rightarrow X$ as follows.

Given $a=(a_{1},a_{2},\dots)$ and $b=(b_{1},b_{2},\dots)$ 
from $X$ we define 
$$
a*b=(a_{1}*b_{1},a_{2}*b_{2},\dots)\in X\ \ \
{\rm and}\ \ \ a^{-1}=(a_{1}^{-1},a_{2}^{-1},\dots)\in X
$$
where for any $i\in\N$,
$$
(a_{i}*b_{i})(t)=\left\{
\begin{array}{ccc}
2a_{i}(2t)&{\rm if}&0<t\leq \pi\\
2b_{i}(2t-2\pi)&{\rm if}&\pi<t\leq 2\pi
\end{array}
\right.
$$
and
$$
a_{i}^{-1}(t)=-a_{i}(2\pi-t)\ ,\ \ \ 0<t\leq 2\pi\ .
$$

From Proposition \ref{integr} it follows that equation (\ref{e2}) 
corresponding to $a\in X$ determines a {\em universal center}, if all 
non-constant basic iterated integrals vanish at $a$, i.e., 
$I_{i_{1},\dots,i_{k}}(a)=0$ for all possible
$i_{1},\dots, i_{k}, k\in\N$.

We say that $a,b\in X$ are {\em equivalent} (written, $a\sim b$) if
$a*b^{-1}\in {\cal U}$.
\begin{Proposition}\label{equivrelat}
The relation $\sim$ is reflexive,
symmetric and transitive so that $X$ partitions into mutually disjoint
equivalence classes.
\end{Proposition}
By $G(X)$ we denote the set of all the equivalence classes, and by
$[a]$ the class containing $a\in X$. The
element $a$ will be called a {\em representative} of the class $[a]$.
\begin{Proposition}\label{oper}
The following relations hold:
\begin{itemize}
\item[{\rm (1)}]
if $a\sim a'$ and $b\sim b'$ then $a*b\sim a'*b'$;
\item[{\rm (2)}]
$(a*b)*c\sim a*(b*c)$;
\item[{\rm (3)}]
$a*0\sim 0*a$;
\item[{\rm (4)}]
$a*a^{-1}\sim 0\sim a^{-1}*a$ .
\end{itemize}
\end{Proposition}

From (1) it follows that the class $[a*b]$ depends only on classes $[a]$ and
$[b]$. Therefore the formula
$$
[a]\cdot [b]=[a*b]
$$
determines a multiplication in $G(X)$. Also, from properties (2), (3) and (4)
we obtain that $G(X)$ is a group with respect to this multiplication.
The class $[0]=e$ is the unit of the group $G(X)$ and the class $[a^{-1}]$
is the inverse to the class $[a]$.
\begin{Proposition}\label{sepit}
If $a,b\in X$ are such that $a\sim b$, then $I(a)=I(b)$ for
any iterated integral $I$. 
\end{Proposition}
Thus every such $I$ determines a function
$\widehat I: G(X)\rightarrow\Co$, $\widehat I([a]):=I(a)$, which also will be
referred to as an iterated integral on $G(X)$. The algebra of iterated
integrals on $G(X)$ will be denoted by ${\cal I}(G(X))$.

Next we introduce a topology $\tau$ on $G(X)$ as follows.

Given $g\in G(X)$ the base of open neighbourhoods at $g$ consists of all
possible sets of the form
$$
\{h\in G(X)\ :\ |\widehat I_{1}(h)-\widehat I_{1}(g)|<\epsilon,\dots,
|\widehat I_{k}(h)-\widehat I_{k}(g)|<\epsilon\}
$$
where $I_{1},\dots, I_{k}$ are basic iterated integrals and 
$\epsilon$ is a positive number. 
This is the {\em weakest topology} in which all functions $\widehat I$ are
continuous. 

We say that a group $G$ is {\em residually torsion free nilpotent} if the 
set of finite-dimensional unipotent representations separates elements of 
$G$. (For instance, a free group with at most countable number of 
generators satisfies this property.)

The following result describes certain properties of $G(X)$.
\begin{Th}\label{group}
\begin{itemize}
\item[{\rm (1)}]
The family ${\cal I}(G(X))$ separates points on $G(X)$. 
\item[{\rm (2)}]
$(G(X), \cdot,\tau)$ is a separable topological group.
\item[{\rm (3)}]
$G(X)$ is metrizable and is the union of an increasing sequence of
compact subsets.
\item[{\rm (4)}]
$G(X)$ is contractible to a point, arcwise connected, locally simply
and arcwise connected.
\item[{\rm (5)}]
$G(X)$ is residually torsion free nilpotent.
\end{itemize}
\end{Th}
{\bf 3.2.}
Let $G_{c}[[r]]$ be the set of complex power series of the form
$f(r)=r+\sum_{k=1}^{\infty}d_{k}r^{k+1}$ each convergent in some open 
neighbourhood of $0\in\Co$. It is easy to see that $G_{c}[[r]]$ is
a group with the multiplication defined by the composition 
of series. Consider the functions $d_{k}:G_{c}[[r]]\rightarrow\Co$ such
that $d_{k}(f)$ equals the $k+1$-st coefficient of the Taylor expansion of
$f$ at $0$. Let $\tau'$ be the weakest topology on $G_{c}[[r]]$ in
which all functions $d_{k}$ are continuous. Then  we prove
\begin{Proposition}\label{ser1}
$(G_{c}[[r]],\circ,\tau')$ is a separable topological group. 
The family of functions $d_{k}$ separates points on $G_{c}[[r]]$.
Moreover, $G_{c}[[r]]$ is contractible to a point, arcwise connected,
locally simply and arcwise connected, and residually torsion free nilpotent.
\end{Proposition}

Next we will confine our attention to the first return map $P(a)$ 
of equation (\ref{e2}) corresponding  to $a\in X$. 
It is clear that $P(a)\in G_{c}[[r]]$. The following property of the
map $P$ is the core of our construction.
\begin{Th}\label{comp1}
The identity 
$$
P(a*b)=P(b)\circ P(a)
$$
holds for any $a,b\in X$.
\end{Th}

From here and the fact that $P(a)(r)\equiv r$ for any $a\in {\cal U}$
we obtain that there exists a map $\widehat P:G(X)\rightarrow G_{c}[[r]]$
such that $\widehat P([a]):=P(a)$. Then we will prove
\begin{Th}\label{comp2}
The map $\widehat P:G(X)\rightarrow G_{c}[[r]]$ is a surjective homomorphism
of topological groups. The kernel $\widehat {\cal C}\subset G(X)$ of
$\widehat P$ corresponds to the image of the center set 
${\cal C}\subset X$ in $G(X)$. Moreover, $\widehat {\cal C}$ is a
non-trivial, normal, closed, contractible to a point, arcwise connected, 
locally simply and arcwise connected subgroup of $G(X)$.
\end{Th}

Let us consider the quotient group $Q(X):=G(X)/\widehat {\cal C}$.
By $\pi: G(X)\rightarrow Q(X)$ we denote the quotient homomorphism.
According to  Proposition \ref{sepit}, the Taylor coefficients $c_{n}$ of the 
first return map $P$ (see Theorem \ref{center}) satisfy 
$c_{n}(a)=c_{n}(b)$ for $a\sim b$. Therefore they can be
viewed as functions on $G(X)$ (written, $\widehat c_{n}$).
\begin{Proposition}\label{cn}
Every function $\widehat c_{n}$ is constant
on fibres of the map $\pi$ and therefore determines a function
$\overline{c}_{n}:Q(X)\rightarrow\Co$ such that 
$\overline{c}_{n}\circ\pi=\widehat c_{n}$.
\end{Proposition}

Let $\tau''$ be the weakest topology on
$Q(X)$ in which all functions $\overline c_{n}$ are continuous. 
\begin{Th}\label{comp3}
$Q(X)$ equipped with $\tau''$ is a topological group. Moreover,
the homomorphism $\widehat P$ determines an isomorphism 
$\overline{P}:Q(X)\rightarrow G_{c}[[r]]$ defined by 
$\widehat P=\overline{P}\circ\pi$.
\end{Th}

Now, given $f(r)=r+\sum_{k=1}^{\infty}d_{k}r^{k+1}\in G_{c}[[r]]$ we
define $a(f)=(a_{1},a_{2},\dots)\in X$ from 
the identity of formal power series
$$
\sum_{i=1}^{\infty}a_{i}(x)t^{i+1}=
\frac{\sum_{k=1}^{\infty}
(d_{k}/2\pi)t^{k+1}}
{1+\sum_{k=1}^{\infty}(k+1)d_{k}(1-x/2\pi)t^{k}}\ \ \ 
{\rm for}\ \ \ x\in (0,2\pi]
$$
and further extended by periodicity. Let $[a(f)]\in G(X)$ be the image of 
$a(f)$. By $T$ we denote the map $f\mapsto [a(f)]$. Then we have
\begin{Th}\label{decomp}
The map $T:G_{c}[[r]]\rightarrow G(X)$ is a continuous embedding such that
$\widehat P\circ T=id$. Moreover, the map 
$\widetilde T: G_{c}[[r]]\times\widehat {\cal C}\rightarrow G(X)$ defined
by $\widetilde T(f,g):=T(f)\cdot g$ is a homeomorphism.
\end{Th}
{\bf 3.3.} In this part we study some naturally defined subgroups
of $G(X)$.

Let us consider the subspace $X_{s}\subset X$ consisting of
$a=(a_{1},a_{2},\dots)\in X$ such that every $a_{i}$ is a piecewise 
continuous complex-valued function on $S^{1}$. Clearly, for any 
$a,b\in X_{s}$ we have $a*b\in X_{s}$ and $a^{-1}\in X_{s}$. 
Let $G_{s}(X)$ be the image
of $X_{s}$ in $G(X)$. Then $G_{s}(X)$ is a subgroup of $G(X)$. Similarly,
let $X_{a}\subset X$ be the subspace of elements $a=(a_{1},a_{2},\dots)\in X$
such that every $a_{i}$ is a piecewise analytic complex-valued
function on $S^{1}$. Then the image $G_{a}(X)$ in $G(X)$ is a subgroup, 
as well. Also, $G_{a}(X)\subset G_{s}(X)$. In what follows we consider
these groups with relative topologies induced by $\tau$. We will prove 
\begin{Proposition}\label{subgroup1}
$G_{s}(X)$ and $G_{a}(X)$ satisfy properties similar to (1),(2),(4),(5) 
of Theorem \ref{group}.
\end{Proposition}

Let $\widehat {\cal C}_{s}:=\widehat {\cal C}\cap G_{s}(X)$ and 
$\widehat {\cal C}_{a}:=\widehat {\cal C}\cap G_{a}(X)$. 
Then $\widehat {\cal C}_{s}$ and 
$\widehat {\cal C}_{a}$ 
are closed normal subgroups of the corresponding groups. 
\begin{Th}\label{subgroup2}
\begin{itemize}
\item[{\rm (1)}]
$\widehat P|_{G_{s}(X)}:G_{s}(X)\rightarrow G_{c}[[r]]$
and $\widehat P|_{G_{a}(X)}:G_{a}(X)\rightarrow G_{c}[[r]]$ are surjective
homomorphisms of the topological groups whose kernels are 
$\widehat {\cal C}_{s}$ and $\widehat {\cal C}_{a}$, respectively.
\item[{\rm (2)}]
The quotient groups $G_{s}(X)/\widehat {\cal C}_{s}$ and 
$G_{a}(X)/\widehat {\cal C}_{a}$
are isomorphic to $Q(X)$.
\item[{\rm (3)}]
The range of the map $T$ from Theorem \ref{decomp} is a subset of $G_{a}(X)$,
thus, maps $\widetilde T_{s}:G_{c}[[r]]\times\widehat {\cal C}_{s}
\rightarrow G_{s}(X)$, $\widetilde T_{s}(f,g):=T(f)\cdot g$, and 
$\widetilde T_{a}:G_{c}[[r]]\times\widehat {\cal C}_{a}\rightarrow G_{a}(X)$, 
$\widetilde T_{a}(f,g):=T(f)\cdot g$, are homeomorphisms.
\end{itemize}
\end{Th}

We say that a subset $S\subset G(X)$ is {\em defined over a field} 
$\F$, $\Q\subseteq\F\subseteq\Co$, if $\widehat I(s)\in\F$ for every 
$s\in S$ and any basic iterated integral $\widehat I$ on $G(X)$.

Let $X_{\Re}\subset X$ be the subspace of elements
$a=(a_{1},a_{2},\dots)\in X$ such that every $a_{i}$ is a real-valued
function.  Then the corresponding subgroup
$G_{\Re}(X)\subset G(X)$ is defined over $\Re$. Also, we consider
subspaces $X_{s,\Re}:=X_{\Re}\cap X_{s}$ and $X_{a,\Re}:=X_{\Re}\cap X_{a}$
of real-valued piecewise continuous and piecewise analytic sequences of
functions, respectively. By $G_{s,\Re}(X)$ and $G_{a,\Re}(X)$
we denote the subgroups of $G_{\Re}(X)$ corresponding to $X_{s,\Re}$ and
$X_{a,\Re}$. By definition, they are defined over $\Re$.
In what follows we consider these groups with topologies induced by
$\tau$. Let $G_{c,\Re}[[r]]\subset G_{c}[[r]]$ be the subgroup of convergent 
power series with real Taylor coefficients. We consider $G_{c,\Re}[[r]]$
with the topology induced by $\tau'$. Also, we set 
$\widehat {\cal C}_{\Re}:=\widehat {\cal C}\cap G_{\Re}(X)$,
$\widehat {\cal C}_{s,\Re}:=\widehat {\cal C}_{s}\cap G_{\Re}(X)$, 
$\widehat {\cal C}_{a, \Re}:=\widehat {\cal C}_{a}\cap G_{\Re}(X)$, 
$Q_{\Re}(X):=G_{\Re}(X)/\widehat {\cal C}_{\Re}$.
All the above groups defined over $\Re$ will be called the
{\em real analogs} of the corresponding groups defined over $\Co$.
\begin{Th}\label{subgroup3}
The results similar to those of Proposition \ref{subgroup1} and Theorem 
\ref{subgroup2} are valid for the real analogs of the corresponding groups.
\end{Th}

Let us consider another construction of subgroups of $G(X)$ defined
over various number fields that reveals certain arithmetic properties
of the center set.

Given a sequence $b:=\{b_{k}\}_{k\geq 1}\subset 
L^{\infty}(S^{1})$ such that $\sup_{k,x\in S^{1}}|b_{k}(x)|\leq 1$,
let $\hat b=(\hat b_{1},\hat b_{2},\dots)\in X$ be the set of 
all possible monomials $b_{i_{1}}^{\alpha_{1}}\cdots b_{i_{n}}^{\alpha_{n}}$,
$\alpha_{1},\dots,\alpha_{n}\in\Z_{+}$, $i_{1},\dots,i_{n},n\in \N$, 
arranged in a sequence. Suppose $[\hat b]\in G(X)$ is defined
over a field $\F$. Given a number field $\F'\supseteq\F$,
let ${\cal A}_{\F'}(b)$ be the algebra over $\F'$ generated
by $b$, that is $f\in {\cal A}_{\F'}(b)$ if there exist 
$\hat b_{i_{1}},\dots, \hat b_{i_{k}}$ such
that $f=l(\hat b_{i_{1}},\dots,\hat b_{i_{k}})$, where
$l$ is a linear function in $\hat b_{i_{1}},\dots, \hat b_{i_{k}}$ 
with coefficients from $\F'$. Let $X_{\F'}(b)\subset X$ be the subset of
sequences $a=(a_{1},a_{2},\dots)$ such that every 
$a_{i}\in {\cal A}_{\F'}(b)$. By $G_{\F'}(b)\subset G(X)$ we denote
the subgroup generated by images of all possible elements from 
$X_{\F'}(b)$.
\begin{Proposition}\label{number1}
The subgroup $G_{\F'}(b)$ is defined over the field $\F'$.
\end{Proposition}
 
Let $\widehat {\cal C}_{\F'}(b):=\widehat {\cal C}\cap G_{\F'}(b)$.
Suppose $\sigma:\F'\rightarrow\F'$ is a field automorphism over $\F$.
\begin{Proposition}\label{number2}
The automorphism $\sigma$ generates a group automorphism \penalty-10000
$\widetilde\sigma: G_{\F'}(b)\rightarrow G_{\F'}(b)$ continuous in
the discrete topology of $G_{\F'}(b)$ such that \penalty-10000
$\widetilde\sigma(\widehat {\cal C}_{\F'}(b))=\widehat {\cal C}_{\F'}(b)$.
\end{Proposition}

We say that $P:[0,1]\rightarrow X$ is an {\em $X_{\F'}(b)$-polynomial map}
if $P(t)\in X$ for any $t\in [0,1]$, and 
$P(t)=c_{1}^{\epsilon_{1}}(t)*\cdots*c_{k}^{\epsilon_{k}}(t)$ where every
$c_{i}(t)=(c_{1i}(t),c_{2i}(t),\dots)$ with $c_{ki}(t)$ polynomials
in $t$ whose coefficients are elements from 
${\cal A}_{\F'}(b)$ and every $\epsilon_{l}$, $1\leq l\leq k$, equals
either $-1$ or $1$. The induced map $[P]:[0,1]\rightarrow G(X)$, 
$[P](t):=[P(t)]$, will be called {\em an $X_{\F'}(b)$-polynomial arc}. 

Clearly, $[P](t)\in G_{\F'}(b)$ for any 
$t\in\F'\cap [0,1]$. In the sequel we will show that 
$[P]:[0,1]\rightarrow G(X)$ is a continuous map (for 
$G_{\F'}(b)$ considered in the relative topology induced by $\tau$).
\begin{Th}\label{number3}
Suppose $[P](0)\not\in\widehat {\cal C}$ and $[P](\xi)\in\widehat {\cal C}$
for some $\xi\in (0,1]$. Then $\xi$ is algebraic over $\F'$. 
\end{Th}
\begin{E}\label{aritm}
{\rm (1) Suppose $b=\{b_{k}\}_{k\geq 1}\subset L^{\infty}(S^{1})$ where 
$b_{2k-1}(x)=\sin kx$
and $b_{2k}(x)=\cos kx$, $k=1,2,\dots$. Then 
$[\hat b]\in G(X)$ is defined over the field $\Q(\pi)$.\\
(2) Suppose $b=\{b_{k}\}_{k\geq 1}\subset L^{\infty}(S^{1})$ where every
$b_{k}(x):=\chi_{[0,1]}(x)\cdot x^{k-1}$ for $x\in (0,2\pi]$ and then 
extended by periodicity. Here $\chi_{[0,1]}$ is the characteristic
function of $[0,1]$. Then $[\hat b]\in G(X)$ is defined over $\Q$.}
\end{E}
\begin{R}\label{finsubgroup}
{\rm (1) Let $X_{n}$ be the set of elements 
$a=(a_{1},\dots,a_{n},0,0,\dots)\in X$. The image 
of $X_{n}$ in $G(X)$ forms a subgroup $G(X_{n})$ with
properties similar to those of Theorem \ref{group}.
The projection $p_{n}:X\rightarrow X_{n}$,
$p_{n}((a_{1},a_{2},\dots)):=(a_{1},\dots,a_{n},0,0,\dots)$, induces a
surjective homomorphism $p_{n *}:G(X)\rightarrow G(X_{n})$. In particular,
$G(X_{n})$ is isomorphic (in the category of topological groups) to the
semidirect product of $G(X_{n})$ and the normal subgroup 
$Ker(p_{n *})\subset G(X)$.
Set $\widehat {\cal C}_{n}:=\widehat {\cal C}\cap G(X_{n})$. 
From the results of Blinov [Bl] (see Remark \ref{che}) it follows, e.g., that
$\widehat {\cal C}_{n}$ is a non-trivial subgroup of $G(X_{n})$. However,
it is a very delicate and difficult problem to describe explicitly
the elements of $\widehat {\cal C}_{n}$. \\
(2) As before, let $X_{i}$ be the space of all coefficients $a_{i}$ 
from (\ref{e2}). We set $X_{u}:=\prod_{i\geq 1}X_{i}$. Then clearly,
$X\subset X_{u}$. The vector space $X_{u}$ parameterizes
all equations of type (\ref{e2}) with right-hand sides being formal power
series in $v$. Also, one can determine a multiplication $*$ on $X_{u}$
similarly to that of on $X$. Considering an equivalence relation on $X_{u}$
similar to that of on $X$, as a quotient set we get a topological group 
$G(X_{u})$ such that $G(X)\subset G(X_{u})$. One can show that $G(X_{u})$
shares most of the properties similar to those of the group $G(X)$ (see
Sections 3.1 and 3.2). In particular, the center set 
$\widehat {\cal C}_{u}\subset G(X_{u})$ is the kernel of the natural
homomorphism $\widehat P_{u}:G(X_{u})\rightarrow G[[r]]$ corresponding
to formal solutions of (\ref{e2}). Here $G[[r]]$ is the group of
formal complex power series of the form $x+\sum_{i=1}^{\infty}d_{i}r^{i+1}$.
Sometimes to get more information about equation (\ref{e2}) it is
easier to work with $G(X_{u})$ instead of $G(X)$, in particular, because
we don't require any special estimates there.}
\end{R}

The results of Section 3 will be proved in Sections 9-11.
\sect{\hspace*{-1em}. Proofs of Results of Section 2.1.}
In this section we prove Theorems \ref{te1}, \ref{center},
\ref{param} and Corollary \ref{dimens}.\\
{\bf 4.1.} We use the following definitions (see e.g. [Ru, Ch. 7] for
possible references).
\begin{D}\label{lip}
Let $I\subset\Re$ be a closed interval.
A map $f:I\rightarrow\Re^{n}$ is called Lipschitz if 
$dist(f(x),f(y))\leq C|x-y|$ for any $x,y\in I$.
Here $dist(\cdot,\cdot)$ denotes the
Euclidean distance on $\Re^{n}$. A map $f:J\rightarrow\Re^{n}$
defined on an open interval $J\subset\Re$ is called locally Lipschitz if
it is Lipschitz on each closed subinterval of $J$.

Let $mes(V)$ denote the Lebesgue measure of $V\subset\Re$. We say that
$f:I\rightarrow\Re$ is absolutely continuous if $mes(f(E))=0$ 
provided that $mes(E)=0$. 
\end{D}

Any Lipschitz $f:I\rightarrow\Re$ is absolutely continuous. Also, such $f$
is differentiable almost everywhere (Lebesgue) and its derivative
belongs to $L^{\infty}(I)$. The first integral of a function from 
$L^{\infty}(I)$
is a Lipschitz function. Now, if $f$ is absolutely continuous and
$f'$ is zero almost everywhere, then $f$ is constant. Therefore derivative
of antiderivative of any function from $L^{\infty}(I)$ almost everywhere
coincides with this function. Moreover, any Lipschitz $f$ is 
an antiderivative of $f'$. Therefore the Newton-Leibnitz formula
$\int_{a}^{b}f'(t)dt=f(b)-f(a)$ holds. Similar results are valid for
complex-valued Lipschitz functions. These facts are basic when we solve
systems of linear differential equations with 
$L^{\infty}$-coefficients by Picard iteration.\\
{\bf 4.2.} Let us consider equation (\ref{e2}). To make sure that 
it has a Lipschitz solution defined on $[0,2\pi]$ we assume
that there is a constant $l<\infty$, such that for any $i$
\begin{equation}\label{est1}
\sup_{x\in [0,2\pi]}|a_{i}(x)|\leq l^{i}\ .
\end{equation}

First observe that the substitution 
$r(x)=t\cdot v(x)$ transforms (\ref{e2}) to
\begin{equation}\label{e21}
v'=\sum_{i=1}^{\infty}a_{i}(x)t^{i}v^{i+1}
\end{equation}
and that equation (\ref{e2}) determines a center if and only if 
(\ref{e21}) determines a center for any sufficiently small $t$.
Multiplying (\ref{e21}) by $v^{k-1}$ we obtain
\begin{equation}\label{e22}
(v^{k})'=\sum_{i=1}^{\infty}(ka_{i}(x)t^{i})v^{i+k}
\end{equation}
Let $V$ be the linear space spanned by vectors
$e_{i}=(0,\dots,0,1,0,\dots)$ with $1$ at the $i^{{\rm th}}$ place.
We set $y_{i}=v^{i}e_{i}$ and $Y=(y_{1},y_{2},\dots)$. Combining
equations (\ref{e22}) for all $k$ we obtain a formal
linear system
\begin{equation}\label{e23}
Y'=(\sum_{i=1}^{\infty}A_{i}a_{i}(x)t^{i})Y
\end{equation}
where $A_{i}:V\rightarrow V$ are linear operators. Let $\Co[[z]]$ be the 
algebra of formal power series with complex coefficients
$$
f(z)=\sum_{k=0}^{\infty}c_{k}z^{k}\ .
$$
We identify $V$ with $\Co[[z]]$ so that $e_{n}$ coincides with $z^{n-1}$. 
Let $D:\Co[[z]]\rightarrow\Co[[z]]$ be the differentiation operator
$$
(Df)(z):=\sum_{k=0}^{\infty}(k+1)c_{k+1}z^{k}
$$
and $L:\Co[[z]]\rightarrow\Co[[z]]$ be the left translation operator
$$
(Lf)(z):=\sum_{k=0}^{\infty}c_{k+1}z^{k}\ .
$$
Then we have
\begin{Proposition}\label{pr1}
$$
A_{i}=DL^{i-1}\ .
$$
Thus (\ref{e23}) can be written as
\begin{equation}\label{e24}
Y'=(\sum_{i=1}^{\infty}a_{i}(x)t^{i}DL^{i-1})Y\ .
\end{equation}
\end{Proposition}
{\bf Proof.} Applying $DL^{i-1}$ to the function $z^{k}$ we have 0 if 
$k\leq i-1$ and otherwise
$$
DL^{i-1}(z^{k})=D(z^{k-i+1})=(k-i+1)z^{k-i}\ .
$$
Now, by definition, $A_{i}=(a_{k,l}^{i})$ is an infinite matrix such that
$a_{s,i+s}^{i}=s$ for any $s\in\Z_{+}$ and $a_{k,l}^{i}=0$ otherwise.
Thus we have 
$$
A_{i}(e_{k+1})=\sum a_{j,k+1}^{i}e_{j}=(k-i+1)e_{k-i+1}.\ \ \ \ \ \Box
$$

Further, by ${\cal A}(D,L)[[t]]$ we denote the associative algebra 
of formal power series
$$
R(t)=\sum_{k=0}^{\infty}p_{k}(D,L,I)t^{k}
$$
where $p_{i}(D,L,I)$ are polynomials with complex coefficients in
variables $D$, $L$ and $I$, where  
$I:\Co[[z]]\rightarrow\Co[[z]]$ is the identity operator. Let $G(D,L)[[t]]$ be
the group of invertible elements of ${\cal A}(D,L)[[t]]$, and
$\rho':\Z\rightarrow G(D,L)[[t]]$ be the monodromy of (\ref{e24}).
We prove that 
\begin{Lm}\label{ident}
There is an isomorphism $\Phi:{\cal A}(D,L)[[t]]\rightarrow
{\cal A}(S_{1},S_{2})[[t]]$ uniquely defined by $\Phi(D)=S_{1}$ and
$\Phi(L)=S_{2}$. In particular, for the monodromy $\rho$ from 
Theorem \ref{te1} we have $\rho=\Phi\circ\rho'$. 
\end{Lm}
{\bf Proof.} Let ${\cal A}(D,L)$ be the associative algebra of 
polynomials with complex coefficients in variables $D,L$ and $I$.
By ${\cal A}_{0}(D,L)$ we
denote the subalgebra of ${\cal A}(D,L)$ generated by $D$ and $L$ only. Let
$\Co(I)\subset {\cal A}(D,L)$ be the one-dimensional central subalgebra 
generated by $I$. Then ${\cal A}(D,L)={\cal A}_{0}(D,L)\oplus\Co(I)$.
First, we prove the following basic relation for $D$ and $L$.
\begin{Proposition}\label{pr3}
$$
[D,L]:=DL-LD=-L^{2}\ .
$$
\end{Proposition}
{\bf Proof.} It suffices to check the identity for elements
$z^{k}$ with $k\geq 2$.  Then we have
$$
DL(z^{k})=D(z^{k-1})=(k-1)z^{k-2},\ \ \ LD(z^{k})=L(kz^{k-1})=kz^{k-2}.
$$
Therefore
$$
[D,L](z^{k})=(k-1)z^{k-2}-kz^{k-2}=-z^{k-2}=-L^{2}(z^{k}).\ \ \ \ \ \Box
$$

Let ${\cal A}(X_{1},X_{2})$ be the algebra of polynomials with complex 
coefficients in free non-commutative variables $X_{1},X_{2}$. Then there 
exists a
canonical homomorphism $\phi:{\cal A}(X_{1},X_{2})\rightarrow 
{\cal A}_{0}(D,L)$
uniquely defined by the relations $\phi(X_{1})=D$, $\phi(X_{2})=L$.
From Proposition \ref{pr3} it follows that 
$[X_{1},X_{2}]+X_{2}^{2}\in Ker(\phi)$. 
Let $J\subset {\cal A}(X_{1},X_{2})$
be the two-sided ideal generated by $[X_{1},X_{2}]+X_{2}^{2}$. Let
${\cal A}(S_{1},S_{2}):={\cal A}(X_{1},X_{2})/J$ be the quotient algebra and
$\gamma:{\cal A}(X_{1},X_{2})\rightarrow {\cal A}(S_{1},S_{2})$ be the 
quotient homomorphism uniquely defined by $\gamma(X_{1})=S_{1}$ and
$\gamma(X_{2})=S_{2}$. Clearly, there is a homomorphism $\widetilde\phi:
{\cal A}(S_{1},S_{2})\rightarrow {\cal A}_{0}(D,L)$
such that $\phi=\widetilde\phi\circ\gamma$.
\begin{Proposition}\label{pr4}
$\widetilde\phi:{\cal A}(S_{1},S_{2})\rightarrow {\cal A}_{0}(D,L)$ is an 
isomorphism.
\end{Proposition}
{\bf Proof.} Obviously,
any polynomial in $S_{1}$ and $S_{2}$ can be 
transformed with the help of the basic relation to a polynomial in the 
canonical form: 
$P(S_{1},S_{2})=
\sum_{1\leq i,j\leq n}a_{ij}S_{1}^{i}S_{2}^{j}$,
$a_{ij}\in\Co$. Suppose, to the contrary, that $Ker(\widetilde\phi)\neq 0$.
This means that there is a $P(S_{1},S_{2})$ whose coefficients
in the canonical form not all are zeros such that $P(D,L)=0$.
Let us write $P(S_{1},S_{2})=
\sum_{k=1}^{n}F_{k}(S_{1},S_{2})$, where
$F_{k}(S_{1},S_{2})=\sum_{i=0}^{k}a_{ik-i}S_{1}^{i}S_{2}^{k-i}$ is the
homogeneous component of degree $k$. Using the identity
$$
D^{i}L^{j}(z^{s})=D^{i}(z^{s-j})=(s-j)(s-j-1)\cdots (s-j-i+1)z^{s-i-j}
$$
we obtain for $s\geq n$
$$
0=P(D,L)(z^{s})=\sum_{k=1}^{n}[a_{0k}+\sum_{j=0}^{k-1}
a_{k-j j}(s-j)\cdots (s-k+1)]z^{s-k}\ .
$$
In particular, for any $1\leq k\leq n$ and any sufficiently big positive
integers $s$ we have
$$
a_{0k}+\sum_{j=0}^{k-1}a_{k-j j}(s-j)\cdots (s-k+1)=0\ .
$$
But then the same identity is valid for any $s\in\Re$.  Substituting in this
identity $s=k-1$ we get $a_{0k}=0$. Then 
dividing this identity by $(s-k+1)$ and substituting $s=k-2$ we obtain
$a_{1k-1}=0$. Proceeding similarly we finally get $a_{ij}=0$ for all 
$i,j$. This means that $P(S_{1},S_{2})=0$, the
contradiction which proves the proposition.
\ \ \ \ \ $\Box$\\

Now, in order to finish the proof of the lemma it remains to determine 
$$
\Phi(\sum_{k=0}^{\infty}p_{k}(D,L,I)t^{k})=
\sum_{k=0}^{\infty}p_{k}(S_{1},S_{2},I)t^{k},\ \ \ \
\sum_{k=0}^{\infty}p_{k}(D,L,I)t^{k}\in {\cal A}(D,L)[[t]]\ .
$$
Then Proposition \ref{pr4} implies that $\Phi$ is an algebraic isomorphism.
The factorization $\rho=\Phi\circ\rho'$ then easily follows from 
Picard's iteration formula for solutions of (\ref{e24}) and (\ref{formal}).

The proof of Lemma \ref{ident} is complete.\ \ \ \ \ $\Box$

From this lemma it follows that it suffices to prove Theorem \ref{te1} for
equation (\ref{e24}) and its monodromy $\rho'$.\\
{\bf 4.3.} Let $\omega(x):=\sum_{i=1}^{\infty}a_{i}(x)t^{i}DL^{i-1}$. 
The fundamental solution $Y$ of (\ref{e24}) can be obtained by the 
Picard method of successive approximations
\begin{equation}\label{e25}
Y(x)=I+\sum_{n=1}^{\infty}
\int\cdots\int_{0\leq s_{1}\leq\cdots\leq s_{n}\leq x}
\omega(s_{n})\cdots\omega(s_{1})\ ds_{n}\cdots ds_{1}\ .
\end{equation}
Let us define the functions
\begin{equation}\label{itform}
I_{i_{1},\dots,i_{k}}(x):=
\int\cdots\int_{0\leq s_{1}\leq\cdots\leq s_{k}\leq x}
a_{i_{k}}(s_{k})\cdots a_{i_{1}}(s_{1})\ ds_{k}\cdots ds_{1}\ .
\end{equation}
Then it is easy to see that
\begin{equation}\label{e26}
Y(x)=\sum_{i=0}^{\infty}p_{i}(x;D,L)t^{i}
\end{equation}
where $p_{0}=I$ and 
$p_{i}(x;D,L)=\sum_{i_{1}+\dots + i_{k}=i}I_{i_{1},\dots,i_{k}}(x)
(DL^{i_{k}-1})\cdots (DL^{i_{1}-1})$
is a homogeneous polynomial of degree $i$ in $D$, $L$ whose 
coefficients are locally Lipschitz functions in $x\in\Re$.

Now, the solution of (\ref{e24}) with initial value 
$Y(0)=\sum_{i=0}^{\infty}y_{i}z^{i}\in\Co[[z]]$ is defined by
$$
Y(x;z,t,Y(0)):=Y(x)\cdot Y(0)\ .
$$
From (\ref{e26}) it follows that
$$
Y(x;z,t,Y(0))=\sum_{i=0}^{\infty}y_{i}(x;z)t^{i}
$$
where each $y_{i}(x;z)$ is a formal power series in $z$ whose
coefficients are locally \penalty-10000
Lipschitz functions in $x\in\Re$.  Also, it follows from the results
of Section 2.1 that the solution of (\ref{e24}) of this
form with the initial value $Y(0)$ is unique.

Further we study the question when $Y(x;z,t,Y(0))$ is 
a continuous function in $x,z$ and $t$ for $|z|,|t|$ small enough.
\begin{Th}\label{te2}
Let $Y(0)=\sum_{i=0}^{\infty}y_{i}z^{i}$ be such that $|y_{i}|<r^{i+1}$
and $r<\frac{e^{-2\pi}}{2\tilde l}$ where $\tilde l:=\max\{1,l\}$.
Then the solution $Y(x;z,t,Y(0))$ of (\ref{e24}) 
is a function continuous in $x\in [0,2\pi]$ and holomorphic in 
$(t,z)\in\Di_{1/2}\times\Di_{1/2}$. Here $\Di_{R}:=\{z\in\Co\ :\ |z|<R\}$.
\end{Th}
{\bf Proof.} First, let us introduce some notation. 
Let $\Re[[z]]$ be the algebra of formal power series with real coefficients, 
and $C_{+}\subset \Re[[z]$ be the cone of series with non-negative 
coefficients.
For $y_{1},y_{2}\in\Re[[z]]$ we write $y_{1}\leq y_{2}$ if 
$y_{2}-y_{1}\in C_{+}$. Further, let 
${\cal A}_{\Re}(D,L)[[t]]\subset {\cal A}(D,L)[[t]]$ be the 
subalgebra of elements $R(t)=\sum_{k=0}^{\infty}p_{k}(D,L,I)t^{k}$ 
with $p_{k}(x,y,z)$ real polynomials. 
Let $R_{i}(t)=\sum_{k=0}^{\infty}p_{ik}(D,L,I)t^{k}
\in {\cal A}_{\Re}(D,L)[[t]]$, $i=1,2$. We write $R_{1}\leq R_{2}$ if 
$p_{1i}(D,L,I)(f)\leq p_{2i}(D,L,I)(f)$ for any $f\in C_{+}$. Finally,
for any $R(t)=\sum_{k=0}^{\infty}p_{k}(D,L,I)t^{k}\in {\cal A}(D,L)[[t]]$ we
define $|R|(t)=\sum_{k=0}^{\infty}\tilde p_{k}(D,L,I)t^{k}\in 
A_{\Re}(D,L)[[t]]$ such that each coefficient of $\tilde p_{k}$ is the 
modulus of the corresponding coefficient of $p_{k}$. Below we use also
that $D,L$ and $I$ map $C_{+}$ to itself.

Let $\omega(x)$ be the form in (\ref{e24}). Then
$$
|\omega(x)|=\sum_{i=1}^{\infty}|a_{i}(x)|DL^{i-1}t^{i}\ .
$$
Further, according to (\ref{est1}),
$$
|\omega(x)|\leq\widetilde\omega:=
\sum_{i=1}^{\infty}l^{i}\cdot DL^{i-1}t^{i}\ .
$$

Let us consider the equation
\begin{equation}\label{modul}
Y'=\widetilde\omega\cdot Y\ .
\end{equation}
By $\widetilde Y(x;z,t)=\sum_{i=0}^{\infty}\widetilde y_{i}(x;z)t^{i}$
we denote the solution of (\ref{modul}) with the initial value
$\widetilde Y(0)=\sum_{i=0}^{\infty}r^{i+1}z^{i}$. Let 
$Y(x;z,t,Y(0))=\sum_{i=0}^{\infty}y_{i}(x;z)t^{i}$ be the solution
of (\ref{e24}) with the initial value $Y(0)$ satisfying the hypothesis of the
theorem. Suppose $y_{i}(x;z)=\sum_{k=0}^{\infty}w_{ki}(x)z^{k}$ and
$\widetilde y_{i}(x;z)=\sum_{k=0}^{\infty}\widetilde w_{ki}(x)z^{k}$. Then
from the estimate $|\omega(x)|\leq\widetilde\omega$ and from the
formula for the fundamental solution (\ref{e26}) we obtain that
$$
|w_{ki}(x)|\leq\widetilde w_{ki}(x)\ .
$$
Thus in order to prove the theorem it suffices to produce good estimates
of each $\widetilde w_{ki}$.

Notice that (\ref{modul}) can be obtained from the equation
\begin{equation}\label{modul1}
v'=\sum_{i=1}^{\infty}l^{i} t^{i}v^{i+1}
\end{equation}
similarly as we derived equation (\ref{e24}) from (\ref{e21}).
Then if $v(x;t)$ is the solution of (\ref{modul1}) with initial value
$v(0;t)=r$, the solution of (\ref{modul}) with the initial value 
$\widetilde Y(0)=\sum_{i=0}^{\infty}r^{i+1}z^{i}$ is 
$\widetilde Y(x;z,t)=\sum_{i=0}^{\infty}v(x;t)^{i+1}z^{i}$.
Now from the Gronwall-type inequality we obtain that
$$
|v(x;t)|\leq re^{x},\ \ \ x\in [0,2\pi],
$$
provided that $|t|\leq 1$ and $|v(x;t)|\leq\frac{1}{2\tilde l}$.
In particular, 
this is true if $r\leq\frac{e^{-2\pi}}{2\tilde l}$. Also, for this estimate
of $r$ from Picard iteration we
get that $v(x,t)$ is a holomorphic function in $t\in\Di_{1}$ for any
fixed $x\in [0,2\pi]$. Let $v(x;t)^{i+1}=\sum_{k=0}^{\infty}c_{ki}(x)t^{k}$.
Then from the Cauchy estimates for derivatives of $v(x;t)$
we have $|c_{ki}(x)|\leq\frac{1}{(2\tilde l)^{i+1}}$, $x\in [0,2\pi]$. Also,
from above it
follows that $\widetilde y_{k}(x;z)=\sum_{i=0}^{\infty}c_{ki}(x)z^{i}$, 
and so, $c_{ki}(x)=\widetilde w_{ki}(x)\in\Re_{+}$. Finally we obtain
that
$$
|w_{ki}(x)|\leq\widetilde w_{ki}(x)\leq\frac{1}{(2\tilde l)^{i+1}}\ \ \
{\rm for\ any}\ \ \ x\in [0,2\pi]\ .
$$
This immediately implies that for any 
$(t,z)\in\Di_{1/2}\times\Di_{1/2}$
and $x\in [0,2\pi]$ the series determining $Y(x;z,t,Y(0))$ converges 
absolutely and uniformly. Thus $Y(x;z,t,Y(0))$ is continuous in $x$ and 
holomorphic in $(t,z)\in\Di_{1/2}\times\Di_{1/2}$.

The proof of the theorem is complete.\ \ \ \ \ $\Box$\\
{\bf 4.4. Proof of Theorem \ref{te1}.} 
Let $I_{k}\subset\Co[[t]]$ be the ideal generated by
series $f(t)=\sum_{j=k+1}^{\infty}f_{j}t^{j}$. Let
$C_{k}:=\Co[[t]]/I_{k}$ be the quotient algebra and 
$\pi_{k}:\Co[[t]]\rightarrow C_{k}$ the quotient homomorphism.
We identify elements of $C_{k}$ with polynomials
$q(\hat t)=\sum_{i=0}^{k}q_{i}\hat t^{i}$, $\hat t:=\pi_{k}(t)$.
 Let $J_{k}\subset {\cal A}(D,L)[[t]$
be the two-sided ideal generated by elements
$R(t)=\sum_{j=k+1}^{\infty}r_{j}(D,L,I)t^{j}$
where $r_{j}$ are polynomials in $D,L$ and $I$ with complex coefficients. Let
$A_{k}:={\cal A}(D,L)[[t]]/J_{k}$ be the quotient algebra and 
$\rho_{k}:{\cal A}(D,L)[[t]]\rightarrow A_{k}$ the quotient homomorphism. 
Algebra $A_{k}$ is naturally
isomorphic to the algebra ${\cal A}(D,L)\otimes_{\Co}C_{k}$. 
Thus we can identify elements of $A_{k}$ with polynomials
$Q(\hat t)=\sum_{j=0}^{k}q_{j}(D,L,I)\hat t^{j}$.
In particular, for any $R\in A_{k}$ and $f\in\Co[[z]]$ the element 
$R\cdot f$ belongs to $\Co[[z]]\otimes_{\Co}C_{k}$.

For the fundamental solution $Y(x)$ of (\ref{e24}) given by formula
(\ref{e26}) let us determine $Y_{k}(x):=\rho_{k}(Y(x;\hat t,D,L))=
\sum_{i=0}^{k}p_{i}(x;D,L)\hat t^{i}$. Clearly we have
\begin{Lm}\label{reduct}
$Y(2\pi)\equiv I$ if and only if $Y_{k}(2\pi)\equiv I$ for any $k$.
\ \ \ \ \ $\Box$
\end{Lm}

Let 
$$
\omega_{k}(x):=\rho_{k}(\omega(x))=
\sum_{i=1}^{k}a_{i}(x)\hat t^{i}DL^{i-1}
$$
be the image of the form $\omega$ from (\ref{e24}). 
Here $\omega_{k}(x)\in A_{k}$ for any fixed $x$.
Let us consider the equation
\begin{equation}\label{homom}
Y'=\omega_{k}\cdot Y \ .
\end{equation}
From Picard iteration it follows that the fundamental solution of 
(\ref{homom}) is $Y_{k}(x)$. Then 
$Y(x;z,\hat t,Y(0))=Y_{k}(x)\cdot Y(0)=\sum_{i=0}^{k}y_{i}(x;z)\hat t^{i}$ 
is the solution of (\ref{homom}) with initial value $Y(0)\in\Co[[z]]$.
Here $Y(x;z,\hat t,Y(0))\in\Co[[z]]\otimes_{\Co}C_{k}$ 
for any fixed $x$. Let $G(A_{k})$ be the group of invertible
elements of $A_{k}$. Then monodromy $\rho_{k}:\Z\rightarrow G(A_{k})$ of
(\ref{homom}) is defined by $\rho_{k}(n):=Y_{k}(2\pi n)$.
Our purpose is to prove 
\begin{Proposition}\label{prop2}
Suppose that equation (\ref{e21}) determines a center.
Then for any $k$ the monodromy $\rho_{k}$ of (\ref{homom}) is trivial.
\end{Proposition}
{\bf Proof.} 
Let ${\cal P}_{s}\subset\Co[[z]]$ be the space of 
complex polynomials of degree at most $s$. Clearly,
$D,L:{\cal P}_{s}\rightarrow {\cal P}_{s-1}$. Below we denote
$D_{s}=D|_{{\cal P}_{s}}$, $L_{s}=L|_{{\cal P}_{s}}$,
$I_{s}=I|_{{\cal P}_{s}}$, and
$Y_{k,s}(x):=Y_{k}(x)|_{{\cal P}_{s}}$. Then
$$
Y_{k,s}(x)=\sum_{i=0}^{k}p_{i}(x;D_{s},L_{s})\hat t^{i}\ .
$$
Note that $Y_{k,s}(x)$ is the fundamental solution of the equation
\begin{equation}\label{restr}
Y'=(\sum_{i=1}^{k}a_{i}(x)\hat t^{i}D_{s}L_{s}^{i-1})\cdot Y\ .
\end{equation}
The form in this equation is an element of
${\cal A}(D_{s},L_{s})\otimes_{\Co} C_{k}$.
Then the solution of (\ref{restr}) with an initial value 
$Y(0)\in {\cal P}_{s}$ is
$Y_{k,s}(x)\cdot Y(0)=\sum_{i=0}^{k}y_{i}(x;z)\hat t^{i}$. For any fixed $x$
this solution belongs to ${\cal P}_{s}\otimes_{\Co}C_{k}$.

Next, recall that the solution $Y(x;z,t,r)$ of (\ref{e24})
with initial value $Y(0)=\sum_{i=0}^{\infty}r^{i+1}z^{i}$ is given by
$$
Y(x;z,t,r)=Y(x)\cdot Y(0)=\sum_{i=0}^{\infty}y(x;t,r)^{i+1}z^{i}
$$
where $y(x;t,r)$ is the solution of (\ref{e21}) with initial value
$y(0;t,r)=r$.
Moreover, $y(x;t,r)=\sum_{i=1}^{\infty}y_{i}(x)t^{i-1}r^{i}$, and the 
series converges absolutely and uniformly for $x,t\in [0,2\pi]$ provided
that $|r|$ is small enough. In particular,  
$y(x;t,r)$ is Lipschitz in $x\in [0,2\pi]$ and analytic in $t$ and $r$ for
$|t|$, $|r|$ small enough.
If equation (\ref{e21}) determines a center then 
$y(2\pi;t,r)=r$ for any sufficiently small $t$ and $r$. This implies that 
$Y(2\pi;z,t,r)=Y(0;t,z,r)=Y(0)$ as formal power series in $t$, $r$ and $z$.
Let us apply to 
$Y(x;z,t,r)$ the homomorphism $\pi_{k}:\Co[[t]]\rightarrow C_{k}$. Then
we have
$$
\begin{array}{c}
\displaystyle
Y_{k}(x;z,\hat t,r):=\pi_{k}(Y(x;z,t,r))=\pi_{k}(Y(x)\cdot Y(0))=
\\
\displaystyle
\sum_{i=0}^{k}p_{i}(x;D,L)(Y(0))\cdot \hat t^{i}=
\rho_{k}(Y(x))\cdot Y(0)=Y_{k}(x)\cdot Y(0)\ .
\end{array}
$$
Therefore $Y_{k}(x;z,\hat t,r)=
\sum_{i=0}^{\infty}(\pi_{k}(y(x;t,r)))^{i+1}z^{i}$
is the solution of (\ref{homom}) satisfying $Y_{k}(0;z,\hat t,r)=Y(0)$.
From the center condition we also have $Y_{k}(2\pi;z,\hat t,r)=Y(0)$. 
But from the series expansion of $y(x;t,r)$ it follows that
$$
Y_{k}(x;z,\hat t,r)=\sum_{j=1}^{k+1}\widetilde y_{j}(x;\hat t,z)r^{j}
$$
where $\widetilde y_{j}(x;\hat t,z)\in C_{k}$ for any fixed $x,z$.
In particular, $\widetilde y_{j}(2\pi;\hat t,z)=\widetilde y_{j}(0;\hat t,z)$
for any $j$. Let us now apply to $Y_{k}(x;z,\hat t,r)$ the homomorphism
$\pi_{s+1}:\Co[[r]]\rightarrow C_{s+1}$. Then we have
$$
\begin{array}{c}
\displaystyle
Y_{k,s}(x;z,\hat t,\hat r):=\pi_{s+1}(Y_{k}(x;z,\hat t,r))=
\pi_{s+1}(Y_{k}(x)\cdot Y(0))=\\
\displaystyle
Y_{k}(x)\cdot\pi_{s+1}(Y(0))=Y_{k}(x)(\sum_{i=0}^{s}\hat r^{i+1}z^{i})
=Y_{k,s}(x)(\sum_{i=0}^{s}\hat r^{i+1}z^{i})\ .
\end{array}
$$
Here $\hat r:=\pi_{s+1}(r)\in C_{s+1}$.
This means that $Y_{k,s}(x;z,\hat t,\hat r)$ is the solution of equation 
(\ref{restr})
with the initial value $Y_{s}(0):=\sum_{i=0}^{s}\hat r^{i+1}z^{i}$.
Moreover, from the center condition we also have 
$Y_{k,s}(2\pi;z,\hat t,\hat r)=Y_{s}(0)$.

Let us identify the element $z^{i-1}$ with the vector 
$e_{i}=(0,\dots,0,1,0,\dots,0)$ with $1$ at the $i^{\rm th}$ place 
$(1\leq i\leq s+1)$. Then $D_{s},L_{s}$ can be thought of as
nilpotent matrices $A,B$ acting on $\Co^{s+1}$. Now we can write 
(\ref{restr}) as
\begin{equation}\label{restr1}
Y'=(\sum_{i=1}^{k}a_{i}(x)\hat t^{i}AB^{i-1})\cdot Y\ .
\end{equation}
Notice that in this identification $Y_{s}(0)$ coincides with the vector 
$v(\hat r)=\sum_{i=1}^{s+1}\hat r^{i}e_{i}$, the fundamental solution
$Y_{k,s}(x)$ of (\ref{restr}) coincides with the
fundamental solution $Y_{k,s}'(x)$ of (\ref{restr1}),
and $Y_{k,s}(x;z,\hat t,\hat r)$ coincides
with the solution $Y_{k,s}(x;\hat t,v(\hat r))$ of (\ref{restr1}) with
initial value $v(\hat r)$. In particular, from the center condition
we obtain that 
\begin{equation}\label{ident1}
Y_{k,s}(2\pi;\hat t,v(\hat r))=v(\hat r)\ .
\end{equation}
Here we can think of $v(\hat r)$ as a complex-valued function in variable
$\hat r$. Let $r_{1},\dots,r_{s+1}$ be pairwise distinct complex numbers. By 
$V(r_{1},\dots,r_{s+1})$ we denote the
matrix whose columns are $v(r_{1}),\dots,v(r_{s+1})$. Note that
$V(r_{1},\dots,r_{s+1})$ is invertible, because the determinant
of this matrix is the Vandermonde determinant. Now from (\ref{ident1})
we have
$$
V(r_{1},\dots,r_{s+1})=Y_{k,s}'(2\pi)\cdot V(r_{1},\dots,r_{s+1})\ .
$$
This implies that $Y_{k,s}'(2\pi)=Y_{k,s}(2\pi)=I_{s}$. Thus
the monodromy of (\ref{restr}) is trivial.

Now, for any $f(z)=\sum_{i=0}^{\infty}c_{i}z^{i}\in\Co[[z]]$ we write 
$f(z)=f_{s}(z)+r_{s}(z)$ where $f_{s}(z)=\sum_{i=0}^{s}c_{i}z^{i}$.
Let $Y(x;\hat t,f(z))$ be the solution of (\ref{homom}) with initial 
value $f(z)$. Then we have
$$
\begin{array}{c}
\displaystyle
Y(2\pi;\hat t,f(z))-f(z)=(Y_{k}(2\pi)-I)(f(z))=\\
\\
(Y_{k}(2\pi)-I)(f_{s}(z))+(Y_{k}(2\pi)-I)(r_{s}(z))=\\
\\
\displaystyle
(Y_{k,s}(2\pi)-I_{s})(f_{s}(z))+
(Y_{k}(2\pi)-I)(r_{s}(z))=\\
\\
\displaystyle
(Y_{k}(2\pi)-I)(r_{s}(z))=
\sum_{i=1}^{k}p_{i}(2\pi;D,L)(r_{s}(z))\cdot \hat t^{i}\ .
\end{array}
$$
Since each $p_{i}$ is a homogeneous polynomial of degree $i$ in $D$ and
$L$, for $s$ big enough we have
$$
p_{i}(2\pi;D,L)(r_{s}(z))=\sum_{j=s-k}^{\infty}c_{ij}(f)\cdot z^{j}\ .
$$
Here $c_{ij}(f)\in\Co$ depends
on $f$ and $p_{i}$ only. As before, below
$I_{l}\subset\Co[[z]]$ denotes the ideal generated by elements
$f(z)=\sum_{j=l+1}^{\infty}f_{j}z^{j}$. Now the above arguments imply
that 
$$
Y(2\pi;\hat t,f(z))-f(z)\in I_{s-k-1}\otimes_{\Co}C_{k}\ \ \ 
{\rm for\ any}\ s\ .
$$
Finally, using the fact that $\cap_{j=1}^{\infty}I_{j}=0$ we obtain that
$Y(2\pi;\hat t,f(z))=f(z)$ for any $f(z)\in\Co[[z]]$. This implies that
$p_{i}(2\pi;D,L)(f(z))=0$ for any $f(z)\in\Co[[z]]$ and any $1\leq i\leq k$.
Then from the arguments used in the proof of Proposition \ref{pr4} it follows
that $p_{i}(2\pi;D,L)=0$ as an element of ${\cal A}_{0}(D,L)$. 
In particular, $Y_{k}(2\pi)=I$. 

This completes the proof of the proposition.\ \ \ \ \ $\Box$

From Proposition \ref{prop2} and Lemma \ref{reduct} we conclude that
if the original equation (\ref{e21}) determines a center then the monodromy
$\rho'$ of (\ref{e24}) is trivial. 

Conversely, suppose that the monodromy $\rho'$ of
equation (\ref{e24}) is trivial. Let us consider the solution of this
equation with the initial value $Y(0)=\sum_{i=0}^{\infty}r^{i+1}z^{i}$.
Then this solution is 
$Y(x)\cdot Y(0)=\sum_{i=0}^{\infty}y(x;t,r)^{i+1}z^{i}$ where $Y(x)$ is
the fundamental solution of (\ref{e24}) and
$y(x;t,r)$ is the solution of (\ref{e21}) with $y(0;t,r)=r$. Also, for
$|r|<\frac{e^{-2\pi}}{2\tilde l}$ and $|t|<\frac{1}{2}$ from 
Theorem \ref{te2} 
it follows that $y(x;t,r)$ is a continuous function in $x$ analytic in
$t$. Now, from the identity for formal power series $Y(2\pi)=I$ we obtain 
that $y(2\pi;t,r)=r$ for any $r$ and $t$ small enough. 

This completes the proof of the theorem.\ \ \ \ \ $\Box$ \\
{\bf 4.5. Proof of Theorem \ref{center}.} Recall that according to formula
(\ref{e26}) the fundamental solution of (\ref{e24}) is
$Y(x):=\sum_{n=0}^{\infty}p_{n}(x;D,L)t^{n}$,
where $p_{0}=I$ and
$p_{n}(x;D,L)=\sum_{i_{1}+\dots + i_{k}=n}I_{i_{1},\dots, i_{k}}(x)
(DL^{i_{k}-1})\cdots (DL^{i_{1}-1})$.
Now, the solution of (\ref{e24}) with the initial value 
$\sum_{i=0}^{\infty}r^{i+1}z^{i}$ is 
$Y(x)\cdot (\sum_{i=0}^{\infty}r^{i+1}z^{i})=\sum_{i=0}^{\infty}
y(x;t,r)^{i+1}z^{i}$ where $y(x;t,r)$
is the solution of (\ref{e21}) with $y(0;t,r)=r$. Also,
the triviality of the monodromy $\rho$ is equivalent to the equation
$y(2\pi;t,r)-r=0$. This means that in order to find 
equations determining the center set for (\ref{e2}) we should apply each
$p_{n}(2\pi;D,L)$ to $\sum_{i=0}^{\infty}r^{i+1}z^{i}$, then substitute
$z=0$ and equate the latter to 0. Now, it is clear that we obtain the same 
equations if we apply $p_{n}(2\pi;D,L)$ to $z^{n}$. 
From here and (\ref{e26}) we obtain for $i_{1}+\dots +i_{k}=n$,
$$
c_{i_{1},\dots,i_{k}}:=(DL^{i_{k}-1})\cdots (DL^{i_{1}-1})(z^{n})=
(n-i_{1}+1)(n-i_{1}-i_{2}+1)\cdots 1\ .
$$
Also, the expression $P(r):=y(2\pi;1,r)$ is the first return map for
equation (\ref{e2}). From above it follows that
$$
P(r)=r+\sum_{n=1}^{\infty}c_{n}r^{n+1}\ \ \ {\rm where}\ \ \
c_{n}=\sum_{i_{1}+\dots +i_{k}=n}c_{i_{1},\dots,i_{k}}I_{i_{1},\dots,i_{k}}\ .
$$

This proves the required result.\ \ \ \ \ $\Box$\\
{\bf 4.6. Proof of Corollary \ref{dimens}.} The proof of part (a) follows
directly from the above formula for the first return map. Let us prove
part (b).

By $\pi_{k}:X\rightarrow\prod_{j=1}^{k}X_{j}$ we denote the natural projection
to the first $k$ coordinates. According to part (a) the set ${\cal C}_{n}$
is given by equations $I_{k}(a)=-f_{k}(a)$, $1\leq k\leq n$, where 
$I_{k}$ depends on $a_{k}$ and $f_{k}$ depends on $a_{1},\dots,a_{k-1}$. 
Let $\widetilde {\cal C}_{n}$ be the subset of $\prod_{j=1}^{n}X_{j}$
defined by these equations. Then clearly 
${\cal C}_{n}=\pi_{n}^{-1}(\widetilde {\cal C}_{n})$. Thus in order
to prove the corollary it suffices to
prove that $\widetilde {\cal C}_{n}\subset\prod_{j=1}^{n}X_{j}$ is a closed 
complex submanifold of codimension $n$ containing $0$. We prove this
by induction on $n$. 

If $n=1$, then $\widetilde {\cal C}_{1}$ is determined
by the equation $I_{1}(a):=\int_{0}^{2\pi}a_{1}(s)\ ds=0$. Since
$I_{1}$ is a continuous linear functional, 
$\widetilde {\cal C}_{1}$ coincides with the complex hyperplane
$Ker(I_{1})\subset X_{1}$. Thus the
required statement holds. Suppose now that the statement is proved for
$k-1$, let us prove it for $k$. 
It is easy to see that 
$\widetilde {\cal C}_{k}$ is a subset of 
$\widetilde {\cal C}_{k-1}\times X_{k}$ determined by the equation
$I_{k}(a)=-f_{k}(a)$.
Since $I_{k}|_{X_{k}}$ is a continuous linear functional, we 
can decompose $X_{k}=E_{k}\oplus l_{k}$, where $E_{k}=Ker(I_{k}|_{X_{k}})$ and
$l_{k}$ is the one-dimensional vector space generated by a vector
$e_{k}$ such that $I_{k}(e_{k})=1$. Further, let us consider the set 
$R_{k}:=\widetilde {\cal C}_{k-1}\times E_{k}$.
It is clear that $R_{k}\subset Ker(I_{k})$ and 
$\widetilde {\cal C}_{k-1}\times X_{k}=R_{k}\oplus l_{k}$.
In particular, for any $a\in \widetilde {\cal C}_{k-1}\times X_{k}$ we have
$a=w+ve_{k}$ for some $w\in R_{k}$, $v\in\Co$.  Thus in the definition of 
$\widetilde {\cal C}_{k}$ we have
$$
v=I_{k}(a)=-f_{k}(a)=-f_{k}(w)\ .
$$
This shows that $\widetilde {\cal C}_{k}$ is the graph of the function 
$-f_{k}: R_{k}\rightarrow\Co$. But according to the induction hypothesis,
$R_{k}\subset (\prod_{j=1}^{k-1}X_{j})\times E_{k}$ is 
a closed complex submanifold of codimension $(k-1)$ containing 0. From here
it follows that 
$\widetilde {\cal C}_{k}\subset \widetilde {\cal C}_{k-1}\times X_{k}$
is a closed complex submanifold of codimension 1 containing 0. This
implies the required result and completes the proof of part (b).

Let us prove part (c). 

It follows
from the definition that the linear term of the Taylor expansion of $c_{n}$ at
$0\in X$ equals $I_{n}$. It remains to note
that the tangent space to ${\cal C}_{n}$ at $0\in X$ is determined as the set
of common zeros of the linear terms of $c_{1},\dots, c_{n}$, that is,
it is given by equations $I_{1}(a)=\dots=I_{n}(a)=0$.\ \ \ \ \ $\Box$\\
{\bf 4.7. Proof of Theorem \ref{param}.} Let $v(x;r;a)$, $x\in [0,2\pi]$,
be the Lipschitz solution of equation (\ref{e2}) corresponding to 
$a=(a_{1}, a_{2},\dots)\in X$ with initial value $v(0;r;a)=r$. From  the 
Picard method of successive approximations it follows that 
\begin{equation}\label{sol1}
v(x;r;a)=r+\sum_{i=1}^{\infty}v_{i}(x;a)r^{i+1}
\end{equation}
where each $v_{i}(x;a)$ is a  Lipschitz function on $[0,2\pi]$ and the
series converges uniformly in the domain $0\leq x\leq 2\pi$, 
$|r|\leq\widetilde r$ for a sufficiently small positive $\widetilde r$. 
Assuming that $a\in {\cal C}$ we get $v_{i}(0;a)=v_{i}(2\pi;a)=0$ for any $i$.
Now the inverse function theorem implies that there is a function
$$
u(x;\rho;a):=\rho+\sum_{i=1}^{\infty}u_{i}(x;a)\rho^{i+1}
$$
where each $u_{i}(x;a)$ is a $2\pi$-periodic Lipschitz function, 
$u_{i}(0;a)=0$, and the
series converges uniformly in the domain $-\infty<x<\infty$, 
$|\rho|\leq\widetilde\rho$ for a sufficiently small positive $\widetilde\rho$,
such that
$$
u(x;v(x;r;a);a)\equiv r
$$
for all sufficiently small $r$. Differentiating this equation in $x$ and
then expressing $v'(x;r;a)$ we obtain
$$
v'(x;r;a)=-\frac{\sum_{k=1}^{\infty}u_{k}'(x;a)v(x;r;a)^{k+1}}
{1+\sum_{k=1}^{\infty}(k+1)u_{k}(x;a)v(x;r;a)^{k}}\ .
$$
This and the equality
$$
v'(x;r;a)=\sum_{i=1}^{\infty}a_{i}(x)v(x;r;a)^{i+1}
$$
imply the required identity of power series.

Conversely, from the identity of Theorem \ref{param} arguing in the reverse
order we can find a formal solution of (\ref{e2}) corresponding to
$a\in X$ of the form (\ref{sol1}) where $v_{i}$ are polynomials in
$u_{1},\dots, u_{i}, v_{1},\dots, v_{i-1}$ without constant terms. In 
particular, each $v_{i}$ is a $2\pi$-periodic Lipschitz function
equals 0 at 0.
But since $a\in X$, there is also a usual solution of the corresponding
equation (\ref{e2}) given by (\ref{sol1}). Hence these two solutions
must coincide. This implies the identity $v(0;r;a)=v(2\pi;r;a)=r$ for any
sufficiently small $r$.

The proof of the theorem is complete.\ \ \ \ \ $\Box$
{\sect{\hspace*{-1em}. Proofs of Results of Sections 2.2 and 2.3.}
In this section we prove Corollary \ref{univ} and Propositions \ref{integr}
and \ref{integr1}.\\
{\bf Proof of Corollary \ref{univ}.} The proof follows directly from
the factorization $\rho=\phi\circ\widetilde\rho$.\ \ \ \ \ $\Box$\\
{\bf Proof of Proposition \ref{integr}.} We set 
$R_{i}:=t^{i}X_{1}X_{2}^{i-1}$.
Let ${\cal A}$ be the algebra of polynomials in 
$X_{1},X_{2},\dots$ with complex coefficients.
The main point of the proof is the following
\begin{Lm}\label{free}
Suppose $P(R_{1},\dots,R_{k})=\sum_{1\leq k\leq n}
a_{i_{1},\dots i_{k}}R_{i_{1}}\cdots R_{i_{k}}\in
{\cal A}$ is zero. Then all $a_{i_{1},\dots i_{k}}$ are zeros.
\end{Lm}
{\bf Proof.} 
The maximal integer $l=\sum_{s=1}^{k}i_{s}$ for indices in the expansion of
$P$ will be called the degree of $P$. We prove the lemma by induction on $l$.

For $l=1$ we have $P=c_{1}R_{1}$ and the statement is obvious. Suppose that we
already proved the lemma for $l\leq s$, let us prove it for $l=s+1$.
We can write $P=P_{1}+\dots+P_{k}$ where $P_{i}$ contains all terms of $P$ 
with $R_{i}$ on the right. Since $X_{1},X_{2}$ are free 
non-commutative variables, we have $P_{1}=\dots=P_{k}=0$. Now 
$P_{i}=Q_{i}\cdot R_{i}$ where $Q_{i}$ is a polynomial of degree $\leq s$
satisfying the hypothesis of the lemma. Thus by the induction hypothesis
all coefficients of every $P_{i}$ are zeros. \ \ \ \ \ $\Box$

Further, it follows 
from formula (\ref{e3'}) for the fundamental solution $F(x)$ of
equation (\ref{formal1}) that 
$F(2\pi)=\sum_{i=0}^{\infty}q_{i}(2\pi,X_{1},X_{2})$ where $q_{0}=I$ and
others $q_{i}$ are polynomials in $R_{1},R_{2},\dots$ of
degree $i$ (as in Lemma \ref{free}). Moreover, every iterated
integral 
$$
\int\cdots\int_{0\leq s_{1}\leq\cdots\leq s_{k}\leq 2\pi}
a_{i_{k}}(s_{k})\cdots a_{i_{1}}(s_{1})\ ds_{k}\cdots ds_{1},\ \ \
\sum_{s=1}^{k}i_{s}=i\ ,
$$
is a coefficient before one of monomials of $q_{i}$. 
Now, the triviality of the 
monodromy
is equivalent to the fact that all $q_{i}\in {\cal A}$, $i\geq 1$, are zeros.
According to Lemma \ref{free} this, in turn, is
equivalent to the fact that all coefficients of all $q_{i}$, $i\geq 1$,
are zeros.

This completes the proof of the proposition.\ \ \ \ \ $\Box$\\
{\bf Proof of Proposition \ref{integr1}.}
The proof of the proposition repeats word-for-word the proof of Proposition
\ref{integr}. We leave the details to the reader.\ \ \ \ \ $\Box$
{\sect{\hspace*{-1em}. Proof of Theorem 2.9.}
{\bf 6.1.} We refer to the book of Hirzebruch [Hi] for an exposition 
about fibre bundles.

Let $U\subset\Co^{n}$ be a domain containing $\widehat\Gamma_{n}$. 
Since $\widehat\Gamma_{n}$
is polynomially convex, for any domain $V\supset\widehat\Gamma_{n}$ there
is a Stein domain $V'$ such that $\widehat\Gamma_{n}\subset V'\subset V$. 
Thus without loss of
generality we may assume that $U$ is Stein. Let $\pi_{1}(U)$ be the
fundamental group of $U$. The universal covering 
$p:\widetilde U\rightarrow U$ is a discrete bundle on $U$ with the fibre
$\pi_{1}(U)$. It is defined on an  open 
acyclic cover ${\cal U}=(U_{i})$ of $U$ by a locally
constant cocycle $\{c_{ij}\}\in Z^{1}({\cal U},\pi_{1}(U))$. By definition,
$c_{ij}:U_{i}\cap U_{j}\rightarrow\pi_{1}(U)$ is a constant map such that
$c_{ij}c_{jk}=c_{ik}$ on $U_{i}\cap U_{j}\cap U_{k}$. Then $\widetilde U$
is biholomorphic to the quotient space of $\sqcup_{j}U_{j}\times\pi_{1}(U)$
by the equivalence relation: 
$$
U_{j}\times\pi_{1}(U)\ni z\times g\sim z\times g\cdot c_{ij}^{-1}\in 
U_{i}\times\pi_{1}(U)\ .
$$

Let $l_{U}:=l_{2}(\pi_{1}(U))$ be the Hilbert space of complex-valued 
sequences on $\pi_{1}(U)$ with the $l_{2}$-norm. By $L(l_{U})$ we denote
the Banach space of bounded complex linear operators 
$l_{U}\rightarrow l_{U}$, and by $GL(l_{U})$ the group of invertible 
operators from $L(l_{U})$.
Let us define the homomorphism 
$\xi:\pi_{1}(U)\rightarrow Iso(l_{U})$ by
$$
\xi(g)(v)(x):=v(x\cdot g^{-1}),\ \ \ g,x\in\pi_{1}(U),\ v\in 
l_{U}\ .
$$
Here $Iso(l_{U})\subset GL(l_{U})$ is the subgroup of unitary operators. 
It is clear that $\xi$ is a faithful representation, that is,
$Ker(\xi)=\{1\}$.

Next, let us construct the holomorphic Banach vector bundle $E_{\xi}$ on 
$U$ with fibre 
$l_{U}$ associated with $\xi$. It is defined as the quotient of 
$\sqcup_{j}U_{j}\times l_{U}$ by the equivalence relation
$$
U_{j}\times l_{U}\ni x\times w\sim x\times\xi(c_{ij})(w)
\in U_{i}\times l_{U}\ .
$$
Since $l_{U}$ is an infinite dimensional Hilbert space, the group
$GL(l_{U})$ is contractible (see [K]). In particular, $E_{\xi}$ is a
topologically trivial Banach vector bundle. Now according to the result of 
Bungart [B] applied to the Stein manifold $U$ we obtain that
{\em $E_{\xi}$ is a holomorphically trivial Banach vector bundle on $U$}.
This means that there is a family of holomorphic functions
$F_{i}:U_{i}\rightarrow GL(l_{U})$ satisfying
$$
F_{i}^{-1}(z)\cdot F_{j}(z)=\xi(c_{ij})\ \ \ {\rm on}\ \ \ U_{i}\cap U_{j}\ .
$$
Note that cocycle $\{\xi(c_{ij})\}$ is locally constant, i.e.
$d(\xi(c_{ij}))=0$ for any $i,j$. Then we define the global holomorphic
1-form $\omega$ on $X$ with values in $L(l_{U})$ by the formula
$$
\omega|_{U_{i}}:=dF_{i}\cdot F_{i}^{-1}\ .
$$
Clearly $\omega$ satisfies the Frobenius condition 
$d\omega-\omega\wedge\omega=0$, which is equivalent to the fact that equation
\begin{equation}\label{monod}
dF=\omega\cdot F
\end{equation}
is locally solvable on $U$. Let $p^{*}\omega$ be the pullback of $\omega$
to $\widetilde U$. Consider the pullback to $\widetilde U$ of 
(\ref{monod})
\begin{equation}\label{monod1}
dF=p^{*}\omega\cdot F\ .
\end{equation}
Then from simply connectedness of $\widetilde U$ it follows that 
(\ref{monod1}) has a global holomorphic solution
$F:\widetilde U\rightarrow GL(l_{U})$. On the connected component
$U_{i}\times\{g\}$ of the open set 
$p^{-1}(U_{i}):=\sqcup_{g\in\pi_{1}(U)}\ U_{i}\times\{g\}\subset\widetilde U$
the solution $F$ equals $p^{*}F_{i}\cdot\xi(g^{-1})$.

Another way to obtain a solution of (\ref{monod1}) is by using Picard 
iteration.
To this end, we fix a point $z_{0}\in\widetilde U$ and for any piecewise 
smooth path $\gamma$ joining $z_{0}$ with a point $z$ we set
\begin{equation}\label{iterat}
S(z)=I+\sum_{i=1}^{\infty}\int\cdots\int_{\gamma}\omega\cdots\omega\ .
\end{equation}
The $k$-th term of this series is the $k$-iterated integral of
$\omega$ over $\gamma$ defined similarly to (\ref{e3'}). It is well known
(cf. [Na]) that the series converges uniformly on compacts in 
$\widetilde U$ to a $GL(l_{U})$-valued function $S$ satisfying 
(\ref{monod1}), and $S(z_{0})=I$. Therefore 
$S(z)=F(z)\cdot F^{-1}(z_{0})$ where $F$ is the solution of 
(\ref{monod1}) constructed above.
Moreover, by the definition of $F$,
$$
S(gz)=S(z)\cdot\xi(g^{-1})\ \ \ {\rm for\ any}\ \ \ g\in\pi_{1}(U),\ 
z\in\widetilde U\ .
$$
In particular, for $z=z_{0}$ we have $S(gz_{0})=\xi(g^{-1})$.

Let $\gamma_{g}\subset U$ be a loop based at $p(z_{0})$ which represents
an element $g\in\pi_{1}(U,p(z_{0}))$. Then the path on $\widetilde U$
with the origin at $z_{0}$ and the endpoint at $gz_{0}$ represents the lifting
of $\gamma_{g}$ to $\widetilde U$. Now the correspondence
$\gamma_{g}\mapsto\xi(g^{-1})$ determines a homomorphism 
$\widetilde\xi:\pi_{1}(U,p(z_{0}))\rightarrow Aut(GL(l_{U}))$ defined by
$$
\widetilde\xi(g)(H)=H\cdot\xi(g^{-1}),\ \ \ g\in\pi_{1}(U,p(z_{0})),\
H\in GL(l_{U})\ .
$$
By the definition of $\xi$ we obtain that $Ker(\widetilde\xi)=\{1\}$.
\\
{\bf 6.2.} For basic facts of complex analysis in domains of
holomorphy see, e.g. the book of Henkin and Leiterer [HL].

Let us write the form $\omega$ in (\ref{monod}) as 
$$
\omega(z)=\sum_{j=1}^{n}R_{j}(z)\ dz_{j}
$$
where $z_{1},\dots,z_{n}$ are standard coordinates on $\Co^{n}$, and
$R_{1},\dots,R_{n}$ are holomorphic $L(l_{U})$-valued functions on
$U$. Let $V\subset\subset U$ be a domain containing $\widehat\Gamma$.
Then there is a Weil polynomial polyhedron 
$$
D=\{z\in U\ :\ |P_{k}(z)|<1,\ P_{k}\in\Co[z_{1},\dots,z_{n}],\ 
k=1,\dots, N\}\subset V
$$
containing $\widehat\Gamma_{n}$. Note that for $L(l_{U})$-valued holomorphic
functions defined in an open neighbourhood of $V$ we have exactly the
same integral representations as for scalar holomorphic functions. In 
particular, any such function $f$ inside $D$ can be represented as a finite 
sum of holomorphic $L(l_{U})$-valued functions $f_{J}$ such that each
$$
f_{J}(z)=\int_{\sigma_{J}}f(\xi)K_{J}(\xi,z)\eta(\xi)
$$ 
where $\eta(\xi)=d\xi_{1}\wedge\dots d\xi_{n}$, $K_{J}(\xi,z)$ is 
a rational holomorphic function in $z$ and $\xi$, and $\sigma_{J}$ is
an $n$-dimensional part of the Silov boundary of $D$. Moreover, one
can uniformly approximate on compacts in $D$ each $K_{J}(\xi,z)$ by 
holomorphic functions $K_{Ji}(\xi,z)$, $i=1,\dots,$ which are holomorphic 
polynomials in $z$. In particular, we obtain the Oka-Weil approximation
result.
\begin{Proposition}\label{approx}
Any $L(l_{U})$-valued function holomorphic in an open neighbourhood 
$V\subset\subset U$ of $\widehat\Gamma_{n}$ can be uniformly approximated on 
$\widehat\Gamma_{n}$ by $L(l_{U})$-valued holomorphic polynomials.\ \ \ \ \
$\Box$
\end{Proposition}

From this proposition it follows that there is a sequence of 
polynomial 1-forms
$$
\omega_{i}(z)=\sum_{j=1}^{n}R_{ji}(z)\ dz_{j}
$$
such that $R_{ji}(z)=\sum_{0\leq |\alpha|\leq d_{ji}}R_{ji,\alpha}z^{\alpha}$,
$R_{ji,\alpha}\in L(l_{U})$, and 
$$
\lim_{i\to\infty}\sup_{z\in\widehat\Gamma}||R_{j}(z)-R_{ji}(z)||=0\ .
$$
Here $z^{\alpha}:=z^{\alpha_{1}}\cdots z^{\alpha_{n}}$ for any
$\alpha:=(\alpha_{1},\dots,\alpha_{n})\in (\Z_{+})^{n}$, 
$|\alpha|:=\sum_{i=1}^{n}\alpha_{i}$, and $||\cdot||$ denotes the 
Banach norm on $L(l_{U})$.\\
{\bf 6.3.} Next, we recall the Ree formula [R] for the product of iterated 
integrals. 

A permutation $\sigma$ of $\{1,2,\dots,r+s\}$ is a {\em shuffle} of type 
$(r,s)$ if
$$
\sigma^{-1}(r+1)<\sigma^{-1}(2)<\cdots <\sigma^{-1}(r)
$$
and
$$
\sigma^{-1}(r+1)<\sigma^{-1}(r+2)<\cdots <\sigma^{-1}(r+s)\ .
$$

Let  $f_{1},f_{2},\dots,f_{r+s}\in L^{\infty}([0,t])$. As before, we set
$$
I_{i_{1},\dots,i_{k}}(x):=
\int\cdots\int_{0\leq s_{1}\leq\cdots\leq s_{k}\leq x}
f_{i_{k}}(s_{k})\cdots f_{i_{1}}(s_{1})\ ds_{k}\cdots ds_{1}\
$$
Then
$$
I_{1,2,\dots,r}(x)\cdot I_{r+1,r+2,\dots,r+s}(x)=\sum_{\sigma}
I_{\sigma (1),\sigma (2),\dots,\sigma (r+s)}(x)
$$
where $\sigma$ runs over the shuffles of type $(r,s)$.\\

Going back to the proof of the theorem we set 
$\widetilde R_{j}(x):=A_{n}^{*}(R_{j})(x)$, $j=1,\dots, n$.
Then $\widetilde R_{j}$ are $L(l_{U})$-valued continuous functions on $S^{1}$.
\begin{Lm}\label{le2}
For any positive integers $k$ and $i_{1},\dots,i_{k}$,
$1\leq i_{s}\leq n$, we have
\begin{equation}\label{zero}
\int\cdots\int_{0\leq s_{1}\leq\cdots\leq s_{k}\leq 2\pi}
\widetilde R_{i_{1}}(s_{k})a_{i_{1}}(s_{k})\cdots 
\widetilde R_{i_{k}}(s_{1})a_{i_{k}}(s_{1})
\ ds_{k}\cdots ds_{1}=0\ .
\end{equation}
\end{Lm}
{\bf Proof.} According to Proposition \ref{approx} it suffices to prove
the result for $\widetilde R_{j}$ holomorphic $L(l_{U})$-valued polynomials.
Moreover, in this case it suffices to prove the lemma for scalar 
monomials. Thus without loss of generality we may assume that
$$
\widetilde R_{i_{l}}(s_{l})=\widetilde a_{1}(s_{l})^{\alpha_{1, i_{l}}}\cdots
\widetilde a_{n}(s_{l})^{\alpha_{n, i_{l}}}\ .
$$
Here $\widetilde a_{i}$ is the antiderivative of $a_{i}$. 
Now the required result
follows from the fact that by the Ree formula each integral 
on the left-hand-side of (\ref{zero}) can be written as a finite sum of  
integrals of the form
$$
\int\cdots\int_{0\leq s_{1}\leq\cdots\leq s_{r}\leq 2\pi}
a_{j_{r}}(s_{r})\cdots a_{j_{1}}(s_{1})
\ ds_{r}\cdots ds_{1}\ 
$$
where $r=\sum_{k=1}^{n}(\alpha_{k,i_{l}}+1)$, and from equations
(\ref{vanish}). 
The calculation is straightforward and we leave it as an exercise for the 
readers. \ \ \ \ \ $\Box$\\
{\bf 6.4.} Let us finish the proof of the theorem.

Let $A_{n}:S^{1}\rightarrow U$ represent an element $h\in\pi_{1}(U,p(z_{0}))$.
Let $\widetilde A_{n}:S^{1}\rightarrow\widetilde U$ be the lifting of $A_{n}$
with the origin at $z_{0}$. Then its endpoint is $hz_{0}$. 
Let $S(z)$ be
the solution of (\ref{monod1}) from Section 4.1. Then 
$S(hz_{0})=\xi(h^{-1})$. On the other hand, this value can be obtained
by formula (\ref{iterat}) where the integrals are taken over the path
$\widetilde A_{n}:S^{1}\rightarrow\widetilde U$. To calculate these integrals
we take the pullback of the form $p^{*}\omega$ by $\widetilde A_{n}$,
and then calculate usual iterated integrals of 
$A_{n}^{*}\omega:=\widetilde A_{n}^{*}(p^{*}\omega)$ over $[0,2\pi]$. But each
of such iterated integrals is zero by Lemma \ref{le2}. This means
that $\xi(h^{-1})=I$. Because $\xi$ is a faithful representation the
latter implies that $h=1\in\pi_{1}(U,p(z_{0}))$.

The proof of the theorem is complete.\ \ \ \ \ $\Box$
\sect{\hspace*{-1em}. Proofs of Theorem 2.12 and Corollaries 2.10,
2.15 and 2.17.}
{\bf Proof of Corollary \ref{coro1}.} Suppose $\Gamma_{n}$ is 
triangulable,
$\widehat \Gamma_{n}=\Gamma_{n}$, and $\Gamma_{n}$ satisfies the
hypothesis of Theorem \ref{contr}. From the triangulability of $\Gamma_{n}$ it
follows that $\Gamma_{n}$ is homeomorphic to a one-dimensional
simplicial complex (i.e., a finite graph). Moreover, by Borsuk's theorem
(for the references see, e.g. [Hu]) there is an open connected neighbourhood
$U$ of $\Gamma_{n}$ and a retraction $r:U\rightarrow\Gamma_{n}$. 
Let $i:\Gamma_{n}\hookrightarrow U$ be the embedding. By 
$i_{*}:\pi_{1}(\Gamma_{n})\rightarrow\pi_{1}(U)$ and 
$r_{*}:\pi_{1}(U)\rightarrow\pi_{1}(\Gamma_{n})$ we denote the corresponding
induced homomorphisms of fundamental groups. Then 
$r_{*}\circ i_{*}:\pi_{1}(\Gamma_{n})\rightarrow\pi_{1}(\Gamma_{n})$ is the
identity homomorphism.
Now, it follows
 from the condition $\widehat\Gamma_{n}=\Gamma_{n}$ and Theorem 
\ref{contr} that the path $i\circ A_{n}:S^{1}\rightarrow U$ represents 
$1\in\pi_{1}(U)$. Thus from above we obtain that 
$A_{n}:S^{1}\rightarrow\Gamma_{n}$
represents $1\in\pi_{1}(\Gamma_{n})$. This means that $A_{n}$ is 
contractible to a point inside $\Gamma_{n}$.

Further, let $p:\Gamma_{nu}\rightarrow\Gamma_{n}$ be the universal covering of
$\Gamma_{n}$. Then $\Gamma_{nu}$ can be thought of as an infinite tree. The
group $\pi_{1}(\Gamma_{n})$ (which is a free group with a finite number of
generators) acts discretely on $\Gamma_{nu}$.
According to the covering homotopy theorem 
(see, e.g. [Hu, Ch.III]), there is a map $A_{2n}:S^{1}\rightarrow\Gamma_{nu}$
such that $A_{n}=p\circ A_{2n}$. The image 
$G_{n}:=A_{2n}(S^{1})\subset\Gamma_{nu}$ is 
a connected compact subset. Therefore $G_{n}$ is homeomorphic to a finite 
tree. Finally, we set $A_{1n}:=p|_{G_{n}}$. 

Conversely, if $A_{n}$ admits the factorization of
Corollary \ref{coro1} then $A_{2n}:S^{1}\rightarrow G_{n}$ is contractible 
to a point inside $G_{n}$. Therefore 
$A_{n}:S^{1}\rightarrow\Gamma_{n}$ is contractible to
a point inside $\Gamma_{n}$.\ \ \ \ \ $\Box$\\
{\bf Proof of Theorem \ref{pr1'}.} Suppose $\Gamma_{n}$ is
Lipschitz triangulable, $\widehat\Gamma_{n}=\Gamma_{n}$, and 
$A_{n}:S^{1}\rightarrow\Gamma_{n}$ is contractible in $\Gamma_{n}$ to a point.
By the covering homotopy theorem, there is a covering 
$s:\widetilde U\rightarrow U$ and a Lipschitz embedding
$i:\Gamma_{nu}\hookrightarrow\widetilde U$ such that $p=s\circ i$.
Without loss of generality we may consider $U$ as a submanifold of some
$\Re^{N_{n}}$ so that $i(\Gamma_{nu})\subset\Re^{N_{n}}$ has an exhaustion by 
Lipschitz triangulable compact subsets. It follows from the Whitney embedding 
theorem. We identify $\Gamma_{nu}$ with $i(\Gamma_{nu})$. 

Now, the contractibility of $A_{n}$ implies that there is a map 
$A_{2n}:S^{1}\rightarrow\Gamma_{nu}\subset\Re^{N_{n}}$ such that 
$A_{n}=s\circ A_{2n}$.
Since $A_{n}$ is Lipschitz, the map $A_{2n}$ is also Lipschitz by the 
definition
of $\widetilde U$. Let us consider the equation on $U$
\begin{equation}\label{lift}
dF=\eta\cdot F,\ \ \ {\rm where}\ \ \ 
\eta:=\sum_{i=1}^{n}t^{i}X_{1}X_{2}^{i-1}\ dz_{i}\ .
\end{equation}
This equation is not locally solvable. 
However, its lifting $dF=s^{*}\eta\cdot F$ 
restricted to $\Gamma_{nu}$ is locally solvable. This follows from the 
fact that $\Gamma_{nu}$ outside a countable discrete set
of points is the disjoint union of one-dimensional Lipschitz manifolds
$M_{i}$, $i=1,2,\dots$ 
(the consequence of Lipschitz triangulability). Moreover, each $M_{i}$
is the image of the interval $(0,1)$ under some Lipschitz map
$g_{i}:(0,1)\rightarrow\Re^{N_{n}}$ so that 
$g_{i}^{-1}:M_{i}\rightarrow (0,1)$
is locally Lipschitz. Now because $\Gamma_{nu}$ is a tree, we can solve the 
equation $dF=s^{*}\eta\cdot F$ on $\Gamma_{nu}$ by Picard
iteration (as in (\ref{iterat})). Indeed, since $g_{i}$ is Lipschitz, we
can solve the pullback of this equation by $g_{i}$ to $(0,1)$. Then
we transfer this solution to $M_{i}$ by $g_{i}^{-1}$ to obtain a solution
of $dF=s^{*}\eta\cdot F$ on $M_{i}$. Finally, we
apply the Picard method successively on the edges of $\Gamma_{nu}$ to glue 
the local solutions to a global one. Let us denote this solution by
$F$. Then $F(x)$, $x\in\Gamma_{nu}$, is an element of the group 
$G(X_{1},X_{2})[[t]]$ (see Introduction). Moreover, since $A_{n}^{-1}(x)$,
$x\in\Gamma_{n}$, is countable, the coefficients of the series
determining $A_{2n}^{*}F$ are locally Lipschitz
on an open subset $O=S^{1}\setminus T$ of $S^{1}$ where
$T\subset S^{1}$ is countable. So
$A_{2n}^{*}F$ has the derivative almost everywhere on $S^{1}$. Then clearly
$$
\frac{d(A_{2n}^{*}F)}{dx}=(A_{2n}^{*}(s^{*}\eta))\cdot A_{2n}^{*}F=
(\sum_{i=1}^{n}a_{i}(x)t^{i}X_{1}X_{2}^{i-1})\cdot A_{2n}^{*}F\ .
$$
Without loss of generality we may assume that $A_{2n}^{*}F(0)=I$.
On the other hand, there is a global solution $H$ of the
above equation (which is equation (\ref{neq})) on $[0,2\pi)$ such that
$H(0)=I$, and such that coefficients of $H$ in the series expansion are 
Lipschitz functions.
In particular, coefficients in the series expansion of 
$H^{-1}\cdot A_{2n}^{*}F$ are locally Lipschitz on $O$. Thus we have 
$$
\frac{d(H^{-1}\cdot A_{2n}^{*}F)}{dx}=0\ \ \  {\rm almost\ everywhere\ on}
\ \ \  S^{1}.
$$
Since $H^{-1}\cdot A_{2n}^{*}F$ is continuous and 
$(H^{-1}\cdot A_{2n}^{*}F)(0)=I$, this and the results of Section 2.1 imply 
that $H(x)=A_{2n}^{*}F(x)$ on $[0,2\pi)$. But
$A_{2n}^{*}F(2\pi)=A_{2n}^{*}F(0)$ showing that the monodromy 
$\widetilde\rho_{n}$ of (\ref{neq}) is trivial.

The proof of the theorem is complete.\ \ \ \ \ \ $\Box$\\
{\bf Proof of Corollary \ref{coro2}.} 
{\bf (A)}\ \
Let $\nu :\widetilde X\rightarrow X$ be the normalization of $X$.
Then $\widetilde X$ is a 
non-compact (possibly disconnected) complex Riemann surface. Since 
$\widetilde A_{n}^{-1}(x)$ is finite for any $x\in\widetilde A(R)$, the 
definition of the normalization implies that there is a continuous map
$\widetilde A_{1n}:R\rightarrow \widetilde X$ such that 
$\widetilde A_{n}=\nu\circ\widetilde A_{1n}$. In particular, the image of
$R$ belongs to a connected component of $\widetilde X$. Thus 
$\Gamma_{n}:=A_{n}(S^{1})$ belongs to an irreducible component of $X$. 
In what
follows without loss of generality we may assume that $X$ itself is 
irreducible.
\begin{Lm}\label{polin}
Suppose that $X$ satisfies condition (1) of Corollary \ref{coro2}. Then
\penalty-10000 $\widehat\Gamma_{n}\subset X$.
\end{Lm}
{\bf Proof.} By the definition of $U$ we obtain that
$\widehat\Gamma_{n}$ belongs to a 
polynomially convex compact $K_{j}$. Since $X$ is a closed subspace of the 
Stein domain $U$, it follows (see e.g. [GR, Ch.5, Sect.4]) that $X$ is the 
set of common zeros of a family $\{f_{i}\}_{i\in I}$ of holomorphic on $U$ 
functions. 
Since $K_{j}$ is polynomially convex, by the Oka-Weil approximation theorem
each $f_{i}$ can be uniformly approximated in a
small open neighbourhood $O\subset\subset U$ of $K_{j}$ 
by holomorphic polynomials. Assume, to the
contrary, that there exists a $z\in\widehat\Gamma_{n}$ such that 
$z\not\in X\cap K_{j}$. Then there is an index $i\in I$ and an $\epsilon>0$ 
such that $|f_{i}(z)|>\epsilon$. Moreover, the polynomial 
approximation of $f_{i}$ produces a holomorphic polynomial
$p_{i}$ such that 
$$
\max_{K_{j}\cap X}|p_{i}|\leq\frac{\epsilon}{2}\ \ \ {\rm but}\ \ \
|p_{i}(z)|>\frac{\epsilon}{2}\ . 
$$
This contradicts to the fact that 
$z\in\widehat\Gamma_{n}\subset\widehat{K_{j}\cap X}$.\ \ \ \ \ $\Box$

Next, there is an open connected neighbourhood $Y\subset\subset X$ of
the compact connected set $\widehat\Gamma_{n}\subset X$.
Since $Y$ is an analytic space, by the Lojasiewicz
theorem [Lo] it is triangulated. In particular, there is an open 
connected neighbourhood $O_{Y}\subset\subset U$ of $Y$ and a retraction
$r:O_{Y}\rightarrow Y$. Then from the triviality of the monodromy
$\widetilde\rho_{n}$, as in the proof of Corollary \ref{coro1}, we
obtain under hypothesis (1) of Corollary \ref{coro2} that the path
$A_{n}:S^{1}\rightarrow Y$ is contractible in $Y$ to a point. 

Further, for one-dimensional Stein spaces $X$ and $\widetilde X$
the normalization map $\nu :\widetilde X\rightarrow X$ induces an embedding 
homomorphism of the fundamental groups. In particular, 
$\widetilde A_{1n}:S^{1}\rightarrow\widetilde X$ is contractible to a point in
$\widetilde X$. Let $p:\widetilde X_{u}\rightarrow\widetilde X$ be the 
universal covering. Then by the covering homotopy theorem there is a 
covering map
$A_{2n}:S^{1}\rightarrow\widetilde X_{u}$ such that 
$\widetilde A_{1n}=p\circ A_{2n}$.
Since $\widetilde X$ is hyperbolic, $\widetilde X_{u}$ is conformally 
equivalent to the unit disk $\Di$. Moreover, 
$A_{2n}(S^{1})\subset\Di$ is compact. Therefore
we can choose a domain $D\subset\subset\Di$ such that $A_{2n}(S_{1})\subset D$
and $D$ is conformally equivalent to $\Di$. Now we set
$A_{1n}:=(\nu\circ p)|_{D}$. 

To finish the first part of the proof let us check 
that $A_{2n}$ is locally Lipschitz outside a finite set in $S^{1}$.

Let $X_{s}$ be a finite set of singular points of $X$
containing in $\Gamma_{n}$. By our hypotheses, the preimage of each point of
$A_{n}:S^{1}\rightarrow X$ is finite. 
Thus $Y:=A_{n}^{-1}(X_{s})\subset S^{1}$ is finite. Now by 
the definition
of the normalization map, there is an open connected neighbourhood 
$V\subset X$ of $\Gamma_{n}$ such that $V\setminus X_{s}$ is smooth, and there
is the holomorphic inverse map 
$\nu^{-1}:V\setminus X_{s}\rightarrow\widetilde X$.
Further, $p:\widetilde X_{u}\rightarrow\widetilde X$ is a locally 
biholomorphic map. Finally,
the lifted map $A_{2n}$ outside $Y$ locally is the composite 
$p^{-1}\circ\nu^{-1}\circ A_{n}$. Since $A_{n}$ is Lipschitz, this implies 
that $A_{2n}$ is locally Lipschitz outside $Y$. Note that if 
$\widetilde A_{n}:R\rightarrow X$ is a holomorphic map of an open annulus
containing $S^{1}$, then the lifted map 
$\widetilde A_{1n}:R\rightarrow\widetilde X$ is holomorphic, as well. 
Moreover, since $\widetilde A_{1n}:S^{1}\rightarrow\widetilde X$ is 
contractible to a point, the covering
homotopy theorem implies that there is a holomorphic map
$A_{2n}':R\rightarrow\widetilde X_{u}$ with $A_{2n}'|_{S^{1}}=A_{2n}$ which 
covers $\widetilde A_{2n}$. In particular, there is an open annulus 
$R_{1}\subset R$ containing $S^{1}$ such that $A_{2n}'(R_{1})\subset D$. 
This proves the statement (2) of Remark \ref{rema2}.

Conversely, suppose $A_{n}=A_{1n}\circ A_{2n}$ where $A_{n}$, $A_{1n}$ 
and $A_{2n}$ satisfy hypotheses of Corollary \ref{coro2}. Let us check that 
the monodromy $\widetilde\rho_{n}$ of the corresponding equation (\ref{neq}) 
is trivial.

As in the proof of Theorem \ref{pr1'} we consider equation (\ref{lift})
defined on $U$. Its lifting to $\Di$ by $A_{1n}$ is locally solvable, because
the lifted form $A_{1n}^{*}\eta$ is a holomorphic 1-form on $\Di$.
Since $\Di$ is simply connected, the lifted equation
$dF=A_{1n}^{*}\eta\cdot F$ has a global holomorphic solution
$F$ on $U$. (This solution can be obtained by Picard iteration.)
Let us consider $A_{2n}^{*}F$ on $S^{1}$. By the hypothesis, 
it is continuous, 
locally Lipschitz outside a finite set $Y$. Therefore $A_{2n}^{*}F$ has the
derivative almost everywhere on $S^{1}\setminus Y$ such that
$$
\frac{d(A_{2n}^{*}F)}{dx}=(A_{2n}^{*}(A_{1n}^{*}\eta))\cdot A_{2n}^{*}F=
(\sum_{i=1}^{n}a_{i}(x)t^{i}X_{1}X_{2}^{i-1})\cdot A_{2n}^{*}F\ .
$$
From here as in the proof of Theorem \ref{pr1'}
we obtain that $A_{2n}^{*}F$ is Lipschitz and the monodromy 
$\widetilde\rho_{n}$ of (\ref{neq}) is trivial.

The proof of part (1) of the corollary is complete.\\
{\bf (B)}\ \  Suppose now that $X$ satisfies hypothesis (2).
As above there is a lifting $\widetilde A_{1n}:R\rightarrow\widetilde X$
of $A_{n}:R\rightarrow X$. We will show that 
$\widetilde A_{1n}:S^{1}\rightarrow\widetilde X$ is contractible to a point. 
Then the further proof repeats literally the proof presented in (A).

First, note that groups $\pi_{1}(X)$ and $\pi_{1}(\widetilde X)$ are free  
(because any 
one-dimensional Stein space is homotopically equivalent to a one-dimensional 
CW-complex [H]). Then hypothesis (2) implies that $\pi_{1}(X)$ has a finite 
number of generators. Further,
$\nu_{*}:H_{1}(\widetilde X,\Z)\rightarrow H_{1}(X,\Z)$ is an embedding 
(because $\nu$ is one-to-one outside a discrete set). This implies that
$\pi_{1}(\widetilde X)$ is finitely generated, as well.
Let $\Omega_{p}$ denote the set of holomorphic 1-forms with
polynomial coefficients on $\Co^{n}$. By  
$\nu^{*}\Omega_{p}$ we denote its pullback by $\nu$ to
$\widetilde X$, and by $\Omega(\widetilde X)$ the set of all holomorphic 
1-forms on $\widetilde X$. Because $\widetilde X$ is a Stein manifold, 
the de Rham 1-cohomology
group $H^{1}(\widetilde X,\Co)$ can be computed by $\Omega(\widetilde X)$. 
Namely, for each $\delta\in H^{1}(\widetilde X,\Co)$ there is an 
$\omega\in\Omega(\widetilde X)$ such that
\begin{equation}\label{1cogo}
\delta(\gamma)=\int_{\gamma}\omega\ \ \ {\rm for\ any}\ \ \ 
\gamma\in H_{1}(\widetilde X,\Co)\ .
\end{equation}
Since $\nu_{*}:H_{1}(\widetilde X,\Co)\rightarrow H_{1}(X,\Co)$ is an 
injection and $dim_{\Co}H_{1}(X,\Co)<\infty$,
for any $\delta\in H^{1}(\widetilde X,\Co)$ there is a 
$\xi\in H^{1}(X,\Co)$ such
that $\nu^{*}\xi=\delta$. In particular, hypothesis (2) implies that each
$\delta\in H^{1}(\widetilde X,\Co)$ can be defined by (\ref{1cogo}) with 
$\omega\in\nu^{*}\Omega_{p}$. Thus for any $\omega\in\Omega(\widetilde X)$ 
there is an $\omega'\in\nu^{*}\Omega_{p}$ such that $\omega-\omega'=df$ for a
holomorphic $f\in {\cal O}(\widetilde X)$. (This means that the embedding
$\nu^{*}\Omega_{p}\hookrightarrow\Omega(\widetilde X)$ is a 
quasi-isomorphism.)
The rest of the proof can be deduced
from Sullivan's theory of minimal models of commutative differential
graded algebras (see [Su]). For the sake of completeness we present the
sketch of the proof.

Let $N_{n}\subset GL_{n}(\Co)$ be the complex Lie subgroup of 
upper triangular unipotent matrices. By $n_{n}\subset gl_{n}(\Co)$
we denote the corresponding Lie algebras. For any homomorphism
$\rho:\pi_{1}(\widetilde X)\rightarrow N_{n}$ by $V_{\rho}$ we denote the
flat vector bundle on $\widetilde X$ constructed by $\rho$ 
(for definitions see e.g. [KN, Ch. II]). Since $N_{n}$ is contractible to a
point,  $V_{\rho}$ is topologically
trivial. Further, since $\widetilde X$ is Stein,  by the Grauert theorem [Gr]
$V_{\rho}$ is also holomorphically trivial. Therefore it is determined by 
a holomorphic flat connection on the trivial bundle 
$\widetilde X\times\Co^{n}$, that is, by an $n_{n}$-valued holomorphic 
1-form $\omega$ on $\widetilde X$.
(Note that $\widetilde X$ is one-dimensional and so the Frobenius condition
$d\omega-\omega\wedge\omega=0$ is automatically fulfilled.) The main point
of the proof is
\begin{Lm}\label{gauge}
There are a holomorphic function $F:\widetilde X\rightarrow N_{n}$ and
a holomorphic $n_{n}$-valued 1-form $\eta$ with entries from
$\nu^{*}\Omega_{p}$ such that
$$
F^{-1}\cdot\omega\cdot F-F^{-1}\cdot dF=\eta\ .
$$
\end{Lm}
{\bf Proof.} We can write $\omega=(\omega_{kr})$ as 
$\omega_{1}+(\omega-\omega_{1})$ where $\omega_{1}=(\omega_{1,kr})$, 
$\omega_{1,k 1+k}=\omega_{k 1+k}$, $k=1,\dots,n-1$, and 
$\omega_{1,kr}=0$ otherwise. 
According to hypothesis (2) and the above
arguments, there is a holomorphic matrix 1-form $\eta_{1}=(\eta_{1,kr})$
such that $\eta_{1,kr}=0$ if $r-k\neq 1$, all 
$\eta_{1,kr}\in\nu^{*}\Omega_{p}$ and $\omega_{1}-\eta_{1}=df_{1}$,
where $f_{1}=(f_{1,kr})$ is a holomorphic matrix such that
$f_{1,kr}=0$ if $r-k\neq 1$. We set $F_{1}=I_{n}+f_{1}$ where
$I_{n}\in GL_{n}(\Co)$ is the unit matrix. Then
$$
\omega_{1}'=F_{1}^{-1}\cdot\omega\cdot F_{1}-F_{1}^{-1}\cdot dF_{1}
$$
is such that $\omega_{1,k 1+k}'=\eta_{1,k 1+k}\in\nu^{*}\Omega_{p}$ for any
$k$. We will continue this process. Next, we find a holomorphic matrix
$f_{2}=(f_{2,kr})$ such that $f_{2,kr}=0$ if $r-k\neq 2$, and
$\omega_{1}'-df_{2}$ has two diagonals after the main one which consist
of elements from $\nu^{*}\Omega_{p}$.  Thus for $F_{2}=I_{n}+f_{2}$ we have
that
$$
\omega_{2}'=F_{2}^{-1}\cdot\omega_{1}'\cdot F_{2}-F_{2}^{-1}\cdot dF_{2}
$$ 
has two diagonals after the main one with entries from 
$\nu^{*}\Omega_{p}$ etc.
After $n-1$ steps we obtain the required matrix $\eta:=\omega_{n-1}'$
with entries from $\nu^{*}\Omega_{p}$. It remains to set
$F:=F_{1}\cdot F_{2}\cdots F_{n-1}$.\ \ \ \ \ \ $\Box$

Let $\widetilde\rho:\pi_{1}(\widetilde X)\rightarrow N_{n}$ be the 
representation
constructed by the flat connection $\eta$. Since by Lemma \ref{gauge}
connections $\eta$ and $\omega$ are $d$-gauge equivalent, representations
$\rho$ and $\widetilde\rho$ are conjugate (see also, e.g. [O, Sec. 4,5]).
Namely, there is a matrix $C\in N_{n}$ such that 
$C^{-1}\cdot\rho\cdot C=\widetilde\rho$. In particular,
$Ker(\rho)=Ker(\widetilde\rho)$. Now as in the proof of Theorem \ref{contr}
we obtain that if the monodromy of equation (\ref{neq}) is trivial 
the element $\gamma$ representing the path
$\widetilde A_{1n}:S^{1}\rightarrow\widetilde X$ belongs to 
$Ker(\widetilde\rho)$
(because $\widetilde\rho$ is defined by Picard iteration applied to
$\eta$). Let $G$ be the intersection of kernels of all homomorphisms
$\rho:\pi_{1}(\widetilde X)\rightarrow N_{n}$ for all $n$. Then the above
argument shows that $\gamma\in G$. But $\pi_{1}(\widetilde X)$ is a free group
with a finite number of generators. In particular, it is residually torsion
free nilpotent, meaning that $G=\{1\}$. This shows that 
$\widetilde A_{1n}:S^{1}\rightarrow\widetilde X$ is contractible. 

As we mentioned above the further proof repeats word-for-word
the proof given in (A).\ \ \ \ \ \ $\Box$\\
{\bf Proof of Corollary \ref{lorentz}.} 
For the basic facts from the Algebraic Geometry see e.g. the book of
Shafarevich [Sh].

In what follows $\Co^{*}:=\Co\setminus\{0\}$. If the coefficients
$a_{1},\dots,a_{n}$ in (\ref{neq}) are trigonometric polynomials the map
$A_{n}:S^{1}\rightarrow\Co^{n}$ can be extended to a holomorphic map
$\widetilde A_{n}:\Co^{*}\rightarrow\Co^{n}$ defined by Laurent
polynomials. If all components of $A_{n}$ are holomorphic polynomials then
the statement of the corollary is trivially fulfilled. Thus we may assume
without loss of generality that $\widetilde A_{n}$ cannot be extended to
$\Co$.

Let $\Co\P^{n}$ denote the complex projective space. Then
$\Co^{n}=\Co\P^{n}\setminus H$, where $H$ is the
hyperplane at infinity. Now, Zariski closure $\overline X$ of the image 
$X=\widetilde A_{n}(\Co^{*})$ is a (possibly singular) rational
curve and $X=\overline X\setminus H\subset\Co^{n}$ is a closed algebraic 
subvariety. Further, if the monodromy $\widetilde\rho_{n}$
of (\ref{neq}) is trivial, then from Corollary \ref{coro2} it follows that
$A_{n}:S^{1}\rightarrow X$ is contractible to a point. Let
$\nu:\Co\P^{1}\rightarrow\overline X$ be the normalization of 
$\overline X$ and $\widetilde X=\nu^{-1}(X)\subset\Co\P^{1}$ be the 
normalization of $X$. Then there is an algebraic covering map 
$\widetilde A_{1n}:\Co^{*}\rightarrow\widetilde X$ such that
$\widetilde A_{n}=\nu\circ\widetilde A_{1n}$. In particular,
$\widetilde A_{1n}:S^{1}\rightarrow\widetilde X$ is contractible to a point.
But since $\widetilde A_{1n}:\Co^{*}\rightarrow\widetilde X$ is a 
finite proper surjective map, the image of the homomorphism 
$(\widetilde A_{1n})_{*}:\pi_{1}(\Co^{*})\rightarrow
\pi_{1}(\widetilde X)$ is a subgroup of a finite index in 
$\pi_{1}(\widetilde X)$. Now the contractibility of 
$\widetilde A_{1n}:S^{1}\rightarrow\widetilde X$ implies that
$\pi_{1}(\widetilde X)=\{1\}$. Thus $X\cong\Co\subset\Co\P^{1}$. Since both
maps $\nu:\widetilde X\rightarrow\Co^{n}$ and 
$\widetilde A_{1n}:\Co^{*}\rightarrow\Co$ are algebraic,
the latter implies that there are polynomials 
$p_{1},\dots,p_{n}\in\Co[z]$ and a Laurent polynomial $q$ such that 
$\nu(z)=(p_{1}(z),\dots,p_{n}(z))$, $z\in\Co$, and
$\widetilde A_{1n}(z)=q(z)$, $z\in\Co^{*}$.
Thus we have $\widetilde a_{i}(x)=p_{i}(q(x))$,
$x\in S^{1}$, $1\leq i\leq n$ as required. 

Conversely, the above
factorization of $A_{n}$ implies that the monodromy $\widetilde\rho_{n}$ of
the corresponding equation (\ref{neq}) is trivial. This can be shown exactly
as in the proof of part (A) of Corollary \ref{coro2} and is based on the 
fact that $\widetilde X=\Co$ is contratcible to a point.\ \ \ \ \ $\Box$
\sect{\hspace*{-1em}. Proofs of Theorems 2.18, 2.20 and Corollaries
2.21 and 2.22.}
{\bf Proof of Theorem \ref{tef1}.}
According to the hypothesis any curve $\Gamma_{k}\subset\Co^{k}$, $k\geq n$,
is piecewise smooth. Let $S_{k}$ be the finite set consisting of 
singularities of $\Gamma_{k}$ and critical values of $A_{k}$.
Without loss of generality we may assume 
that each $S_{k}$ is not empty. By $p_{k}:\Gamma_{k+1}\rightarrow\Gamma_{k}$
we denote the map induced by the projection
$\Co^{k+1}\rightarrow\Co^{k}$ to the first $k$ coordinates. Then
$p_{k}:\Gamma_{k+1}\rightarrow\Gamma_{k}$, $k\geq n$, is a finite 
surjective map. It is clear that $S_{k+1}\subset p_{k}^{-1}(S_{k})$. Then
from finiteness it follows that there is a number $K\geq n$ such that
for all $k\geq K$, all finite sets $A_{k}^{-1}(S_{k})\subset S^{1}$ coincide
with a finite set, say, $Y\subset S^{1}$. Then $S^{1}\setminus Y$ is the
disjoint union of  a finite number of open intervals 
$I_{1},\dots, I_{l}$ such that each $I_{s}$ is
diffeomorphic by $A_{k}$ to one of the smooth connected components of
$\Gamma_{k}\setminus S_{k}$. Moreover, we have 
$p_{k}^{-1}(S_{k})=S_{k+1}$ and $p_{k}:S_{k+1}\rightarrow S_{k}$ is a
bijection. Based on this and using finiteness of each $p_{k}$ we can 
find a number $N\geq K$ such that for each $k\geq N$ we obtain that 
$p_{k}:\Gamma_{k+1}\rightarrow\Gamma_{k}$ is a homeomorphism (in fact, it
is even diffeomorphism outside sets $S_{k+1}$ and $S_{k}$).

Suppose now that $\widehat\Gamma_{k}=\Gamma_{k}$ for any $k$, and that
the monodromy $\widetilde\rho_{N}$ of equation (\ref{neq}) is trivial.
Then according to Corollary \ref{coro1} and Theorem \ref{pr1'} this
is equivalent to the contractibility of $A_{N}:S^{1}\rightarrow\Gamma_{N}$
to a point. But it follows from above that 
$p_{N}\circ p_{N+1}\circ\cdots\circ p_{k}:\Gamma_{k}\rightarrow\Gamma_{N}$ 
is a homeomorphism for any $k\geq N$. Moreover, 
$A_{N}=p_{N}\circ p_{N+1}\circ\cdots\circ p_{k}\circ A_{k}$. This implies
that $A_{k}:S^{1}\rightarrow\Gamma_{k}$, $k\geq N$, is contractible to a
point. Then from Corollary \ref{coro1} and Theorem \ref{pr1'} we obtain
that the monodromy $\widetilde\rho_{k}$ is trivial for any $k\geq N$.
In particular, according to Proposition \ref{integr1},
the monodromy $\widetilde\rho$ of (\ref{formal1}) is trivial.

The proof of the theorem is complete.\ \ \ \ \ $\Box$\\
{\bf Proof of Theorem \ref{tef2}.}
According to Proposition \ref{integr} the monodromy of the equation 
$F'(x)=(\sum_{i=1}^{k}b_{i}'(x)t^{i}X_{1}X_{2}^{i-1})\cdot F(x)$ 
is trivial if and only if
\begin{equation}\label{bint}
\int\cdots\int_{0\leq s_{1}\leq\cdots\leq s_{k}\leq 2\pi}
b_{i_{k}}'(s_{k})\cdots b_{i_{1}}'(s_{1})\ ds_{k}\cdots ds_{1}
\end{equation}
is zero for any positive integers $i_{1},\dots, i_{k}$ and any $k$.
Since each $a_{i}$ is the uniform limit of functions of the
form $\sum_{j=1}^{k}p_{ji}(b_{1},\dots,b_{k})\cdot b_{j}'$ where
$p_{ji}\in\Co[z_{1},\dots,z_{k}]$ are holomorphic polynomials,
from the Ree formula (see 
Section 4.3 in the proof of Theorem \ref{contr}) it follows 
directly that each iterated integral
\begin{equation}\label{aint}
\int\cdots\int_{0\leq s_{1}\leq\cdots\leq s_{l}\leq 2\pi}
a_{i_{l}}(s_{l})\cdots a_{i_{1}}(s_{1})\ ds_{l}\cdots ds_{1}
\end{equation}
is the limit of integrals of the form (\ref{bint}). This implies
that all integrals in (\ref{aint}) are zeros. Then Proposition
\ref{integr} implies that the monodromy $\widetilde\rho$ of the 
corresponding equation (\ref{formal1}) is trivial.\ \ \ \ \ $\Box$\\
{\bf Proof of Corollary \ref{planar}.}
Passing to polar coordinates we obtain the equation
$$
\frac{dr}{d\phi}=\frac{P}{1+Q}r\ \ \ \ {\rm where}\ \ \ 
P(r,\phi)=\frac{\partial A(r,\phi)}{\partial\phi}\ ,
\ \ \ Q(r,\phi)=B(r,\phi)\ .
$$ 
Let us write $H(r,\phi)=h(\phi)r^{k}$ where $h$ is a trigonometric polynomial
of degree $k$. Then we have
$P(r,\phi)=\sum_{i=1}^{\infty}a_{i}h^{i-1}(\phi)h'(\phi)r^{ki}$,
$Q(r,\phi)=\sum_{i=0}^{\infty}b_{i}h^{i}(\phi)r^{ki}$, where
$a_{i},b_{i}\in\Co$.
Now, if we expand the right-hand side of the above differential equation
as a series in $r$, then each coefficient of this series will be
of the form $(p\circ h)\cdot h'$ where $p\in\Co[z]$. In particular, the
first integrals of these coefficients satisfy the hypotheses of 
Corollary \ref{lorentz}. From this corollary it follows that the
required vector field determines a center.\ \ \ \ \ $\Box$\\
{\bf Proof of Corollary \ref{sim}.}
Passing to polar coordinates we obtain an equation
$$
\frac{dr}{d\phi}=\frac{P}{1+Q}r\ \ \ {\rm where}\ \ \
P(r,-\phi)=-P(r,\phi),\ \ Q(r,-\phi)=Q(r,\phi)\ .
$$
In particular $P(r,\phi)=\sum_{j=1}^{\infty}p_{j}(\phi)r^{j}$ where
$p_{j}$ are odd trigonometric polynomials, and 
$Q(r,\phi)=\sum_{j=0}^{\infty}q_{j}(\phi)r^{j}$ where $q_{j}$ are even 
trigonometric polynomials. Since each $p_{j}$ has the Fourier series
expansion by functions $\sin(n\phi)$, $n=1,2,\dots$, the first integrals
$\widetilde p_{j}(\phi)=\int_{0}^{\phi}p_{j}(s)ds$ of each $p_{j}$ are 
even trigonometric polynomials. Thus $\widetilde p_{j}$ and $q_{j}$ are 
polynomials in $\cos\phi$. Now
the required result follows from Corollary \ref{lorentz}.\ \ \ \ \ $\Box$
{\sect{\hspace*{-1em}. Proofs of Results of Section 3.1.}
In this section we prove Propositions \ref{equivrelat}, \ref{oper},
\ref{sepit} and Theorem \ref{group}.\\
{\bf Proof of Proposition \ref{equivrelat}.} Let us consider equation
({\ref{formal1}) related to $a\in X$ introduced in Section 2.2. By
$\omega_{a}(x):=\sum_{i=1}^{\infty}a_{i}(x)t^{i}X_{1}X_{2}^{i-1}$ we
denote the coefficient of (\ref{formal1}) and by
$\widetilde\rho_{a}:\Z\rightarrow G(X_{1},X_{2})[[t]]$ the
corresponding monodromy. Let $F_{a}:\Re\rightarrow G(X_{1},X_{2})[[t]]$
be the fundamental solution of ({\ref{formal1}) (see Section 2.3). Then
$\widetilde\rho_{a}(n):=F_{a}(2\pi n)=F_{a}(2\pi)^{n}$, $n\in\Z$. 
By the definition for $a*b$, $a,b\in X$, we have 
$$
F_{a*b}'(x):=\left\{
\begin{array}{cll}
2\omega_{a}(2x)\cdot F_{a*b}(x)&{\rm if}&0<x\leq\pi\\
2\omega_{b}(2x-2\pi)\cdot F_{a*b}(x)&{\rm if}&\pi< x\leq 2\pi
\end{array}
\right. 
$$
From here we have 
$$
F_{a*b}(x)=\left\{
\begin{array}{cll}
F_{a}(2x)&{\rm if}&0<x\leq\pi\\
F_{b}(2x-2\pi)\cdot F_{a}(2\pi)&{\rm if}&\pi<x\leq 2\pi
\end{array}
\right.
$$
In particular, 
\begin{equation}\label{91}
\widetilde\rho_{a*b}(1):=F_{a*b}(2\pi)=F_{b}(2\pi)\cdot F_{a}(2\pi)=
\widetilde\rho_{b}(1)\cdot\widetilde\rho_{a}(1)\ .
\end{equation}
Similarly, for $b^{-1}\in X$ we obtain that
\begin{equation}\label{92}
\widetilde\rho_{b^{-1}}(1)=(\widetilde\rho_{b}(1))^{-1}\ .
\end{equation}

Suppose now that $a,b\in X$ and $a\sim b$. This means that 
$a*b^{-1}\in {\cal U}$ or, equivalently, that
$\widetilde\rho_{a*b^{-1}}(1)=1$.
Then the above formulae imply that $a\sim b$ if and only if
\begin{equation}\label{93} 
\widetilde\rho_{a}(1)=\widetilde\rho_{b}(1)\ .
\end{equation}

Now, the fact that $\sim $ is an equivalence relation on $X$ follows
easily from (\ref{93}).

The proof of the proposition is complete.\ \ \ \ \ $\Box$\\
{\bf Proof of Proposition \ref{oper}.} The proof follows straightforwardly
from equations (\ref{91}), (\ref{92}) and (\ref{93}). We leave the details
to the reader.\ \ \ \ \ $\Box$\\
{\bf Proof of Proposition \ref{sepit}.} Suppose $a,b\in X$ and $a\sim b$.
We need to show that $I(a)=I(b)$ for any iterated integral $I$. Clearly, it
suffices to prove the statement for the basic iterated integrals
(see Section 2.1 for the definition). 

As before, suppose $F_{a}$ and $F_{b}$ are the fundamental solutions of 
equations (\ref{formal1}) corresponding to $a$ and $b$, respectively. It
follows from the formula (\ref{e3'}) that 
$$
F_{a}(2\pi)=\sum_{i=0}^{\infty}q_{ia}(2\pi,X_{1},X_{2}),\ \ \ 
F_{b}(2\pi)=\sum_{i=0}^{\infty}q_{ib}(2\pi,X_{1},X_{2})\ ,
$$
where $q_{0a}=q_{0b}=I$ and others $q_{ia}$, $q_{ib}$ are polynomials in 
$R_{1},R_{2},\dots$ of degree $i$.
Moreover, iterated integrals 
$I_{i_{1},\dots,i_{k}}(a)$ and
$I_{i_{1},\dots,i_{k}}(b)$ with $\sum_{s=1}^{k}i_{s}=i$
are coefficients before the monomial 
$(X_{1}X_{2}^{i_{k}-1})\cdots (X_{1}X_{2}^{i_{1}-1})$ in $q_{ia}$ and 
$q_{ib}$, respectively. 
Since, by the hypothesis, $F_{a}(2\pi)=F_{b}(2\pi)$, Lemma \ref{free}
then implies that $I_{i_{1},\dots,i_{k}}(a)=I_{i_{1},\dots,i_{k}}(b)$ for
all possible integers $i_{1},\dots, i_{k}, k$.

The proposition is proved.\ \ \ \ \ $\Box$\\
{\bf Proof of Theorem \ref{group}.} 
{\bf (1)} Suppose, to the contrary, that there are $g_{1},g_{2}\in G(X)$ such
that $g_{1}\neq g_{2}$ and $\widehat I(g_{1})=\widehat I(g_{2})$ for every
$\widehat I\in {\cal I}(G(X))$. Let $\tilde g_{1},\tilde g_{2}\in X$ be
such that $[\tilde g_{1}]=g_{1}$, $[\tilde g_{2}]=g_{2}$. Then by 
definition we have $I(\tilde g_{1})=I(\tilde g_{2})$ for any \penalty-10000
iterated integral $I$. Then the arguments used in the proof of
Proposition \ref{sepit} show that 
$\widetilde\rho_{\tilde g_{1}}(1)=\widetilde\rho_{\tilde g_{2}}(1)$. This
implies that $\tilde g_{1}\sim\tilde g_{2}$ and so $g_{1}=g_{2}$. 
This contradiction proves (1). 
\\ 
{\bf (2)} First we will show that $G(X)$ equipped with the
topology $\tau$ is a topological group. (Note that $G(X)$ is a
Hausdorff space by the definition of $\tau$.) 
This will be done by showing that operations of multiplication 
$\cdot :G(X)\times G(X)\rightarrow G(X)$ and of taking the inverse
$^{-1}:G(X)\rightarrow G(X)$ are continuous maps. 

Let $(g,h)\in G(X)\times G(X)$ and $U$ be an open neighbourhood
of $gh\in G(X)$ of the form
$$
U=\{s\in G(X):\ \max_{1\leq i\leq k}|\widehat I_{i}(s)-\widehat I_{i}(gh)|<
\epsilon\}
$$
where $I_{1},\dots, I_{k}$ are basic iterated integrals, and $0<\epsilon<1$.
We recall the 
following property of iterated integrals (see e.g. 
[Ha, Proposition 2.9]). For any $g_{1},g_{2}\in G(X)$,
\begin{equation}\label{94}
\widehat I_{i_{1},\dots,i_{k}}(g_{1}g_{2})=
\widehat I_{i_{1},\dots, i_{k}}(g_{1})+
\sum_{j=1}^{k-1}\widehat I_{i_{j+1},\dots,i_{k}}(g_{1})\cdot 
\widehat I_{i_{1},\dots,i_{j}}(g_{2})+
\widehat I_{i_{1},\dots, i_{k}}(g_{2})\ .
\end{equation}
Thus without loss of generality we may assume that
$\widehat I_{l}(g_{1}g_{2})=\sum_{j=1}^{t_{l}}R_{lj}(g_{1})
\cdot S_{lj}(g_{2})$,
$1\leq l\leq k$, where $R_{lj}$ and $S_{lj}$ are basic iterated integrals
on $G(X)$. We set
$$
M_{1}=\max_{j,l}|R_{jl}(g)|,\ \ \ M_{2}=\max_{j,l}|S_{jl}(h)|,\ \ \
M=M_{1}+M_{2}+1,\ \ \ t=\max_{l}\ t_{l}\ .
$$

Let us consider open neighbourhoods $U_{1}$ and $U_{2}$ of $g$ and $h$
defined by
$$
\begin{array}{c}
\displaystyle
U_{1}=\left\{s_{1}\in G(X)\ :\ \max_{l,j}|R_{lj}(g)-R_{lj}(s_{1})|<
\frac{\epsilon}{Mt}\right\}\ ,\\
\\
\displaystyle
U_{2}=\left\{s_{2}\in G(X)\ :\ \max_{l,j}|S_{lj}(h)-S_{lj}(s_{2})|<
\frac{\epsilon}{Mt}\right\}\ .
\end{array}
$$
Here $R_{lj}$ and $S_{lj}$ are from the previous formulae. Then for 
$s_{1}s_{2}$ belonging to the image of $U_{1}\times U_{2}$ under 
multiplication we have (for $1\leq l\leq k$),
$$
\begin{array}{c}
\displaystyle
|\widehat I_{l}(gh)-\widehat I_{l}(s_{1}s_{2}|\leq
\sum_{j=1}^{t_{l}}(|R_{lj}(g)-R_{lj}(s_{1})|\cdot |S_{lj}(h)|+
|R_{lj}(s_{1})|\cdot |S_{lj}(h)-S_{lj}(s_{2})|)\leq
\\
\displaystyle
M_{2}t_{l}\left(\frac{\epsilon}{Mt}\right)+ 
\left(M_{1}+\frac{\epsilon}{Mt}\right)t_{l}
\left(\frac{\epsilon}{Mt}\right)<\epsilon\ .
\end{array}
$$
This shows that the image of $U_{1}\times U_{2}$ is containing in $U$.
Therefore we established that multiplication
$\cdot : G(X)\times G(X)\rightarrow G(X)$ is a continuous map.

The continuity of $^{-1}:G(X)\rightarrow G(X)$ follows directly from the
definition of the topology $\tau$ and from
the formula (see e.g. [Ha, Proposition 2.12])
\begin{equation}\label{95}
\widehat I_{i_{1},\dots,i_{k}}(g^{-1})=(-1)^{k}
\widehat I_{i_{1},\dots,i_{k}}(g)\ ,\ \ \ g\in G(X)\ .
\end{equation}

Also, the fact that $G(X)$ is separable follows from (3), because
every compact metric space is separable and therefore the union of at most
countable number of compact metric spaces is separable, as well.\\
{\bf (3)} Let us prove that $G(X)$ is metrizable.

By definition the set $\{\widehat I_{i_{1},\dots,i_{k}}\}$ of basic 
iterated integrals is countable. Let us arrange the
elements of this set in a sequence $\{J_{n}\}$ of distinct elements.
Now we define
$$
d(g,h)=\sum_{i=1}^{\infty}
\frac{2^{-i}|J_{i}(g)-J_{i}(h)|}{1+|J_{i}(g)-J_{i}(h)|}\ ,\ \ \ g,h\in G(X)\ .
$$ 
It is easy to see that $d$ is a metric in $G(X)$. Let us prove that $d$ is 
compatible with the topology $\tau$ on $G(X)$. Let
$$
B_{r}(g):=\{h\in G(X)\ :\ d(h,g)<r\}\ ,\ \ \ r>0\ ,
$$
be a ball centered at $g$ of radius $r$. Since by  definition for
every fixed $g$ the function $d(g,\cdot)$ is continuous on $(G(X),\tau)$, 
every ball $B_{r}(g)$ is an open subset of $G(X)$.
Suppose 
$$
U=\{h\in G(X):\ \max_{1\leq i\leq k}|J_{l_{i}}(g)-
J_{l_{i}}(h)|<\epsilon\}\ ,\ \ \  0<\epsilon<1\ ,
$$
is an open neighbourhood of $g$. Let $N:=\max_{1\leq i\leq k}l_{i}$. We
set $r:=\frac{\epsilon}{2^{N+1}}$. It is easy to see that if 
$h\in B_{r}(g)$ then
$$
\frac{2^{-l_{i}}|J_{l_{i}}(g)-J_{l_{i}}(h)|}{1+|J_{l_{i}}(g)-J_{l_{i}}(h)|}<r
\ \ \ {\rm for}\ \ \ 1\leq i\leq k\ .
$$
This is equivalent to the inequalities
$$
|J_{l_{i}}(g)-J_{l_{i}}(h)|<\frac{2^{l_{i}}r}{1-2^{l_{i}}r}<
\frac{\epsilon/2}{1-\epsilon/2}<\epsilon\ ,\ \ \ 1\leq i\leq k\ .
$$
Thus $B_{r}(g)\subset U$. This shows that $d$ is compatible with
$\tau$.
\begin{R}\label{notcomplete}
{\rm One can show that the metric space $(G(X),d)$ is not complete.}
\end{R}

Prove now that $G(X)$ is the union of an increasing sequence of compact
sets.

Given a natural number $n$, let $V_{n}\subset X$ consist of elements 
$a=(a_{1},a_{2},\dots)$ such that 
$$
\sup_{x\in S^{1}}|a_{i}(x)|\leq n^{i}\ ,\ \ \ i=1,2,\dots\ .
$$
By $K_{n}$ we denote the image of $V_{n}$ in $G(X)$. 
Clearly, $X=\cup_{n\geq 1}V_{n}$ and \penalty-10000 
$G(X)=\cup_{n\geq 1}K_{n}$. 
Next, let us prove that $K_{n}$ is compact.

Recall that $X_{i}:=L^{\infty}(S^{1})$ is the space of all coefficients $
a_{i}$ from equation (\ref{e2}). We consider every $X_{i}$ in the weak-star
topology. According to the Banach-Alaoglu theorem
the set
$$
V_{ni}:=\{f\in X_{i}\ :\ \sup_{x\in S^{1}}|f(x)|\leq n^{i}\}
$$
is weak-star compact. Also, by Tychonoff's theorem, 
$V_{n}=\prod_{i\geq 1}V_{ni}\subset X$ is compact in the product topology.
Since $K_{n}$ is a metric space, in order to prove that $K_{n}$ is
compact it suffices to show that every infinite sequence 
$\{x_{k}\}\subset K_{n}$ has a limit point in $K_{n}$. 

Given $\{x_{k}\}\subset K_{n}$, let $\{y_{k}\}\subset V_{n}$ be such that
$[y_{k}]=x_{k}$ for every $k$. Since $V_{n}$ is compact, without loss
of generality we may assume that $\lim_{k\to\infty}y_{k}=y\in V_{n}$
(here $V_{n}$ is equipped with the above product topology). Let us
prove that $\lim_{k\to\infty}x_{k}=[y]$. It is easy to see that the latter
is equivalent to the fact that
$$
\lim_{k\to\infty}I_{i_{1},\dots,i_{l}}(x_{k})=I_{i_{1},\dots,i_{l}}(y)
$$
for every basic iterated integral $I_{i_{1},\dots,i_{l}}$. Let
$x_{k}=(x_{1k},x_{2k},\dots)$ and $y=(y_{1},y_{2},\dots)$ where 
$x_{ik},y_{i}\in V_{ni}$ for $i=1,2,\dots$.
Then 
$$
I_{i_{1},\dots,i_{l}}(x_{k})=
\int\cdots\int_{0\leq s_{1}\leq\cdots\leq s_{l}\leq 2\pi}
x_{i_{l}k}(s_{l})\cdots x_{i_{1}k}(s_{1})\ ds_{l}\cdots ds_{1}\ .
$$
Our proof is based on the following result.
\begin{Lm}\label{weak}
Let $\{f_{k}\},\ \{g_{k}\}\subset L^{\infty}([0,2\pi])$ be uniformly
bounded sequences convergent in the weak-star topology to 
$f,g\in L^{\infty}([0,2\pi])$, respectively. Let 
$$
\tilde g_{k}(x)=\int_{0}^{x}g_{k}(s)ds\ \ \ {\rm and}\ \ \ 
\tilde g(x)=\int_{0}^{x}g(s)ds\ ,
\ \ \ 0\leq x\leq 2\pi\ .
$$
Then the sequence $\{f_{k}\tilde g_{k}\}$ converges in the weak-star
topology of $L^{\infty}([0,2\pi])$ to $f\tilde g$.
\end{Lm}
{\bf Proof.} According to the hypothesis, there exists a constant $M<\infty$
such that for any $k$,
$$
\sup_{x\in [0,2\pi]}\max\{|f_{k}(x)|,|g_{k}(x)|\}<M\ .
$$
From here we obtain that $\tilde g_{k}$ converges uniformly on
$[0,2\pi]$ to $\tilde g$, and that for any $k$,
$$
\sup_{x\in [0,2\pi]}\max\{|\tilde g_{k}(x)|,|\tilde g(x)|\}\leq 2\pi M\ .
$$
The latter inequality is obvious. Let us check the first statement. 
Given $\epsilon>0$ choose points $0=t_{1}<t_{2}<\dots <t_{r}=2\pi$ 
from $[0,2\pi]$ such that $t_{i+1}-t_{i}<
\frac{\epsilon}{3M}$ for any $i$.
Since $\{g_{k}\}$ converges in the weak-star topology to $g$, for any
$x\in [0,2\pi]$ we obtain that
$$
\tilde g_{k}(x)=\int_{0}^{2\pi}g_{k}(s)\chi_{[0,x]}(s)ds\rightarrow
\int_{0}^{2\pi}g(s)\chi_{[0,x]}(s)ds=\tilde g(x)\ \ \  {\rm as}\ \ \
k\to\infty\ .
$$
Here $\chi_{[0,x]}$ is the characteristic function of $[0,x]$. In particular,
there exists a natural number $P$ such that for any $p\geq P$ we have 
$$
\max_{1\leq i\leq r}|\tilde g_{p}(t_{i})-\tilde g(t_{i})|<\epsilon/3\ .
$$
Next, for an $x\in [0,2\pi]$ there is a $t_{i}$ such that 
$|x-t_{i}|<\frac{\epsilon}{3M}$. Thus for $p\geq P$ we have
$$
\begin{array}{cc}
\displaystyle
|\tilde g_{p}(x)-\tilde g(x)|\leq  |\tilde g(x)-\tilde g(t_{i})|+
|\tilde g(t_{i})-\tilde g_{p}(t_{i})|+|\tilde g_{p}(t_{i})-\tilde g_{p}(x)|<\\
\\
\displaystyle
M\cdot\left(\frac{\epsilon}{3M}\right)+\frac{\epsilon}{3}+
M\cdot\left(\frac{\epsilon}{3M}\right)=\epsilon\ .
\end{array}
$$
This shows that $\{\tilde g_{k}\}$ converges uniformly on $[0,2\pi]$ to
$\tilde g$.

Now, given $h\in L^{1}([0,2\pi])$ we have
$$
\begin{array}{cc}
\displaystyle
\left|\int_{0}^{2\pi}(f_{k}(s)\tilde g_{k}(s)-f(s)\tilde g(s))
\overline{h(s)}ds\right|\leq
\int_{0}^{2\pi}|f_{k}(s)|\cdot (|\tilde g_{k}(s)-\tilde g(s)|)\cdot |h(s)|ds+
\\
\\
\displaystyle
\left|\int_{0}^{2\pi}(f_{k}(s)-f(s))\tilde g(s)\overline{h(s)}ds\right|\leq
M\cdot\sup_{s\in [0,2\pi]}|\tilde g_{k}(s)-\tilde g(s)|\cdot
\int_{0}^{2\pi}|h(s)|ds+\\
\\
\displaystyle
\left|\int_{0}^{2\pi}(f_{k}(s)-f(s))\tilde g(s)\overline{h(s)}ds\right|\ .
\end{array}
$$
The last two terms in this inequality tend to $0$ as $k$ tends to $\infty$.
For the first term it follows from the uniform convergence of 
$\{\tilde g_{k}\}$ to $\tilde g$, whereas for the second one because
$\overline{\tilde g}h\in L^{1}([0,2\pi])$ and
$\{f_{k}\}$ converges in the weak-star topology to $f$. Thus
$$
\lim_{k\to\infty}\int_{0}^{2\pi}f_{k}(s)\tilde g_{k}(s)\overline{h(s)}ds=
\int_{0}^{2\pi}f(s)\tilde g(s)\overline{h(s)}ds\ \ \ {\rm for\ any}\
\ \ h\in L^{1}([0,2\pi])\ .
$$
This is equivalent to the fact that $f_{k}\tilde g_{k}$ converges in
the weak-star topology to $f\tilde g$.

The proof of the lemma is complete.\ \ \ \ \ $\Box$

Now the statement $\lim_{k\to\infty}I_{i_{1},\dots,i_{l}}(x_{k})=
I_{i_{1},\dots,i_{l}}(y)$
follows by induction on $l$ from Lemma \ref{weak} if we use the fact that 
every $\{x_{i_{r}k}\}_{k\geq 1}$ converges in the weak-star
topology of $L^{\infty}([0,2\pi])$ to $y_{i_{r}}$, $1\leq r\leq l$.
We leave the details to the reader. 

Therefore we have $\lim_{k\to\infty}x_{k}=[y]$. According to the preceding
discussion this shows that $K_{n}$ is compact, and completes the proof of
(3).\\
{\bf (4)} Let us prove first that $G(X)$ is contractible to $1\in G(X)$.
To this end for any $t\in\Re_{+}$ and $a=(a_{1},a_{2},\dots,)\in X$ we
define 
$$
ta:=(ta_{1},ta_{2},\dots)\ .
$$
Clearly, $ta$ satisfies inequalities (\ref{bound}) and so it belongs to $X$.
Also, it is easy to see that if $a\sim b$ then $ta\sim tb$.
Now we define a map $\Phi:[0,1]\times G(X)\rightarrow G(X)$ by the formula
$$
\Phi(t,[a]):=[ta]\ .
$$
(Recall that $[a]$ stands for the image of $a\in X$ in $G(X)$.)
Let us prove that $\Phi$ is a continuous map. Since $[0,1]\times G(X)$ and
$G$ are metrizable, it suffices to check that if a sequence
$\{(t_{k},g_{k})\}\subset [0,1]\times G(X)$ converges to $(t,g)\in G(X)$,
then $\{\Phi(t_{k},g_{k})\}$ converges to $\Phi(t,g)$. In turn, the
latter is equivalent to the fact that
$$
\lim_{k\to\infty}I_{i_{1},\dots,i_{l}}(t_{k}a_{k})=I_{i_{1},\dots,i_{l}}(ta)
$$
for every basic iterated integral $I_{i_{1},\dots,i_{l}}$. Here
$a_{k}, a\in X$ are such that $[a_{k}]=g_{k}$ and $[a]=g$. Now, by the
definition we have
$$
I_{i_{1},\dots,i_{l}}(t_{k}a_{k})=t_{k}^{l}I_{i_{1},\dots,i_{l}}(a_{k})\ \ \
{\rm and}\ \ \ I_{i_{1},\dots,i_{l}}(ta)=t^{l}I_{i_{1},\dots,i_{l}}(a)\ .
$$
According to our hypothesis $t_{k}^{l}\to t^{l}$ and
$I_{i_{1},\dots,i_{l}}(a_{k})\to I_{i_{1},\dots,i_{l}}(a)$ as $k\to\infty$,
which implies the required result. 

Further, observe that $\Phi(1,\cdot)=id$ and $\Phi(0,\cdot)$ maps
$G(X)$ to $\{1\}$. Thus $\Phi$ determines the required contraction.

Since $G(X)$ is a topological group, to prove that $G(X)$ is arcwise 
connected it suffices to join every $g\in G(X)$ by a continuous path with 
$1\in G(X)$. But then such a path is given by the formula
$$
\gamma(t):=\Phi(t,g)\ ,\ \ \  t\in [0,1]\ .
$$

Let us prove now that $G(X)$ is locally simply and arcwise connected.
As before, it suffices to prove that every open neighbourhood
of $1\in G(X)$ of the form
$$
U=\{g\in G(X)\ :\ \max_{1\leq i\leq l}|\widehat I_{i}(g)|<\epsilon\}\ ,
\ \ \ \epsilon>0\ ,
$$
where $I_{1},\dots, I_{l}$ are non-constant basic iterated integrals,
is contractible inside $U$ to $1$.  Now for every $i$, $1\leq i\leq l$,
there is a natural number
$s_{i}$ such that $I_{i}(ta)=t^{s_{i}}I_{i}(a)$ for every 
$(t,a)\in [0,1]\times X$. In particular, for $g\in U$ we have
$$
\max_{1\leq i\leq l}|\widehat I_{i}(\Phi(t,g))|=
\max_{1\leq i\leq l}t^{s_{i}}|\widehat I_{i}(g)|<
\epsilon\cdot \max_{1\leq i\leq l}t^{s_{i}}\leq\epsilon\ .
$$
This shows that $\Phi$ maps $[0,1]\times U$ into $U$, and so $U$ is 
contractible inside $U$ to $1$.

The proof of (4) is complete.\\
{\bf (5)} Let us prove that $G(X)$ is residually torsion free nilpotent.
Let us consider a subset $H\subset G(X_{1},X_{2})[[t]]$ (see Section 2.2 for
the definition) consisting of series of the form
$$
s=\sum_{i=0}^{\infty}p_{i}(X_{1},X_{2})t^{i}
$$
where $p_{0}=I$ and $p_{i}$ are homogeneous polynomials of degree $i$ in
free non-commutative variables $X_{1}$ and $X_{2}$ with complex
coefficients. It is easy to see that $H$ is a group with respect to the
multiplication in $G(X_{1},X_{2})[[t]]$. Further, according to
equations (\ref{91}), (\ref{92}) and (\ref{93}) there exists a natural
injective homomorphism $\phi: G(X)\rightarrow H$ defined by the 
formula
$$
\phi([a]):=\widetilde\rho_{a}(1)\ ,\ \ \ a\in X\ ,
$$
where $1$ is the unit in $\Z$ and $\widetilde\rho_{a}$ is the monodromy
of equation (\ref{formal1}) corresponding to $a\in X$. Thus it suffices
to prove that $H$ is a residually torsion free nilpotent group.

First, observe that the number of monomials of degree $i$ in variables 
$X_{1},X_{2}$ is finite. Thus the vector space (over $\Co$) of homogeneous
polynomials of degree $i$ in $X_{1},X_{2}$ with complex coefficients is
finite-dimensional. Let us consider a subset $H_{n}\subset H$ consisting
of series of the form
$$
s=I+\sum_{i=n+1}^{\infty}p_{i}(X_{1},X_{2})t^{i}\ .
$$
\begin{Lm}\label{l92}
$H_{n}$ is a normal subgroup of $H$.
\end{Lm}
{\bf Proof.} It easy to see that $s_{1}s_{2}\in H_{n}$ for any 
$s_{1},s_{2}\in H_{n}$. Thus  $H_{n}\subset H$ is a subgroup.

Let $s=I+\sum_{i=1}^{\infty}p_{i}(X_{1},X_{2})t^{i}\in H$ and
$s^{-1}=I+\sum_{i=1}^{\infty}\widetilde p_{i}(X_{1},X_{2})t^{i}$.
From the identity $ss^{-1}=I$ it follows easily that in the product
$$
l_{n}=I+\sum_{j=1}^{\infty}q_{j}(X_{1},X_{2})t^{j}:=
(I+\sum_{i=1}^{n}p_{i}(X_{1},X_{2})t^{i})\cdot
(I+\sum_{i=1}^{n}\widetilde p_{i}(X_{1},X_{2})t^{i})\ ,
$$
$q_{1}=\cdots =q_{n}=0$. 
Further observe that 
for an $h=I+\sum_{i=n+1}^{\infty}h_{i}(X_{1},X_{2})t^{i}\in H_{n}$,
there exist homogeneous polynomials 
$\widetilde h_{i}(X_{1},X_{2})$ of degree $i$, $n+1\leq i<\infty$, in
$X_{1}, X_{2}$ such that
$$
shs^{-1}=l_{n}+\sum_{i=n+1}^{\infty}\widetilde h_{i}(X_{1},X_{2})t^{i}\ .
$$
This shows that $shs^{-1}\in H_{n}$ and completes the proof of the lemma.
\ \ \ \ \ $\Box$

Now, the quotient group $Q_{n}:=H/H_{n}$ can be identified with the
set of the elements 
$$
s=I+\sum_{i=1}^{n}p_{i}(X_{1},X_{2})\hat t^{i}\ 
$$
where $p_{i}$ are homogeneous polynomials of degree $i$ in $X_{1},X_{2}$
with complex coefficients, and 
$$
\hat t^{i}\cdot\hat t^{j}=
\left\{
\begin{array}{ccc} 
\hat t^{i+j}&{\rm if}&i+j\leq n\\
0&{\rm if}&i+j>n
\end{array}
\right.
$$
Here the multiplication of the coefficients $p_{i}$ of such elements
is the same as in the algebra of non-commutative polynomials.

A straightforward computation shows that $Q_{n}$ is a nilpotent group.
Moreover, $Q_{n}$ is finite-dimensional and has a natural structure of a
complex algebraic group. In particular, $Q_{n}$ admits a 
faithful finite-dimensional unipotent representation \penalty-10000
$t_{n}:Q_{n}\rightarrow GL_{N}(\Co)$. (For the basic results about algebraic
groups see e.g. \penalty-10000 [OV, Chapter 3].)

Let $\pi_{n}:H\rightarrow Q_{n}$ be the quotient homomorphism. 
Observe that the family $\{t_{n}\circ\pi_{n}\}$ of finite-dimensional 
unipotent representations separates elements of $H$. This shows that
$H$ is residually torsion free nilpotent and completes the proof of
the theorem.\ \ \ \ \ $\Box$
{\sect{\hspace*{-1em}. Proofs of Results of Section 3.2.}
In this section we prove Propositions \ref{ser1}, \ref{cn}, and
Theorems \ref{comp1}, \ref{comp2}, \ref{comp3} and \ref{decomp}.\\
{\bf Proof of Proposition \ref{ser1}.} Let $G_{c}[[r]]$ be the
set of complex power series
$f(r)=r+\sum_{k=1}^{\infty}d_{k}(f)r^{k+1}$
each convergent in some open neighbourhood of $0\in\Co$. First, we will
show that $G_{c}[[r]]$ is a group with respect to the operation of 
composition. 

Suppose $f_{1},f_{2}\in G_{c}[[r]]$. Then by definition, there exists
an open disk $D_{i}$ centered at 0 such that $f_{i}$ is a holomorphic
function on $D_{i}$, $i=1,2$. In particular, there exists an open disk
$D_{1}'\subset D_{1}$ centered at 0 such that
$f_{1}(D_{1}')\subset D_{2}$. Therefore $f_{2}\circ f_{1}:D_{1}'\rightarrow
\Co$ is a holomorphic function such that $(f_{2}\circ f_{1})(0)=0$ and
$(f_{2}\circ f_{1})'(0)=1$. This implies that
$f_{2}\circ f_{1}$ has a Taylor series expansion in some open disk
centered at 0. Thus $f_{2}\circ f_{1}\in G_{c}[[r]]$. Next, for
$f\in G_{c}[[r]]$, by the inverse function theorem we obtain that
$f^{-1}$ exists, is holomorphic on an open disk containing 0, and
$f^{-1}(0)=0$, $(f^{-1})'(0)=1$. Hence, $f^{-1}\in G_{c}[[r]]$. 

Further, recall that the topology $\tau'$ on $G_{c}[[r]]$ is the weakest
topology in which all the functions $d_{k}:G_{c}[[r]]\rightarrow\Co$ are
continuous. It is obvious that the family $\{d_{k}\}_{k\geq 1}$ separates
points on $G_{c}[[r]]$.
Also, since the set of these functions is countable, $G_{c}[[r]]$
is metrizable (cf. the proof of Theorem \ref{group} (3)
for a similar argument). Thus in order to prove that $(G_{c}[[r]],\tau')$
is a topological group it suffices to check that
\begin{itemize}
\item[{\rm (1)}] if 
$\{f_{i}\}, \{g_{i}\}\subset G_{c}[[r]]$ converge to $f,g\in G_{c}[[r]]$,
respectively, then $\{f_{i}\circ g_{i}\}$ converges to $f\circ g$.
\item[{\rm (2)}]
if $\{f_{i}\}\subset G_{c}[[r]]$ converges to $f\in G_{c}[[r]]$,
then $\{f_{i}^{-1}\}$ converges to $f^{-1}$.
\end{itemize}
In turn, $\lim_{i\to\infty}h_{i}=h$ in $\tau'$ is equivalent to the fact
that $\lim_{i\to\infty}d_{k}(h_{i})=d_{k}(h)$ for any $k$. Next, note that
$d_{k}(f_{i}\circ g_{i})=P_{k}(d_{1}(f_{i}),\dots,d_{k}(f_{i}),
d_{1}(g_{i}),\dots, d_{k}(g_{i}))$ and $d_{k}(f_{i}^{-1})=
Q_{k}(d_{1}(f_{i}),\dots, d_{k}(f_{i}))$ for any $i$, where
$P_{k}\in\Co[z_{1},\dots,z_{2k}]$ and $Q_{k}\in\Co[z_{1},\dots,z_{k}]$
are complex holomorphic polynomials. This implies (1) and (2).

Further, clear that the subset of polynomials with rational coefficients 
is dense in $(G_{c}[[r]],\tau')$. Therefore $G_{c}[[r]]$ is separable.

Now, the contraction map can be defined by the formula
$$
F(t,f)(r):=\frac{f(tr)}{t},\ \ \ t\in [0,1],\ f\in G_{c}[[r]]\ .
$$
It is easy to check that $F:[0,1]\times G_{c}[[r]]\rightarrow G_{c}[[r]]$
is continuous, and $F(1,\cdot)=id$, $F(0,\cdot)=1$, where $1$ is the unit
of $G_{c}[[r]]$, i.e. the function $f(r)= r$. 

Now, the facts that
$G_{c}[[r]]$ is arcwise connected and locally simply and arcwise connected
can be obtained exactly in the same way as similar statements for
$G(X)$ (cf. the proof of Theorem \ref{group} (4)).

Finally, let us show that $G_{c}[[r]]$ is residually torsion free nilpotent.

To this end, let us introduce a subset $G_{n}\subset G_{c}[[r]]$ consisting
of series of the form
$$
f(r)=r+\sum_{k=n+1}^{\infty}d_{k}r^{k}\ .
$$
Then using the arguments similar to those of Lemma \ref{l92} one can
show that $G_{n}$ is a normal subgroup of $G_{c}[[r]]$. Further, the
quotient group $L_{n}:=G_{c}[[r]]/G_{n}$ can be identified with
the set of elements
$$
f(\widehat r)=\widehat r+\sum_{k=1}^{n}d_{k}\widehat r^{k}
$$
such that
$$
\widehat r^{i}\circ\widehat r^{j}=
\left\{
\begin{array}{ccc}
\widehat r^{ij}&{\rm if}& i\cdot j\leq n\\
0&{\rm if}& i\cdot j>n
\end{array}
\right.
$$
Then $L_{n}$ admits a faithful unipotent representation 
$q_{n}:L_{n}\rightarrow GL_{N}(\Co)$ because, as it is easy to see, $L_{n}$
has a structure of a finite-dimensional nilpotent algebraic group.
If $\pi_{n}:G_{c}[[r]]\rightarrow L_{n}$ is the quotient homomorphism, then,
by definition, the family $\{q_{n}\circ\pi_{n}\}$ separates elements of 
$G_{c}[[r]]$. Thus $G_{c}[[r]]$ is residually torsion free nilpotent.

The proof of the proposition is complete.\ \ \ \ \ $\Box$\\
{\bf Proof of Theorem \ref{comp1}.} We should check that
$P(a*b)=P(b)\circ P(a)$, for $a,b\in X$. As before, let $v(x;r;c)$,
$x\in [0,2\pi]$, be the Lipschitz solution of equation (\ref{e2})
corresponding to $c\in X$ with the initial value $v(0;r;c)=r$.
By the definition of $a*b$ we have 
$$
v(x;r;a*b)=\left\{
\begin{array}{ccc}
v(2x;r;a)&{\rm if}& 0<x\leq\pi\\
v(2x-2\pi;v(2\pi;r;a);b)&{\rm if}& \pi<x\leq 2\pi
\end{array}
\right.
$$
Thus by the definition of the first return map $P$ (see Section 2.1) we
obtain that
$$
P(a*b)(r):=v(2\pi;r;a*b)=v(2\pi;v(2\pi;r;a);b)=P(b)(P(a)(r))\ .
$$
This completes the proof of the theorem.\ \ \ \ \ $\Box$\\
{\bf Proof of Theorem \ref{comp2}.}
The map $\widehat P:G(X)\rightarrow G_{c}[[r]]$ is given by the formula
$$
\widehat P([a])=P(a)\ ,\ \ \ a\in X\ .
$$
It is correctly defined because $P(a)(r)\equiv r$ for any $a\in {\cal U}$.
Also, from Theorem \ref{comp1} it follows that $\widehat P$ is
a homomorphism of groups. Let us prove that $\widehat P$ is continuous.

Suppose $\{g_{n}\}\subset G(X)$ converges to $g\in G(X)$. We must show
that $\{\widehat P(g_{n})\}$ converges to $\widehat P(g)$. (This is
equivalent to the continuity of $\widehat P$ because $G(X)$ and
$G_{c}[[r]]$ are metric spaces.) Also, as it was mentioned before, the
latter is equivalent to the following statement
\begin{itemize}
\item[]
$\lim_{n\to\infty}d_{k}(\widehat P(g_{n}))=d_{k}(\widehat P(g))$
for any $k$, provided that $\lim_{n\to\infty}\widehat I(g_{n})=\widehat I(g)$
for any iterated integral $\widehat I$ on $G(X)$.
\end{itemize}
Now, according to Theorem \ref{center},
$d_{k}(\widehat P(g_{n}))=c_{k}(g_{n}')$, where $g_{n}'\in X$ is such that
$[g_{n}']=g_{n}$, and $c_{k}$ is an iterated integral on $X$. Thus the
above statement is true.
This shows that $\widehat P:G(X)\rightarrow G_{c}[[r]]$ is a homomorphism
of topological groups. 

Let us prove that $\widehat P$ is a surjection. 

Given $f(r)=r+\sum_{k=1}^{\infty}d_{k}r^{k+1}\in G_{c}[[r]]$ we define 
$a(f)=(a_{1},a_{2},\dots)$
by the identity of formal power series
\begin{equation}\label{101}
\sum_{i=1}^{\infty}a_{i}(x)t^{i+1}=
\frac{\sum_{k=1}^{\infty}
(d_{k}/2\pi)t^{k+1}}
{1+\sum_{k=1}^{\infty}(k+1)d_{k}(1-x/2\pi)t^{k}}\ \ \ 
{\rm for}\ \ \ x\in (0,2\pi]
\end{equation}
and further extended by periodicity. By our hypothesis on $f$ it follows
that the expression on the
right-hand side of (\ref{101}) is a bounded, continuous in $x\in (0,2\pi)$ 
function, holomorphic in $t$ for $t$ varying in a small open disk
centered at 0. This and the Cauchy inequalities for derivatives of a 
holomorphic function imply that $a(f)\in X$. Let $v(x;r;a(f))$ be the
Lipschitz solution of equation (\ref{e2}) corresponding to $a(f)\in X$
with the initial value $v(0;r;a(f))=r$. Then from equation (\ref{101})
it follows that for any sufficiently small $r$,
$$
v(x;r;a(f))+\sum_{k=1}^{\infty}d_{k}(1-x/2\pi)
[v(x;r;a(f))]^{k+1}=f(r)\ ,\ \ \ 0\leq x\leq 2\pi\ .
$$
Since $P(a(f))(r):=v(2\pi;r;a(f))$, the latter implies that
$$
P(a(f))(r)=f(r)\ .
$$
So $\widehat P$ is a surjection.

It is clear that $P(a)(r)\equiv r$ if and only if $a$ belongs to
the center set ${\cal C}\subset X$. This shows that 
$\widehat {\cal C}:=Ker(\widehat P)$ coincides with the image of ${\cal C}$
in $G(X)$. Then $\widehat {\cal C}$ is a normal subgroup of $G(X)$. It
is closed because $\widehat P$ is continuous. Since there are examples of
non-universal centers (see e.g. Remark \ref{che}), $\widehat {\cal C}$ is
non-trivial. Now, the contraction of $\widehat {\cal C}$ is defined as
follows.

For $a=(a_{1},a_{2},\dots)\in X$ we define
$a_{t}=(a_{1t},a_{2t},\dots)\in X$ by the formula
$$
a_{it}=t^{i}a_{i}\ ,\ \ \  1\leq i<\infty\ .
$$
Then 
\begin{equation}\label{102}
I_{i_{1},\dots, i_{k}}(a_{t})=t^{i_{1}+\dots+i_{k}}I_{i_{1},\dots, i_{k}}(a)
\end{equation}
for every basic iterated integral $I_{i_{1},\dots, i_{k}}$. This implies
that if $a\sim b$ then $a_{t}\sim b_{t}$. In particular, we can define
a map $F:[0,1]\times G(X)\rightarrow G(X)$ by the formula
$$
F(t,[a]):=[a_{t}]\ ,\ \ \ a\in X\ .
$$
Now from Theorem \ref{center} it follows that
$$
c_{n}(a_{t})=t^{n}c_{n}(a)\ ,\ \ \ a\in X,\ \ n=1,2,\dots\ .
$$
This shows that $F$ maps $\widehat {\cal C}$ to itself. Moreover,
$F(1,\cdot)=id$, $F(0,\cdot)=1$, and, according to (\ref{102}), $F$ is
continuous. Thus $F:[0,1]\times\widehat {\cal C}\rightarrow\widehat {\cal C}$
is the required contraction. Now the statements that $\widehat {\cal C}$
is arcwise connected, locally simply and arcwise connected can be obtained
in a way similar to that of the proof of Theorem \ref{group} (4).
We leave the details to the reader.

The proof of Theorem \ref{comp2} is complete.\ \ \ \ \ $\Box$\\
{\bf Proof of Proposition \ref{cn}.} By definition,
$\widehat c_{n}([a])=c_{n}(a)$, $a\in X$, where $c_{n}$,
$1\leq n<\infty$, are the Taylor coefficients of the first return map $P$ 
(see Theorem \ref{center}). In particular, we have
\begin{equation}\label{103}
\widehat P(s)(r)=r+\sum_{n=1}^{\infty}\widehat c_{n}(s)r^{n+1}\ ,\ \ \
s\in G(X)\ ,
\end{equation}
and the series converges absolutely for sufficiently small $r$.
Suppose $g\in \widehat {\cal C}$. Then, according to Theorem \ref{comp1},
for every $h\in G(X)$ we have $\widehat P(hg)=\widehat P(gh)=\widehat P(h)$.
This is equivalent to the identities
$$
r+\sum_{n=1}^{\infty}\widehat c_{n}(hg)r^{n+1}=
r+\sum_{n=1}^{\infty}\widehat c_{n}(gh)r^{n+1}=
r+\sum_{n=1}^{\infty}\widehat c_{n}(h)r^{n+1}\ .
$$
Thus 
$$
\widehat c_{n}(hg)=\widehat c_{n}(gh)=\widehat c_{n}(h)\ ,\ \ \ 
1\leq n<\infty\ .
$$
Therefore every $\widehat c_{n}$, $1\leq n<\infty$, is constant on fibres of 
the quotient homomorphism $\pi: G(X)\rightarrow Q(X):=G(X)/\widehat {\cal C}$,
and hence determines a function $\overline{c}_{n}:Q(X)\rightarrow\Co$
such that $\overline{c}_{n}\circ\pi=\widehat c_{n}$.

The proof of the proposition is complete.\ \ \ \ \ $\Box$\\
{\bf Proof of Theorem \ref{comp3}.} Let $\tau''$ be the weakest topology on
$Q(X)$ in which all functions $\overline{c}_{n}$ are continuous. According
to formula (\ref{103}) and Proposition \ref{cn}, the homomorphism
$\widehat P$ determines a homomorphism 
$\overline{P}:Q(X)\rightarrow G_{c}[[r]]$ defined by 
$\widehat P=\overline{P}\circ\pi$. In particular, we have
\begin{equation}\label{104}
\overline{P}(q)(r)=r+\sum_{n=1}^{\infty}\overline{c}_{n}(q)r^{n+1}\ ,\ \ \
q\in Q(X)\ ,
\end{equation}
and the series converges absolutely for sufficiently small $r$.
By definition, $\overline{P}:Q(X)\rightarrow G_{c}[[r]]$ is a group 
isomorphism. According to (\ref{104}), $\overline{P}$ is continuous,
because the functions $\overline{c}_{n}$ are continuous on $Q(X)$.
Let us prove now that $\overline{P}^{-1}$ is also continuous. Since
the set of functions $\{\overline{c}_{n}\}$ is countable, $Q(X)$ is 
metrizable (cf. the arguments used in the proof of Theorem \ref{group} (3)).
In particular, it suffices to prove the following statement
\begin{itemize}
\item[]
Suppose $\{f_{k}\}\subset G_{c}[[r]]$ is such that
$\lim_{k\to\infty}f_{k}=f\in G_{c}[[r]]$ in the topology $\tau'$ on
$G_{c}[[r]]$. Then the sequence 
$\{h_{k}:=\overline{P}^{-1}(f_{k})\}\subset Q(X)$ converges in $\tau''$
to $h:=\overline{P}^{-1}(f)$.
\end{itemize}

Let $f_{k}(r)=r+\sum_{n=1}^{\infty}d_{n}(f_{k})r^{n+1}$, $1\leq k<\infty$,
and $f(r)=r+\sum_{n=1}^{\infty}d_{n}(f)r^{n+1}$. Then identity (\ref{104}) 
implies that $\overline{c}_{n}(h_{k})=d_{n}(f_{k})$, and 
$\overline{c}_{n}(h)=d_{n}(h)$ for any $k,n$. Now by the hypothesis
we have $\lim_{k\to\infty}d_{n}(f_{k})=d_{n}(f)$ for any $n$. This and
the above identity imply that 
$\lim_{k\to\infty}\overline{c}_{n}(h_{k})=\overline{c}_{n}(h)$ for any $n$.
The latter is equivalent to $\lim_{k\to\infty}h_{k}=h$ in $\tau''$.
Thus $\overline{P}$ is an isomorphism of topological groups.

The proof of the theorem is complete.\ \ \ \ \ $\Box$\\
{\bf Proof of Theorem \ref{decomp}.} Let $a(f)\in X$,  $f\in G_{c}[[r]]$,
be defined by equation (\ref{101}). By 
$T:G_{c}[[r]]\rightarrow G(X)$ we denote the map $f\mapsto [a(f)]$.
In the proof of Theorem \ref{comp2} we established that
$(\widehat P\circ T)(f)=f$ for any $f\in G_{c}[[r]]$. Let us prove now that
$T$ is continuous. In fact, from (\ref{101}) it follows that
$$
a_{i}(x)=p_{i}(x,d_{1}(f),\dots,d_{i}(f))\ ,\ \ \ 1\leq i<\infty\ ,\ \ 
x\in (0,2\pi]\ ,
$$
where every $p_{i}\in\Co[z_{1},\dots,z_{i+1}]$ is a holomorphic polynomial.
In particular, by definition, for any basic iterated integral
$I_{i_{1},\dots,i_{k}}$ we have
$$
I_{i_{1},\dots,i_{k}}(a(f))=
\widetilde p_{i_{1},\dots,i_{k}}(d_{1}(f),\dots,d_{l}(f))\ ,
$$
where $l=\max\{i_{1},\dots,i_{k}\}$, and 
$\widetilde p_{i_{1},\dots,i_{k}}\in\Co[z_{1},\dots,z_{l}]$ is a holomorphic
polynomial. Now, this implies that 
\begin{itemize}
\item[]
if $\lim_{s\to\infty}f_{s}=f$
on $G_{c}[[r]]$, then 
$\lim_{s\to\infty}\widehat I_{i_{1},\dots,i_{k}}(T(f_{s}))=
\widehat I_{i_{1},\dots,i_{k}}(T(f))$ for every basic iterated integral
$\widehat I_{i_{1},\dots,i_{k}}$ on $G(X)$. 
\end{itemize}
This is equivalent to the
fact that $\lim_{s\to\infty}T(f_{s})=T(f)$ on $(G(X),\tau)$. Thus
$T$ is a continuous embedding.

Further, let us define 
$\widetilde T:G_{c}[[r]]\times\widehat {\cal C}\rightarrow G(X)$ by the 
formula
$$
\widetilde T(f,g)=T(f)\cdot g\ ,\ \ \ f\in G_{c}[[r]]\ ,\ 
g\in\widehat {\cal C}\ .
$$
Since $T$ and $\cdot$ are continuous maps, $\widetilde T$ is continuous.
Let us show that $\widetilde T$ is a bijection. Indeed, if we have
$T(f_{1})\cdot g_{1}=T(f_{2})\cdot g_{2}$ then
$$
f_{1}=\widehat P(T(f_{1})\cdot g_{1})=\widehat P(T(f_{2})\cdot g_{2})=f_{2}
\ .
$$
This implies that $(f_{1},g_{1})=(f_{2},g_{2})$. Thus $\widetilde T$ is
an injection. Now for any $h\in G(X)$ we obviously have
$$
h=\widetilde T(f,g)\ ,\ \ \ {\rm where}\ \ \ f=\widehat P(h)\ \ \
{\rm and}\ \ \ g=(T\circ\widehat P)(h^{-1})\cdot h\ .
$$
Thus $\widetilde T$ is a surjection. Also, the inverse to $\widetilde T$
is defined by the formula
$$
\widetilde T^{-1}(h):=(\widehat P(h), (T\circ\widehat P)(h^{-1})\cdot h)\ ,
\ \ h\in G(X)\ .
$$
It is clearly a continuous map because $\widehat P, T,\cdot$ and $^{-1}$
are continuous. 

Thus $\widetilde T$ is a homeomorphism. This completes the proof of the 
theorem.\ \ \ \ \ $\Box$
{\sect{\hspace*{-1em}. Proofs of Results of Section 3.3.}
In this section we prove Propositions \ref{subgroup1}, \ref{number1},
\ref{number2}, and Theorems \ref{subgroup2}, \ref{subgroup3},
\ref{number3}.\\
{\bf Proof of Proposition \ref{subgroup1}.} The statements 
\begin{itemize}
\item[{\rm (1)}]
the family of iterated integrals separates points on $G_{s}(X)$ and 
$G_{a}(X)$;
\item[{\rm (2)}]
$G_{s}(X)$ and $G_{a}(X)$ are separable topological groups in the 
relative topologies induced by $\tau$;
\item[{\rm (3)}]
$G_{s}(X)$ and $G_{a}(X)$ are residually torsion free nilpotent groups
\end{itemize}
follow directly from the fact that $G_{s}(X)$ and $G_{a}(X)$ are subgroups
of $G(X)$ and $G(X)$ satisfies similar properties.

Let us prove that $G_{s}(X)$ and $G_{a}(X)$ are contractible to a point,
arcwise connected, locally simply and arcwise connected.

We recall that the map $\Phi:[0,1]\times G(X)\rightarrow G(X)$
is defined by the formula
$$
\Phi(t,[a])=[ta]\ ,\ \ \ a\in X\ ,
$$
where $ta=(ta_{1},ta_{2},\dots)$ for $a=(a_{1},a_{2},\dots)$ 
(see the proof of Theorem \ref{group} (4)). Observe that if 
$a\in X_{a}$ then $ta\in X_{a}$ for any $t$. Since by the definition
$X_{a}\subset X_{s}$, from the above we obtain that
$\Phi$ maps $[0,1]\times G_{s}(X)$ into $G_{s}(X)$ and 
$[0,1]\times G_{a}(X)$ into $G_{a}(X)$. This and the properties of $\Phi$
described in the proof of Theorem \ref{group} (4) imply the required
statements.

The proof of the proposition is complete.\ \ \ \ \ $\Box$\\
{\bf Proof of Theorem \ref{subgroup2}.}
From formula (\ref{101}) it follows that for any $f\in G_{c}[[r]]$
the element $a(f)\in X_{a}\subset X_{s}$. Thus the map $T$,
$f\mapsto [a(f)]$, maps $G_{c}[[r]]$ into $G_{a}(X)\subset G_{s}(X)$.
Since $(\widehat P\circ T)(f)=f$, the latter implies that 
$$
\widehat P|_{G_{s}(X)}:G_{s}(X)\rightarrow G_{c}[[r]]\ \ \ {\rm and}
\ \ \ \widehat P|_{G_{a}(X)}:G_{a}(X)\rightarrow G_{c}[[r]]
$$
are surjective homomorphisms of topological groups. By definition their
kernels coincide with $\widehat {\cal C}_{s}$ and $\widehat {\cal C}_{a}$,
respectively. Also, the maps 
$\widetilde T_{s}:G_{c}[[r]]\times\widehat {\cal C}_{s}\rightarrow G_{s}(X)$,
$\widetilde T_{s}(f,g)=T(f)\cdot g$, and 
$\widetilde T_{a}:G_{c}[[r]]\times\widehat {\cal C}_{a}\rightarrow G_{a}(X)$,
$\widetilde T_{a}(f,g)=T(f)\cdot g$ are homeomorphisms (this can be
shown exactly in the same way as for $\widetilde T$ in the proof of
Theorem \ref{decomp}). From here it follows that 
$G_{s}(X)/\widehat {\cal C}_{s}$ and $G_{a}(X)/\widehat {\cal C}_{a}$ are
isomorphic to $G_{c}[[r]]$ which, in turn,  is isomorphic to $Q(X)$.

The proof of the theorem is complete.\ \ \ \ \ $\Box$\\
{\bf Proof of Theorem \ref{subgroup3}.} The proof repeats word-for-word
proofs of Proposition \ref{subgroup1} and Theorem \ref{subgroup2}.
We leave the details to the reader.\ \ \ \ \ $\Box$\\
{\bf Proof of Proposition \ref{number1}.}
From property (\ref{94}) for iterated integrals on $G(X)$ it follows that
it suffices to prove the following statement
\begin{itemize}
\item[]
The image in $G(X)$ of every  $a=(a_{1},a_{2},\dots)\in X_{\F'}(b)$ is defined
over $\F'$.
\end{itemize}
Recall that every $a_{k}$ is a linear function in 
$\hat b_{i_{1}},\dots, \hat b_{i_{s}}$ with coefficients from 
$\F'$. Here $\hat b=(\hat b_{1},\hat b_{2},\dots)\in X$ is such that
\begin{itemize}
\item[(a)]
its components run over the set of all possible monomials determined by 
elements of a sequence $b=\{b_{i}\}\subset L^{\infty}(S^{1})$ satisfying
$\sup_{i,S^{1}}|b_{i}|\leq 1$;
\item[(b)]
$[\hat b]\in G(X)$ is defined over the field $\F (\subseteq\F')$. 
\end{itemize}

Now, the definition of $a_{k}$ implies that every 
$I_{j_{1},\dots,j_{p}}(a)$ is a linear combination with
coefficients from $\F'$ of values of certain basic iterated integrals computed
at $\hat b$. In particular, $\widehat I_{j_{1},\dots,j_{p}}(a)\in\F'$.

This completes the proof of the proposition.\ \ \ \ \ $\Box$\\
{\bf Proof of Proposition \ref{number2}.}
In what follows by $[X_{\F'}(b)]^{-1}$ we denote the set of elements
$a\in X$ such that $a^{-1}\in X_{\F'}(b)$.

Let $\sigma:\F'\rightarrow\F'$ be a field automorphism over $\F$.
For any $a\in X_{\F'}(b)\subset X$ we define 
$\sigma(a), \sigma(a^{-1})\in X_{\F'}(b)$ as follows.

Suppose $a_{i}=\sum_{j=1}^{s(i)}c_{ij}\hat b_{ij}$ where $c_{ij}\in\F'$ and
$\hat b_{ij}$ are components of $\hat b\in X$, $1\leq j\leq s(i)$, $i\in\N$.
Then 
$$
\begin{array}{c}
\displaystyle
\sigma(a):=(\sigma(a_{1}),\sigma(a_{2}),\dots)\ , \ \ \ 
\sigma(a^{-1}):=(\sigma(a_{1}^{-1}),\sigma(a_{2}^{-1}),\dots)\ \ \
{\rm where}\\
\\
\displaystyle
\sigma(a_{i}):=\sum_{j=1}^{s(i)}\sigma(c_{ij})\hat b_{ij}\ ,
\ \ \
\sigma(a_{i}^{-1}):=\sum_{j=1}^{s(i)}\sigma(c_{ij})\hat b_{ij}^{-1}\ , \ \ \
1\leq j\leq s(i)\ ,\ i\in\N\ .
\end{array}
$$
In particular, we have $\sigma(a^{-1})=\sigma(a)^{-1}$.

Further, define the map
$\widetilde\sigma:G_{\F'}(b)\rightarrow G_{\F'}(b)$ by the formula
$$
\widetilde\sigma([a_{1}]\cdots [a_{s}])=
[\sigma(a_{1})]\cdots [\sigma(a_{s})]\ ,\ \ \ a_{1},\dots, a_{s}
\in X_{\F'}(b)\cup [X_{\F'}(b)]^{-1}\ .
$$
Let us show that $\widetilde\sigma$ is correctly defined.

First, suppose $a,a'\in X_{\F'}(b)\cup [X_{\F'}(b)]^{-1}$ are such that 
$a\sim a'$. Then for any basic
iterated integral we have $I(a)=I(a')$. Moreover, the definition for
the components $a_{i}$ and $a_{i}'$ of $a$ and $a'$ and the fact that
$[\hat b]$ is defined over the field $\F$ imply that there exist
$c_{1},\dots,c_{l}\in\F'$ (that coincide with coefficients 
$c_{ij}$, $c_{ij}'$ in the decompositions of certain $a_{i}$ and 
$a_{i}'$) such that $I(a)=p(a)$ and $I(a')=q(a')$ where 
$p,q\in\F[z_{1},\dots,z_{l}]$ are polynomials 
with coefficients from $\F$ and $p(a)$, $q(a')$ are their values at
the point $(c_{1},\dots, c_{l})\in (\F')^{l}$. Note that
for $\sigma(a)$ and $\sigma(a')$ we have $I(\sigma(a))=p(\sigma(a))$ and
$I(\sigma(a'))=q(\sigma(a'))$ with the same $p$ and $q$, where
$p(\sigma(a))$ and $q(\sigma(a'))$ are the values of $p$ and $q$ at
$(\sigma(c_{1}),\dots,\sigma(c_{l}))\in (\F')^{l}$. Since 
$\sigma$ is a field automorphism over $\F$, 
$$
I(\sigma(a))=p(\sigma(a))=\sigma(p(a))=\sigma(q(a'))=q(\sigma(a'))=
I(\sigma(a'))\ .
$$
This shows that $\sigma(a)\sim\sigma(a')$. Thus $[\sigma(a)]$, 
$a\in X_{\F'}(b)\cup [X_{\F'}(b)]^{-1}$, depends only on the equivalence 
class of $a$.

Suppose now $g_{1},\dots, g_{s},h_{1},\dots, h_{l}\in G_{\F'}(b)$ are 
such that  $g_{1}\cdots g_{s}=h_{1}\cdots h_{l}$, and $g_{i}=[g_{i}']$,
$1\leq i\leq s$, $h_{j}=[h_{j}']$, $1\leq j\leq l$, for some
$g_{1}',\dots, g_{s}', h_{1}',\dots,h_{l}'\in 
X_{\F'}(b)\cup [X_{\F'}(b)]^{-1}$. Then for a basic iterated integral
$\widehat I$ by formula (\ref{94}) we obtain that
$$
\widehat I(g_{1}\cdots g_{s})=\sum_{j=1}^{r}R_{j1}(g_{1})\cdots
\cdot R_{js}(g_{s})\ ,\ \ \
\widehat I(h_{1}\cdots h_{l})=\sum_{n=1}^{t}S_{n1}(h_{1})\cdots
\cdot S_{nl}(h_{l})\
$$
where all $R_{ji}$ and $S_{nk}$ are basic iterated integrals on $G(X)$.
Also, according to the above argument, 
$R_{ji}(g_{i})=p_{ji}(g_{i})$ and $S_{nk}=q_{nk}(h_{k})$
where there exists $m\in\N$ such that 
$p_{ji}, q_{nk}\in\F[z_{1},\dots,z_{m}]$ and
$p_{ji}(g_{i})$, $q_{nk}(h_{k})$ are the values of these polynomials at some
$(c_{1},\dots,c_{m})\in (\F')^{m}$ ($1\leq j\leq r$, $1\leq i\leq s$, 
$1\leq n\leq t$, $1\leq k\leq l$). In particular, we have (below we retain the
notation of the first part of the proof)
$$
\begin{array}{c}
\displaystyle
\widehat I(\widetilde\sigma(g_{1})\cdots\widetilde\sigma(g_{s}))=
\sum_{j=1}^{r}R_{j1}(\widetilde\sigma(g_{1}))\cdots 
R_{js}(\widetilde\sigma(g_{s}))=\sum_{j=1}^{r}\sigma(R_{j1}(g_{1}))\cdots
\sigma(R_{js}(g_{s}))=\\
\\
\displaystyle
\sigma(\widehat I(g_{1}\cdots g_{s}))=\sigma(\widehat I(h_{1}\cdots h_{l}))=
\sum_{n=1}^{l}\sigma(R_{n1}(h_{1}))\cdots\sigma(R_{nl}(h_{l}))=\\
\\
\displaystyle
\sum_{n=1}^{l}R_{n1}(\widetilde\sigma(h_{1}))\cdots 
R_{nl}(\widetilde\sigma(h_{l}))=
\widehat I(\widetilde\sigma(h_{1})\cdots\widetilde\sigma(h_{l}))\ .
\end{array}
$$
This shows that $\widetilde\sigma(g_{1})\cdots\widetilde\sigma(g_{s})=
\widetilde\sigma(h_{1})\cdots\widetilde\sigma(h_{l})$. Thus 
$\widetilde\sigma$ is correctly defined. Clearly, $\widetilde\sigma(1)=1$.
Therefore the definition of $\widetilde\sigma$ shows that it is an 
automorphism of $G_{\F'}(X)$. (Here we do not assert that 
$\widetilde\sigma$ is continuous in the topology on $G_{\F'}(X)$.)

Let us proof now that $\widetilde\sigma(\widehat {\cal C}_{\F'}(b))=
\widehat {\cal C}_{\F'}(b)$.

Suppose $g\in \widehat {\cal C}_{\F'}(b)$. Then $\widehat c_{n}(g)=0$,
$1\leq n<\infty$, where $\widehat c_{n}$ are defined in Section 3.2.
According to Theorem \ref{center} every $\widehat c_{n}$ is the sum of basic
iterated integrals with integer coefficients.  Also, according to
the above computation,
$\sigma(\widehat I(h))=\widehat I(\widetilde\sigma(h))$, $h\in G_{\F'}(b)$,
for any basic iterated integral $\widehat I$. Therefore for any $n$ we have
$$
\widehat c_{n}(\widetilde\sigma(g))=\sigma(\widehat c_{n}(g))=0\ .
$$
This equivalent to the fact that 
$\widetilde\sigma(g)\in \widehat {\cal C}_{\F'}(b)$.

The proof of the proposition is complete.\ \ \ \ \ $\Box$\\
{\bf Proof of Theorem \ref{number3}.}
First, let us prove that $[P]:[0,1]\rightarrow G(X)$ is a continuous map.
To this end, we must check that $\widehat I([P](t))$ is a continuous function
on $[0,1]$ for every basic iterated integral $\widehat I$. Now,
according to formulae (\ref{94}), (\ref{95}) it suffices to prove the 
statement for $P(t)=c(t)$ where 
$c(t)=(c_{1}(t),c_{2}(t),\dots)$ with $c_{k}(t)$ polynomials
in $t$ whose coefficients are elements from ${\cal A}_{\F'}(b)$.
But then a simple computation shows that $\widehat I(c(t))\in\Co[t]$.
This implies the continuity of $[P]$.

Now let us prove the theorem. According to Proposition \ref{number1} and
formulae (\ref{94}), (\ref{95}), 
for every function $\widehat c_{n}$ (defined in Section 3.2) the function
$\widehat c_{n}([P](t))\in\F'[t]$. Now, the hypothesis of the
theorem implies that $\widehat c_{n}([P](\xi))=0$ for every $n$, and there
exists $k\in\N$ such that $\widehat c_{k}([P](0))\neq 0$. In particular,
this shows that $\widehat c_{k}([P])(t)$ is a non-trivial polynomial.
Since the coefficients of this polynomial belong to $\F'$ and
$\widehat c_{k}([P](\xi))=0$, the number $\xi$ is algebraic over $\F'$.

This completes the proof of the theorem.\ \ \ \ \ $\Box$


\end{document}